\documentclass[righttag]{article}
\usepackage{amsmath}
\usepackage{graphicx}
\usepackage{amsthm}
\usepackage{amsfonts}
\usepackage{amssymb}
\newtheorem{theorem}{Theorem}[section]

\newtheorem{corollary}[theorem]{Corollary}

\newtheorem{lemma}[theorem]{Lemma}

\newtheorem{proposition}[theorem]{Proposition}
\newtheorem{remark}[theorem]{Remark}

\numberwithin{equation}{section}

\begin{document}

\title{EXISTENCE THEOREMS OF FOLD-MAPS}
\author{YOSHIFUMI ANDO \thanks{This research was partially supported by Grant-in-Aid
for Scientific Research (No. 11640081), Ministry of Education, Science and
Culture, Japan.}}
\date{}
\maketitle
\begin{abstract}
Let $N$ and $P$ be smooth manifolds of dimensions $n$ and $p$ ($n\geq p\geq2$)
respectively. A smooth map having only fold singularities is called a
fold-map. We will study conditions for a continuous map $f:N\rightarrow P$ to
be homotopic to a fold-map from the viewpoint of the homotopy principle. By
certain homotopy principles for fold-maps, we prove that if there exists a
fiberwise epimorphism $TN\oplus\theta_{N}\rightarrow TP$ covering $f$, then
there exists a fold-map homotopic to $f$, where $\theta_{N}$ is the trivial
line bundle. We also give an additional condition for finding a fold-map which
folds only on a finite number of spheres of dimension $p-1$.
\end{abstract}

\section*{Introduction}

Let $N$ and $P$ be smooth ($C^{\infty}$) manifolds of dimensions $n$ and $p$
respectively with $n\geq p\geq2$. We say that a smooth map germ of $(N,x)$
into $(P,y)$ has a singularity of fold type at $x$ if it is written as
$(x_{1},\ldots,x_{n})\mapsto(x_{1},\ldots,x_{p-1},\pm x_{p}^{2}\pm\cdots\pm
x_{n}^{2})$ under suitable local coordinate systems on neighborhoods of $x\in
N$ and $y\in P$. A smooth map $f:N\rightarrow P$ is called a \textit{fold-map}
if it has only fold singularities. We will study conditions for a given
continuous map to be homotopic to a fold-map from the viewpoint of the
homotopy principle. Main results of this paper are the following. Let $TN$,
$TP$ and $\theta_{N}$ be the tangent bundles of $N$, $P$ and the trivial
bundle $N\times\mathbf{R}$ respectively.$\ $

\begin{theorem}
Let $n\geq p\geq2$. Let $f:N\rightarrow P$ be a continuous map. Assume that
there exists a fiberwise epimorphism $h:TN\oplus\theta_{N}\rightarrow TP$
covering $f$. Then there exists a fold-map $g:N\rightarrow P$ homotopic to $f$.
\end{theorem}

If $n-p+1$ is odd, then we can easily prove that the converse also holds
(Lemma 3.1).

Recently Saeki[Sa2] has given a complete characterization of those
$4$-dimensional closed orientable manifolds which admit fold-maps into
$\mathbf{R}^{3}$. The converse of Theorem 0.1 does not hold for the dimension
pair $(n,p)=(4,3)$.

\begin{corollary}
Let $n\geq p\geq2$. Let $f:N\rightarrow P$ be a continuous map. If the
structure group of $TN\oplus\theta_{N}$ is reduced to $O(n-p)\times E_{p+1}$
$($this refers to $E_{p+1}$ for the case $n=p)$ and if $TP$ is stably trivial,
then there exists a fold-map $g:N\rightarrow P$ homotopic to $f$.
\end{corollary}

Here, $O(k)$ is the orthogonal group of degree $k$ and $E_{k}$ is the unit
matrix of degree $k$. We have another existence theorem of fold-maps.

\begin{theorem}
Let $n\geq p\geq2$. Let $f:N\rightarrow P$ be a continuous map. Assume that
there exists a vector bundle $\xi$\ of dimension $n-p$ over $N$ such that $TN$
is stably equivalent to the vector bundle $f^{\ast}(TP)\oplus\xi$. Then $f$ is
homotopic to a fold-map $g:N\rightarrow P$ which folds exactly on the
boundaries of a finite number of \ disjointly embedded disks of dimension $p$
within an embedded disk of dimension $p$ in $N$.
\end{theorem}

\begin{corollary}
Let $n\geq p\geq2$. Let $f:N\rightarrow P$ be a continuous map. If $TN$ and
$TP$ are stably trivial, then there exists a fold-map $g:N\rightarrow P$
homotopic to $f$ which folds exactly on the boundaries of a finite number of
\ disjointly embedded disks of dimension $p$ within an embedded disk of
dimension $p$ in $N$.
\end{corollary}

Corollary 0.4 should be compared with [E2, 5.1D and 5.4]. For example, any
continuous map $f:S^{n}\rightarrow S^{p}$ is homotopic to such a fold-map
(Corollary 3.4).

Theorem 0.1 is a simple consequence of Theorem 0.5 below, while we need the
Phillips Submersion Theorem ([P]) and a further tool in the proof of Theorem
0.3 to apply Theorem 2.4.

We explain the terminology \textit{homotopy principle} used in [G2]. In the
2-jet space $J^{2}(n,p)$, let $\Sigma^{n-p+1,0}(n,p)$ be the subspace of all
jets of germs with fold singularities at the origin and let $\Omega
^{n-p+1,0}(n,p)$ be the union of all jets of regular germs and $\Sigma
^{n-p+1,0}(n,p)$. In the 2-jet space $J^{2}(N,P)$ with projection $\pi_{N}%
^{2}\times\pi_{P}^{2}:J^{2}(N,P)\rightarrow N\times P$, let $\Sigma
^{n-p+1,0}(N,P)$ and $\Omega^{n-p+1,0}(N,P)$ be its subbundles associated with
$\Sigma^{n-p+1,0}(n,p)$ and $\Omega^{n-p+1,0}(n,p)$ respectively. A smooth map
$f:N\rightarrow P$ is a fold-map if and only if the image of $j^{2}f$ is
contained in $\Omega^{n-p+1,0}(N,P)$ and transverse to $\Sigma^{n-p+1,0}%
(N,P)$. Let $C_{\Omega}^{\infty}(N,P)$ denote the space consisting of all
fold-maps equipped with the $C^{\infty}$-topology. Let $\Gamma(N,P)$ denote
the space consisting of all continuous sections of the fiber bundle $\pi
_{N}^{2}|\Omega^{n-p+1,0}(N,P):\Omega^{n-p+1,0}(N,P)\rightarrow N$ equipped
with the compact-open topology. Then there exists a continuous map
\[
j_{\Omega}:C_{\Omega}^{\infty}(N,P)\rightarrow\Gamma(N,P)
\]
defined by $j_{\Omega}(f)=j^{2}f$. It follows from the well-known theorem due
to Gromov[G1] that if $N$ is a connected open manifold, then $j_{\Omega}$ is a
weak homotopy equivalence. This property is called the homotopy principle.

The existence problem of fold-maps has been first dealt with in [T] and [L1]
in dimensions $n\geq p=2$ from a different viewpoint. In [E1] and [E2]
$\grave{\mathrm{E}}\mathrm{lia}\check{\mathrm{s}}\mathrm{berg}$ has proved a
certain homotopy principle for fold-maps in the existence level for closed
smooth manifolds.

We will prove the following theorem in the existence level (see also Theorem
2.4) in $\S2$, where two theorems [G1, 4.1.1 Theorem] and [E2, 4.7 Theorem]
will play important roles.

\begin{theorem}
Let $n\geq p\geq2$. Let $N$ and $P$ be connected manifolds of dimensions $n$
and p respectively with $\partial N=\emptyset$. Let $C$ be a closed subset of
$N$. Let $s$ be a section of $\Gamma(N,P)$ such that there exists a fold-map
$g$ defined on a neighborhood of $C$ into $P$, where $j^{2}g=s$. Then there
exists a fold-map $f:N\rightarrow P$ such that $j^{2}f$ is homotopic to s
relative to $C$ by a homotopy $h_{\lambda}$ in $\Gamma(N,P)$ with $h_{0}=s$
and $h_{1}=j^{2}f$.
\end{theorem}

Most of the results in the case $n=p$ of this paper except for Theorems 0.3,
2.4 and Proposition 5.7 have already been proved in [An3].

We refer to [B-R], [S-S] and its references in low dimensions $(3,2)$ and
$(4,3)$ in another line of investigation concerning the existence problem of
fold-maps of special generic type, which are closely related to the
differentiable structures of manifolds.

In $\S1$ we explain well-known results concerning fold singularities. In $\S2
$ we state Theorem 2.4, which is another type of homotopy principle for
fold-maps, and Proposition 2.5 without proofs, and prove Theorem 0.5 by using
them. In $\S3$ we prove Theorems 0.1 and 0.3 by using Theorems 0.5 and 2.4. In
$\S4$ we prove Theorem 2.4. In $\S5$ we prove Proposition 2.5.

The author would like to thank the referee for his kind and helpful comments,
which improved the paper.

\section{Preliminaries}

Throughout the paper all manifolds are smooth of class $C^{\infty}$. Maps are
basically continuous, but may be smooth (of class $C^{\infty}$) if so stated.
We always work in dimensions $n\geq p\geq1.$

Given a fiber bundle $\pi:E\rightarrow X$ and a subset $C$ in $X,$ we denote
$\pi^{-1}(C)$ by $E|_{C}.$ Let $\pi^{\prime}:F\rightarrow Y$ be another fiber
bundle. A map $\tilde{b}:E\rightarrow F$ is called a fiber map over a map
$b:X\rightarrow Y$ if $\pi^{\prime}\circ\tilde{b}=b\circ\pi$ holds. The
restriction $\tilde{b}|(E|_{C}):E|_{C}\rightarrow F$ (or $F|_{b(C)}$) is
denoted by $\tilde{b}|_{C}$. In particular, for a point $x\in X,$ $E|_{x}$ and
$\tilde{b}|_{x}$ are simply denoted by $E_{x}$ and $\tilde{b}_{x}%
:E_{x}\rightarrow F_{b(x)}$ respectively. When $E$ and $F$ are vector bundles
over $X=Y$, Hom$(E,F)$ denotes the vector bundle over $X$ with fiber
Hom$(E_{x},F_{x})$, $x\in X$, which consists of all homomorphisms
$E_{x}\rightarrow F_{x}$. A fiberwise homomorphism, epimorphism and
monomorphism $E\rightarrow F$ are simply called homomorphism, epimorphism and
monomorphism respectively. The trivial bundle $X\times\mathbf{R}^{k}$ is
denoted by $\theta_{X}^{k}$. In particular, $\theta_{X}^{1}$ is often written
as $\theta_{X}$.

We review well-known results about fold singularities (see [B], [L2]). Let
$J^{k}(N,P)$ denote the $k$-jet space of manifolds $N$ and $P$. Let $\pi
_{N}^{k}$ and $\pi_{P}^{k}$ be the projections mapping a jet to its source and
target respectively. The map $\pi_{N}^{k}\times\pi_{P}^{k}:J^{k}%
(N,P)\rightarrow N\times P$ induces a structure of a fiber bundle with
structure group $L^{k}(p)\times L^{k}(n)$, where $L^{k}(m)$ denotes the group
of all $k$-jets of local diffeomorphisms of $(\mathbf{R}^{m},0)$. The fiber
$(\pi_{N}^{k}\times\pi_{P}^{k})^{-1}(x,y)$ is denoted by $J_{x,y}^{k}(N,P)$.

Let $\pi_{1}^{2}:J^{2}(N,P)\rightarrow J^{1}(N,P)$ be the canonical forgetting
map. Let $\Sigma^{i}(N,P)$ denote the submanifold of $J^{1}(N,P)$\ consisting
of all $1$-jets $z=j_{x}^{1}f$ such that the kernel of $d_{x}f$ is of
dimension $i$. Let $\Omega^{n-p+1}(N,P)$ denote the union of $\Sigma
^{n-p}(N,P)$ and $\Sigma^{n-p+1}(N,P)$ in $J^{1}(N,P)$. We denote $(\pi
_{1}^{2})^{-1}(\Sigma^{i}(N,P))$ by the same symbol $\Sigma^{i}(N,P)$ if there
is no confusion. For a $2$-jet $z=j_{x}^{2}f$ of $\Sigma^{i}(N,P)$, there has
been defined the second intrinsic derivative $d_{x}^{2}f:T_{x}N\rightarrow
\mbox{{\rm Hom}}(\mbox{{\rm Ker}}(d_{x}f),\mbox{{\rm Cok}}(d_{x}f))$. Let
$\Sigma^{i,j}(N,P)$ denote the subbundle of $J^{2}(N,P)$ consisting of all
jets $z=j_{x}^{2}f$ such that $\dim(\mbox
{{\rm Ker}}(d_{x}f))=i$ and $\dim(\mbox
{{\rm Ker}}(d_{x}^{2}f|\mbox{{\rm Ker}}(d_{x}f)))=j$. A jet of $\Sigma
^{n-p+1,0}(N,P)$ will be called a fold jet. Let $\Omega^{n-p+1,0}(N,P)$ denote
the union of $\Sigma^{n-p}(N,P)$ and $\Sigma^{n-p+1,0}(N,P)$ in $J^{2}(N,P)$.
Then $\pi_{N}^{2}\times\pi_{P}^{2}|\Omega^{n-p+1,0}(N,P)$ induces a structure
of an open subbundle of $\pi_{N}^{2}\times\pi_{P}^{2}$. Let $J_{x,y}^{k}%
(N,P)$, $\Sigma_{x,y}^{i}(N,P)$, $\Omega_{x,y}^{n-p+1}(N,P)$, $\Sigma
_{x,y}^{n-p+1,0}(N,P)$, and $\Omega_{x,y}^{n-p+1,0}(N,P)$ denote the
respective fibers over $(x,y)$. In particular, we set $J^{k}(n,p)=J_{0,0}%
^{k}(\mathbf{R}^{n},\mathbf{R}^{p})$, $\Sigma^{i}(n,p)=\Sigma_{0,0}%
^{i}(\mathbf{R}^{n},\mathbf{R}^{p})$, $\Omega^{n-p+1}(n,p)=\Omega
_{0,0}^{n-p+1}(\mathbf{R}^{n},\mathbf{R}^{p})$, $\Sigma^{n-p+1,0}%
(n,p)=\Sigma_{0,0}^{n-p+1,0}(\mathbf{R}^{n},\mathbf{R}^{p})$, and
$\Omega^{n-p+1,0}(n,p)=\Omega_{0,0}^{n-p+1,0}(\mathbf{R}^{n},\mathbf{R}^{p})$.

Let $\pi_{N}$ and $\pi_{P}$ be the projections of $N\times P$ onto $N$ and $P$
respectively.\ We set
\begin{equation}
J^{2}(TN,TP)=\text{\textrm{Hom}}(\pi_{N}^{\ast}(TN),\pi_{P}^{\ast}%
(TP))\oplus\text{\textrm{Hom}}(S^{2}(\pi_{N}^{\ast}(TN)),\pi_{P}^{\ast}(TP))
\end{equation}
over $N\times P$, where $S^{2}(\pi_{N}^{\ast}(TN))$ is the 2-fold symmetric
product of $\pi_{N}^{\ast}(TN)$. If we provide $N$ and $P$ with Riemannian
metrics, then the Levi-Civita connections induce the exponential maps
$\exp_{N,x}:T_{x}N\rightarrow N$ and $\exp_{P,y}:T_{y}P\rightarrow P$. In
dealing with the exponential maps we always consider the so called convex
neighborhoods ([K-N]). We define the bundle map
\begin{equation}
J^{2}(N,P)\mathbf{\rightarrow}J^{2}(TN,TP)\text{ \ \ \ over }N\times P
\end{equation}
by sending $z=j_{x}^{2}f\in J_{x,y}^{2}(N,P)$ to the $2$-jet of $(\exp
_{P,y})^{-1}\circ f\circ\exp_{N,x}$ at $\mathbf{0}\in T_{x}N$, which is
regarded as an element of $J^{2}(T_{x}N,T_{y}P)(=J_{x,y}^{2}(TN,TP))$. The
structure group of $J^{2}(TN,TP)$ is reduced to $O(p)\times O(n)$.

In this paper we often express an element of $J_{x,y}^{2}(N,P)$ as
$(\alpha,\beta)$ for $\alpha\in\mbox{{\rm
Hom}}(T_{x}N,T_{y}P)$ and $\beta\in\mbox
{{\rm Hom}}(S^{2}(T_{x}N),T_{y}P)$. For a subspace $V$ in $T_{y}P,$ let
$pr(V)$ be the orthogonal projection of $T_{y}P$ onto $V$. For an element
$(\alpha,\beta)\in\Sigma_{x,y}^{n-p+1}(N,P)$, let $\beta_{\alpha}$ denote the
homomorphism defined by
\begin{equation}
\beta_{\alpha}=pr(\mathrm{Im}(\alpha)^{\perp})\circ(\beta|S^{2}\mathrm{Ker}%
(\alpha)),
\end{equation}
where the symbol $\perp$ refers to the orthogonal complement. Under the
identification (1.2), $\alpha\in J_{x,y}^{1}(N,P)$ lies in $\Sigma_{x,y}%
^{i}(N,P)$ if and only if dim \textrm{Ker}$(\alpha)=i,$ and $(\alpha,\beta
)\in\Sigma_{x,y}^{n-p+1}(N,P)$ lies in $\Sigma_{x,y}^{n-p+1,0}(N,P)$ if and
only if $\beta_{\alpha}$ is a non-singular quadratic form. Let $\iota$ be an
integer such that $0\leq\iota\leq\lbrack(n-p+1)/2]$ ([$a$] refers to the
greatest integer not exceeding $a$). Let $\Sigma^{n-p+1,0}(n,p)^{\iota}$
denote the subspace which consists of all elements $(\alpha,\beta)\in
\Sigma^{n-p+1,0}(n,p)$ such that the index (the number of negative eigen
values) of $\beta_{\alpha}$ is equal to $\iota$ or $n-p+1-\iota$ depending on
the choice of the orientation of $\mathrm{Im}(\alpha)^{\perp}$. Let
$\Sigma^{n-p+1,0}(N,P)^{\iota}$ denote the subspace of $\Sigma^{n-p+1,0}(N,P)$
associated to $\Sigma^{n-p+1,0}(n,p)^{\iota}$.

By (1.1) and (1.2), $J^{2}(n,p)$ is canonically identified with
\begin{equation}
\mbox{{\rm Hom}}(\mathbf{R}^{n},\mathbf{R}^{p})\oplus\mbox{{\rm
Hom}}(S^{2}\mathbf{R}^{n},\mathbf{R}^{p})
\end{equation}
under the canonical bases of $\mathbf{R}^{n}$ and $\mathbf{R}^{p}$. For a jet
$z=j_{x}^{2}f\in J^{2}(\mathbf{R}^{n},\mathbf{R}^{p})$, we define $\pi_{J}$ by
$\pi_{J}(z)=j_{0}^{2}(\mathit{l}(-f(x))\circ f\circ\mathit{l}(x))$, where
$\mathit{l}(a)$ denotes the parallel translation defined by $\mathit{l}%
(a)(x)=x+a$. Let $\pi_{\Omega}:\Omega^{n-p+1,0}(\mathbf{R}^{n},\mathbf{R}%
^{p})\rightarrow\Omega^{n-p+1,0}(n,p)$ be the restriction of $\pi_{J}$. We
obtain the canonical diffeomorphisms
\begin{align}
\pi_{\mathbf{R}^{n}}^{2}\times\pi_{\mathbf{R}^{p}}^{2}\times\pi_{J}  &
:J^{2}(\mathbf{R}^{n},\mathbf{R}^{p})\rightarrow\mathbf{R}^{n}\times
\mathbf{R}^{p}\times(\mbox{{\rm Hom}}(\mathbf{R}^{n},\mathbf{R}^{p}%
)\oplus\mbox{{\rm
Hom}}(S^{2}\mathbf{R}^{n},\mathbf{R}^{p})),\\
\pi_{\mathbf{R}^{n}}^{2}\times\pi_{\mathbf{R}^{p}}^{2}\times\pi_{\Omega}  &
:\Omega^{n-p+1,0}(\mathbf{R}^{n},\mathbf{R}^{p})\rightarrow\mathbf{R}%
^{n}\times\mathbf{R}^{p}\times\Omega^{n-p+1,0}(n,p).
\end{align}

Furthermore, we always identify $\mbox{{\rm Hom}}(\mathbf{R}^{n}%
,\mathbf{R}^{p})$ with the space $M_{p\times n}$ of all $p\times n$ matrices
and identify $\mbox{{\rm Hom}}(S^{2}\mathbf{R}^{n},\mathbf{R}^{p})$ with the
space of all $p$-tuples of $n\times n$ symmetric matrices throughout the paper.

Next we review the properties of the submanifolds $\Sigma^{n-p+1}(N,P)$ and
$\Sigma^{n-p+1,0}(N,P)$ along the line of [B, $\S7]$. Let $\mathbf{D}^{\prime
}$ denote the induced bundle $(\pi_{N}^{2})^{\ast}(TN)$ over $J^{2}(N,P)$.
Recall the homomorphism
\[
\mathbf{d}^{1}:\mathbf{D}^{\prime}\mathbf{\rightarrow}(\pi_{P}^{2})^{\ast
}(TP)\qquad\mbox{over $J^2(N,P)$},
\]
which maps an element $\mathbf{v=(}z,\mathbf{v}^{\prime})\in\mathbf{D}%
_{z}^{\prime}$ with $z=j_{x}^{2}f$ to $(z,d_{x}f(\mathbf{v}^{\prime}))$. Here
$\mathbf{d}^{1}$ is identified with a section of $\mbox{{\rm Hom}}%
(\mathbf{D}^{\prime},(\pi_{P}^{2})^{\ast}(TP))$ over $J^{2}(N,P)$. Let
$\mathbf{K}$ and $\mathbf{Q}$ be the kernel bundle and the cokernel bundle of
$\mathbf{d}^{1}$ over $\Sigma^{n-p+1}(N,P)$ with $\dim\mathbf{K=}n-p+1$ and
$\dim\mathbf{Q}=1$ respectively. Then we have the second intrinsic derivative
$\mathbf{d}^{2}$ $:\mathbf{D}^{\prime}\rightarrow\mathrm{Hom}(\mathbf{K}%
,\mathbf{Q})$ over $\Sigma^{n-p+1}(N,P)$. Indeed, for $z=j_{x}^{2}f\in
\Sigma^{n-p+1}(N,P)$, the induced homomorphism $((j^{2}f)^{\ast}\mathbf{d}%
^{2})_{x}:((j^{2}f)^{\ast}\mathbf{D}^{\prime})_{x}\rightarrow\mathrm{Hom}%
(((j^{2}f)^{\ast}\mathbf{K})_{x},((j^{2}f)^{\ast}\mathbf{Q})_{x})$ is nothing
but the homomorphism induced from $d_{x}^{2}f:T_{x}N\rightarrow\mathrm{Hom}%
(\mathrm{Ker}(d_{x}f),\mathrm{Cok}(d_{x}f))$ by $(\pi_{N}^{2})^{\ast}$ and
$(\pi_{P}^{2})^{\ast}$. As is explained in [B, p. 412], the second intrinsic
derivative $\mathbf{d}^{2}|\mathbf{K}$ is extended to the epimorphism%
\[
\mathbf{d}^{2}:T(J^{2}(N,P))|_{\Sigma^{n-p+1}(N,P)}\rightarrow\mathrm{Hom}%
(\mathbf{K},\mathbf{Q})\qquad\mbox{over
$\Sigma^{n-p+1}(N,P)$},
\]
where $\mathbf{K}$ is regarded as a subbundle of $T(J^{2}(N,P))|_{\Sigma
^{n-p+1}(N,P)}$. Furthermore, on the subset of $\Sigma^{n-p+1}(N,P)$
consisting of all jets $z=j_{x}^{2}f$ such that $d(j^{2}f)$ is transverse to
$\Sigma^{n-p+1}(N,P)$ at $x$, $\mathbf{D}^{\prime}$ is regarded as a subbundle
of $T(J^{2}(N,P))|_{\Sigma^{n-p+1}(N,P)}$ by [B, Theorem 7.15] (see also
another interpretation of [An3, \S1]). It has been proved in [B, Lemma 7.13]
that there exists an exact sequence%
\[
0\longrightarrow T(\Sigma^{n-p+1}(N,P))\overset{\subset}{\longrightarrow
}T(J^{2}(N,P))|_{\Sigma^{n-p+1}(N,P)}\overset{\mathbf{d}^{2}}{\longrightarrow
}\mathrm{Hom}(\mathbf{K},\mathbf{Q})\longrightarrow0.
\]
Under these notations, a 2-jet $z\in\Sigma^{n-p+1}(N,P)$ lies in
$\Sigma^{n-p+1,0}(N,P)$ if and only if $\mathbf{d}^{2}|\mathbf{K}_{z}$ is an
isomorphism. This implies that $T(\Sigma^{n-p+1}(N,P))_{z}\cap\mathbf{K}%
_{z}=\{0\}$ for any jet $z\in\Sigma^{n-p+1,0}(N,P)$. Hence $\mathbf{K}%
|_{\Sigma^{n-p+1,0}(N,P)}$ and $\mbox
{{\rm Hom}}(\mathbf{K},\mathbf{Q})|_{\Sigma^{n-p+1,0}(N,P)}$ are isomorphic to
the normal bundle of $\Sigma^{n-p+1,0}(N,P)$ in $J^{2}(N,P)$.

Let $C_{\Omega}^{\infty}(N,P)$ and $\Gamma(N,P)$ denote the spaces defined in
\textrm{Int}roduction with the continuous map $j_{\Omega}:C_{\Omega}^{\infty
}(N,P)\rightarrow\Gamma(N,P)$. Let $\Gamma^{tr}(N,P) $ denote the subspace of
$\Gamma(N,P)$ consisting of all sections $s$ such that $s$ is smooth on some
neighborhood of $s^{-1}(\Sigma^{n-p+1,0}(N,P))$ and that $s$ is transverse to
$\Sigma^{n-p+1,0}(N,P)$. Throughout the paper, $s^{-1}(\Sigma^{n-p+1,0}(N,P))$
is denoted by $S(s)$. Let $K(s)$ and $Q(s)$ denote $(s|S(s))^{\ast}%
(\mathbf{K})$ and $(s|S(s))^{\ast}(\mathbf{Q})$ respectively. Let
$d^{1}s:TN\rightarrow(\pi_{P}^{2}\circ s)^{\ast}(TP)$ and $d^{2}%
s:TN|_{S(s)}\rightarrow\mathrm{Hom}(K(s),Q(s))$ over $S(s)$ denote the
homomorphisms induced from $\mathbf{d}^{1}$ and $\mathbf{d}^{2}|\mathbf{D}%
^{\prime}$ by $s$ respectively. The map $d^{2}s$ induces a symmetric
homomorphism $q(s):S^{2}K(s)\rightarrow Q(s)$, which we call the quadratic
form associated with $d^{2}s$.

Then $ds|(TN|_{S(s)}):TN|_{S(s)}\rightarrow T(J^{2}(N,P))|_{\Sigma
^{n-p+1,0}(N,P)}$ canonically induces the homomorphism
\begin{equation}
d(s)_{\thicksim}:TN|_{S(s)}\rightarrow(s|S(s))^{\ast}(T(J^{2}(N,P))|_{\Sigma
^{n-p+1,0}(N,P)}),
\end{equation}
and $\mathbf{d}^{2}:T(J^{2}(N,P))|_{\Sigma^{n-p+1,0}(N,P)}\rightarrow
\mathrm{Hom}(\mathbf{K},\mathbf{Q})$ similarly induces the homomorphism
\begin{equation}
d^{2}{(s)}_{\thicksim}:(s|S(s))^{\ast}(T(J^{2}(N,P))|_{\Sigma^{n-p+1,0}%
(N,P)})\rightarrow\mathrm{Hom}(K(s),Q(s))\text{.}%
\end{equation}
Throughout the paper, $S(s)^{\iota}$ denotes the subset consisting of all
points $c\in S(s)$ such that the index of the quadratic form $q(s)_{c}$ is
either $\iota$ or $n-p+1-\iota$. If $\iota\neq(n-p+1)/2$, then we always
provide $Q(s)_{c}$ with the orientation such that the index of $q(s)_{c}$ is
equal to $\iota$. Furthermore, whenever $TP$ is provided with a metric,
$\mathbf{Q}_{z}$ and $Q(s)_{c}$ are always identified with a line of $(\pi
_{P}^{2})^{\ast}(TP)_{z}$ and a line of $f^{\ast}(TP)_{c}$ respectively.

A homotopy $c_{\lambda}$ with $\lambda\in\lbrack0,1]$ refers to a continuous
map $c$ of $I=[0,1]$ into a space. For example, a homotopy $h_{\lambda}$ in
$\Gamma(N,P)$ relative to a closed subset $C$ of $N$ refers to a continuous
map $h:I\rightarrow\Gamma(N,P)$ such that $h_{\lambda}|C=h_{0}|C$ for any
$\lambda$.

\section{Homotopy principle for fold-maps}

If for any section $s$ of $\Gamma(N,P)$ there exists a fold-map
$f:N\rightarrow P$ such that $j^{2}f$ is homotopic to $s$ by a homotopy in
$\Gamma(N,P)$, then we say that the homotopy principle for fold-maps
\textit{in the existence level} holds. In this section we prove Theorem 0.5 as
a consequence of Theorems 2.2, 2.4 and Proposition 2.5 below. We also prove
another version of Theorem 0.5 as follows.

\begin{theorem}
Let $n\geq p\geq2$. Let $N$, $P$ and $C$ be manifolds and a closed subset of
$N$ respectively and let $g$ be a fold-map defined on a neighborhood of $C$
into $P$. Let $s$ be a section of $\Gamma^{tr}(N,P)$ such that $j^{2}g=s$ on a
neighborhood of $C$ and that $Q(s)$ is a trivial bundle. Then there exists a
fold-map $f:N\rightarrow P$ such that $j^{2}f$ is homotopic to $s$ relative to
$C$ by a homotopy $h_{\lambda}$ in $\Gamma(N,P)$ such that $h_{0}=s$,
$h_{1}=j^{2}f$ and that $Q(j^{2}f)$ is a trivial bundle.
\end{theorem}

If the closure of no component of $N\setminus C$ is compact, then the
assertion of Theorem 0.5 is a direct consequence of [G1, Theorem 4.1.1].
Theorem 0.5 is a relative form of a special case of [An1, Theorem 1]. Its
proof given there was sketchy and the proof of Proposition 2.5 below was
abbreviated, because it seemed a reformulation of Theorem 2.2 below. However,
it turns out that Theorem 0.5 together with Theorem 2.4 is very useful as
Corollaries 0.2 and 0.4 show. This is the reason why we give here its proof in
detail for the case $n>p$. A proof of Theorem 0.5 for the case $n=p$ is given
in [An3, Theorem 4.1].

In the following definition let $S^{0},\ldots,S^{[(n-p+1)/2]}$ be submanifolds
of dimension $p-1$ of $N$, $C$ be a closed subset of $N$ and let $g$ be a
fold-map into $P$ defined on a neighborhood of $C$.

Let $\frak{M}(N,P;S^{0},\ldots,S^{[(n-p+1)/2]},C,g)$ denote the space
consisting of all fold-maps $f:N\rightarrow P$, which is equipped with the
$C^{\infty}$-topology, such that

(1) $S(j^{2}f)^{\iota}=S^{\iota}$ for $0\leq\iota\leq\lbrack(n-p+1)/2], $

(2) $f=g$ on a neighborhood of $C$.

\noindent Let $\frak{m}(N,P;S^{0},\ldots,S^{[(n-p+1)/2]},C,g)$ denote the
space consisting of all homomorphisms $h:TN\rightarrow TP$, which is equipped
with the compact-open topology, such that

(1) if $x\in N\setminus(\cup_{\iota=0}^{[(n-p+1)/2]}S^{\iota})$, then $h_{x}$
is of rank $p$,

(2) for any point $c\in S^{\iota}$, there exists a neighborhood $U_{c}$ of $c$
in $N$ and a fold-map $f_{U_{c}}:U_{c}\rightarrow P$ such that $h|TU_{c}%
=df_{U_{c}}$ and that $S(j^{2}f_{U_{c}})^{\iota}=U_{c}\cap S^{\iota},$

(3) for any point $c\in C$, we have $h_{c}=dg_{c}.$

\noindent Then we have the continuous map
\begin{equation}
\mathbf{d:}\frak{M}(N,P;S^{0},\ldots,S^{[(n-p+1)/2]},C,g)\rightarrow
\frak{m}(N,P;S^{0},\ldots,S^{[(n-p+1)/2]},C,g)
\end{equation}
defined by $\mathbf{d}(f)=df.$

The following theorem due to $\grave{\mathrm{E}}\mathrm{lia}\check{\mathrm{s}%
}\mathrm{berg}$[E2, 4.7 Theorem] will play an important role in the proof of
Theorems 0.5, 2.1 and 2.4.

\begin{theorem}
[{[}E2{]}]Let $n\geq p\geq2$. Let $N$ and $P$ be connected manifolds of
dimensions $n$ and $p$ respectively and let $S^{0},\ldots,S^{[(n-p+1)/2]}$,
$C$ and $g$ be as above. Assume that each connected component of $N\setminus
C$ has non-empty intersection with each $S^{\iota}$ $(0\leq\iota\leq
\lbrack(n-p+1)/2])$. Let $B:TN\rightarrow TP$ be any homomorphism in
$\frak{m}(N,P;S^{0},\ldots,S^{[(n-p+1)/2]},C,g)$.

Then there exist a fold-map $f:N\rightarrow P$ in $\frak{M}(N,P;S^{0}%
,\ldots,S^{[(n-p+1)/2]},C,g)$ and a homotopy of homomorphisms $B_{\lambda
}:TN\rightarrow TP$ in $\frak{m}(N,P;S^{0},\ldots,S^{[(n-p+1)/2]},C,g)$ such
that $B_{0}=B$, $B_{1}=df$ and $B_{\lambda}|_{C}=B|_{C}$ for any $\lambda$.
\end{theorem}

We begin by preparing several notions and results, which are necessary for the
proof of Theorem 0.5. For the fold-map $g$ and the closed subset $C$ in the
statement of Theorem 0.5, we take a closed neighborhood $V(C)$ of $C$ such
that $V(C)$ is an $n$-dimensional submanifold with boundary for a while and
that $g$ is defined on a neighborhood of $V(C)$, where $j^{2}g=s$. Without
loss of generality we may assume that $N\setminus\mathrm{Int}V(C)$ is
nonempty. Take a smooth function $h_{C}:N\rightarrow\lbrack0,1]$ such that
\begin{equation}
\left\{
\begin{array}
[c]{ll}%
h_{C}(x)=1 & \text{for }x\in C,\\
h_{C}(x)=0 & \text{for }x\in N\setminus\mathrm{Int}V(C),\\
0<h_{C}(x)<1 & \text{for }x\in\mathrm{Int}V(C)\setminus C.
\end{array}
\right.
\end{equation}
By the Sard Theorem ([H2]) there is a regular value $r$ of $h_{C}$ with
$0<r<1$. Then $h_{C}^{-1}(r)$ is a submanifold and we set $U(C)=h_{C}%
^{-1}([r,1])$. We decompose $N\setminus\mathrm{Int}U(C)$ into the connected
components, say $L_{1},\ldots,L_{j},\ldots$. It suffices to prove Theorem 0.5
for each $L_{j}\cup$Int$U(C)$. Since $\partial N=\emptyset$, we have that
$N\setminus U(C)$ has empty boundary. If $L_{j}$ is not compact, then Theorem
0.5 holds for $L_{j}\cup$Int$U(C)$ by Gromov's theorem ([G1, Theorem 4.1.1]).
Therefore, it suffices to consider the special case where

(C1) $C$ has closed neighborhoods $U(C)$ and $V(C)$ with $U(C)\subset
$Int$V(C)$,

(C2) $N\setminus V(C)\neq\emptyset$ and $g$ is defined on a neighborhood of
$V(C)$, where $j^{2}g=s$,

(C3) $N\setminus\mathrm{Int}U(C)$ is compact, connected and nonempty,

(C4) $\partial U(C)$ and $\partial V(C)$ are submanifolds of dimension $n-1$,

(C5) there is a smooth function $h_{C}:N\rightarrow\lbrack0,1]$ satisfying
(2.2) such that for a sufficiently small positive real number $\varepsilon$
with $r-2\varepsilon>0$, $r-t\varepsilon$ ($0\leq t\leq2$) are all regular
values of $h_{C}$, and hence $h_{C}^{-1}([r-2\varepsilon,1])$ is contained in
$V(C)$.

We set $V_{t}=h_{C}^{-1}(r-(2-t)\varepsilon)$. There exists an embedding
$e_{V}:V_{0}\times\lbrack0,2]\rightarrow N$ such that $e_{V}(V_{0}\times
t)=V_{t}$. We denote, by $N_{t}$, the union of $N\setminus\mathrm{Int}%
(h_{C}^{-1}([r-2\varepsilon,1]))$ and $e_{V}(V_{0}\times\lbrack0,t])$. In
particular, we have that $N\setminus$Int$N_{2}=U(C)$. Furthermore, we may
assume that

(C6) $s\in\Gamma^{tr}(N,P)$ and $S(s)$ is transverse to $V_{0}$ and $V_{2}$.

\begin{remark}
If $g:(N,c)\rightarrow(P,g(c))$ is a fold-map germ and if $c$ is a fold
singularity, then $d_{c}^{2}g:T_{c}N\rightarrow\mathrm{Hom}(K(j^{2}%
g)_{c},Q(j^{2}g)_{c})$ coincides with $d_{c}^{2}(j^{2}g)$ and is an
epimorphism $($see $\S1)$. Since $K(j^{2}g)_{c}\cap T_{c}(S(j^{2}g))=\{0\}$,
we may regard $K(j^{2}g)$ as the normal bundle of $S(j^{2}g)$ near $c$.
\end{remark}

Let $s$ be a section of $\Gamma^{tr}(N,P)$. Recall the homomorphisms in (1.7)
and (1.8)
\begin{align*}
d(s)_{\thicksim}  &  :TN|_{S(s)}\mathbf{\rightarrow}(s|S(s))^{\ast}%
(T(J^{2}(N,P))|_{\Sigma^{n-p+1,0}(N,P)}),\\
d^{2}(s)_{\thicksim}  &  :(s|S(s))^{\ast}(T(J^{2}(N,P))|_{\Sigma
^{n-p+1,0}(N,P)})\mathbf{\rightarrow}\mathrm{Hom}(K(s),Q(s)).
\end{align*}
Since $K(s)$ is a subbundle of $TN|_{S(s)}$, we have the inclusion map
$i_{K(s)}:K(s)\rightarrow TN|_{S(s)}$. By Remark 2.3 there exists a normal
bundle $\mathcal{N}(s)$ of $S(s)$ in $TN|_{S(s)}$ such that if $c\in
S(s)\setminus$\textrm{Int}$N_{1}$, then $\mathcal{N}(s)_{c}=K(s)_{c}$ and if
$c\in S(s)\cap\text{Int}N_{1/2}$, then $\mathcal{N}(s)_{c}$ is orthogonal to
$T_{c}(S(s)).$ Such a normal bundle is always denoted by $\mathcal{N}(s).$ If
$s\in\Gamma^{tr}(N,P)$, then $d^{2}(s)_{\thicksim}\circ d(s)_{\thicksim}$ is
an epimorphism and ${d^{2}s}|K(s)$ is an isomorphism. Then we have the
monomorphism $\Phi(s):\mathcal{N}(s)\rightarrow TN|_{S(s)}$ defined by
\begin{equation}
\Phi(s)=i_{K(s)}\circ({d^{2}s}|K(s))^{-1}\circ d^{2}(s)_{\thicksim}\circ
d(s)_{\thicksim}|\mathcal{N}(s).
\end{equation}
Let $i_{\mathcal{N}(s)}:\mathcal{N}(s)\rightarrow TN|_{S(s)}$ be the
inclusion. If $f$ is a fold-map, then it follows from the definition of the
intrinsic derivative and Remark 2.3 that $i_{\mathcal{N}(j^{2}f)}=\Phi
(j^{2}f)$, where $\mathcal{N}(j^{2}f)=K(j^{2}f)$. Let \textrm{Mono}%
$(\mathcal{N}(s),TN|_{S(s)})$ denote the subset of \textrm{Hom}$(\mathcal{N}%
(s),TN|_{S(s)})$ consisting of all monomorphisms. Then two monomorphisms
$i_{\mathcal{N}(s)}$ and $\Phi(s)$ are regarded as sections of \textrm{Mono}%
$(\mathcal{N}(s),TN|_{S(s)})$ over $S(s)$.

For a section $s\in\Gamma^{tr}(N,P)$, let $S(s)^{\iota}\setminus
\mathrm{Int}U(C)$ be decomposed into the connected components ${M(s)}%
_{1}^{\iota},\ldots,{M(s)}_{j_{\iota}}^{\iota}$ $(0\leq\iota\leq
\lbrack(n-p+1)/2])$, which may have non-empty boundary. Any one of ${M(s)}%
_{j}^{\iota}$'s will be often denoted by $M$ when we fix it. For each $M,$ the
two sections $i_{\mathcal{N}(s)}|_{M}$ and $\Phi(s)|_{M}$ satisfy
$i_{\mathcal{N}(s)}|_{\partial M}=$ $\Phi(s)|_{\partial M}$, where $\partial
M$ may be empty.

We denote the bundle of the local coefficients $\mathcal{B}(\pi_{p-1}%
(\mathrm{Mono}((\mathcal{N}(s)|_{M})_{c},(TN|_{M})_{c}))),$ $c\in M,$ by
$\mathcal{B}(\pi_{p-1})$, which is a covering space over $M$ with fiber
$\pi_{p-1}(\mathrm{Mono}(\mathcal{N}(s)_{c},T_{c}N))$ defined in [Ste, 30.1].
By the obstruction theory due to [Ste, 36.3], we have the primary difference
$d(i_{\mathcal{N}(s)}|_{M},$ $\Phi(s)|_{M})$ defined in the cohomology group
with local coefficients $H^{p-1}(M,\partial M;\mathcal{B}(\pi_{p-1})).$ Since
\textrm{Mono}$(\mathbf{R}^{n-p+1},\mathbf{R}^{n})$ is identified with
$GL(n)/GL(p-1)$, it follows from [Ste, 38.2] that
\[
\pi_{p-1}(\text{\textrm{Mono}}(\mathbf{R}^{n-p+1},\mathbf{R}^{n}%
))\cong\left\{
\begin{array}
[c]{ll}%
\mathbf{Z} & \text{if }p\text{ is odd or }n=p,\\
\mathbf{Z/}2\mathbf{Z} & \text{if }p\text{ is even and }n>p.
\end{array}
\right.
\]
We should note that if $p$ is even and $n>p$, then the bundle of the local
coefficients is trivial. If $p$ is odd or $n=p$, then we have an integer
$m(s,M)$ such that
\begin{equation}
H^{p-1}(M,\partial M;\mathcal{B}(\pi_{p-1}))\cong\mathbf{Z/}m(s,M)\mathbf{Z}.
\end{equation}
In the following brief proof of this fact we use the notation used in [Ste,
\S31]. Consider a triangulation of $(M,\partial M)$ and let $C^{i}(M,\partial
M;\mathcal{B}(\pi_{p-1}))$ be the $i$-th cochain group. Let $\sigma^{p-2}(M)$
be a $(p-2)$-simplex of $M\setminus\partial M$ which is a common face of two
$(p-1)$-simplexes $\tau_{1}^{p-1}(M)$ and $\tau_{2}^{p-1}(M)$ of $M$. Let
$a_{\sigma}$ and $a_{\tau_{i}}$ be base points of $\sigma^{p-2}(M)$ and
$\tau_{i}^{p-1}(M)$ ($i=1,2$) respectively. Let $\omega_{\sigma}\in
C^{p-2}(M,\partial M;\mathcal{B}(\pi_{p-1}))$ and $\omega_{\tau_{i}}\in
C^{p-1}(M,\partial M;\linebreak \mathcal{B}(\pi_{p-1}))$ be the dual $p-2$
cochain and the dual $p-1$ cochains of $\sigma^{p-2}(M)$ and $\tau_{i}%
^{p-1}(M)$ respectively, which map them to the respective generators of
$\pi_{p-1}(\mathrm{Mono}(\mathcal{N}(s)_{a_{\sigma}},T_{a_{\sigma}}N))$ and
$\pi_{p-1}(\mathrm{Mono}(\mathcal{N}(s)_{a_{\tau_{i}}},T_{a_{\tau_{i}}}N))$
($i=1,2$). Let $\delta:C^{p-2}(M,\partial M;\mathcal{B}(\pi_{p-1}))\rightarrow
C^{p-1}(M,\partial M;\linebreak \mathcal{B}(\pi_{p-1}))$ be the coboundary
homomorphism. Then we have that%
\[
\delta(\omega_{\sigma})(\tau_{i}^{p-1}(M))=\omega_{\sigma}(\partial(\tau
_{i}^{p-1}(M)))
\]
is equal to a generator of $\pi_{p-1}(\mathrm{Mono}(\mathcal{N}(s)_{a_{\sigma
}},T_{a_{\sigma}}N))$ and if $\tau(M)$ is one of the other $p-1$ simplexes,
then $\delta(\omega_{\sigma})(\tau(M))=0$. Therefore, we have $\omega
_{\tau_{1}}=\pm\omega_{\tau_{2}}$ modulo Im$(\delta)$. This implies (2.4).

For the proof of Theorems 0.5 and 2.1 we need the following theorem and proposition.

\begin{theorem}
{Let $n\geq p\geq2$. Let $N$, $P$, }$C$, $s$ and $g$ be ones{ given in Theorem
}$0.5${ which satisfy the assumptions }$(C1)$ to $(C6)$.{ Suppose for each
}$\iota$ {that\ $S(s)^{\iota}\setminus\mathrm{Int}V(C)$\ is a non-empty set
and that $S(s)^{\iota}\setminus\mathrm{Int}U(C)$ is decomposed into the
connected components }${M(s)}_{1}^{\iota},\ldots,{M(s)}_{j_{\iota}}^{\iota}$
$(0\leq\iota\leq\lbrack(n-p+1)/2])$.{\ Assume that }$d(i_{\mathcal{N}%
(s)}|_{_{{M(s)}_{j}^{\iota}}},$ $\Phi(s)|_{_{{M(s)}_{j}^{\iota}}})=0${\ for
all }$j$ and $\iota${. Then there exist a homotopy $s_{\lambda}$ relative to a
neighborhood of }$U(C)${ in $\Gamma(N,P)$ and a fold-map $f:N\rightarrow P$
satisfying the following. }

$(1)${\ $s_{0}=s$ and }$s_{1}=${$j^{2}f$,}

{$(2)$ $S(s_{\lambda})^{\iota}=S(s)^{\iota}$ for any $\lambda$ and }$\iota${,}

$(3)$ if $Q(s)$ is a trivial bundle, then $Q(s_{1})$ is also a trivial bundle. {\ }
\end{theorem}

\begin{proposition}
{Let $n\geq p\geq2$. Let $N$, $P$, }$C$, $s$ and $g$ be ones{ given in Theorem
}$0.5${ which satisfy the assumptions }$(C1)$ to $(C6)$. {Then there exists a
homotopy $s_{\lambda}$ relative to }$V(C)$ {in $\Gamma(N,P)$ with }$s_{0}=s$
{such that if $S(s_{1})^{\iota}\setminus\mathrm{Int}U(C)$\ is decomposed into
the nonempty connected components }${M(s}_{1}{)}_{1}^{\iota},\ldots,{M(s}%
_{1}{)}_{j_{\iota}}^{\iota}$, then

$(1)$ {$s_{1}\in\Gamma^{tr}(N,P)$, }

$(2)$ $d(i_{\mathcal{N}(s_{1})}|_{_{{M(s_{1})}_{j}^{\iota}}},$ $\Phi
(s_{1})|_{_{{M(s_{1})}_{j}^{\iota}}}){=0}${\ for all }$j$ and $\iota$,

$(3)$ ${S(s)}^{\iota}{\setminus V(C)}\neq\emptyset$ if and only if ${S(s}%
_{1}{)}^{\iota}{\setminus V(C)}\neq\emptyset$,

$(4)$ if $Q(s)$ is a trivial bundle, then $Q(s_{1})$ is also a trivial bundle.
\end{proposition}

The proofs of Theorem 2.4 and Proposition 2.5 will be given in $\S4$ and $\S5
$ respectively.

Here we give a proof of Theorems 0.5 and 2.1. In the proof we use the notation
explained above.

\begin{proof}
[Proof of Theorems 0.5 and 2.1.]By the above argument, it suffices to prove
Theorem 0.5 for the case where (C1) to (C6) are satisfied. If there exists an
$\iota$ such that $S(s)^{\iota}\setminus V(C)=\emptyset$, then we have a
neighborhood $\mathbf{U}$ in some local chart of $N$ with
Cl$(\mathbf{U)\subset}N\setminus(V(C)\cup S(s))$ and a homotopy $s_{\lambda
}^{\prime}\in\Gamma(N,P)$ relative to $N\setminus\mathbf{U}$\ with
$s_{0}^{\prime}=s$ and $s_{1}^{\prime}\in\Gamma^{tr}(N,P)$, which satisfies,
for each $\iota$ with $S(s)^{\iota}\setminus V(C)=\emptyset$, the following.

(1) $S(s_{1}^{\prime})^{\iota}\cap\mathbf{U}$ is nonempty and is the boundary
of an embedded $p$-disk within $\mathbf{U}$; these $p$-disks are mutually disjoint.

(2) $Q(s_{1}^{\prime})|_{S(s_{1}^{\prime})^{\iota}\cap\mathbf{U}}$ and the
normal bundle $\mathcal{N}_{S(s_{1}^{\prime})^{\iota}\cap\mathbf{U}}$ of
$S(s_{1}^{\prime})^{\iota}\cap\mathbf{U}$\ are extended to bundles over the
respective $p$-disks. In particular, they are trivial bundles.

It follows from Proposition 2.5 for $s_{1}^{\prime}$ that there exists a
homotopy $s_{\lambda}^{\prime\prime}$ relative to a neighborhood of $V(C)$ in
$\Gamma(N,P)$ with $s_{0}^{\prime\prime}=s_{1}^{\prime}$ such that
$s_{1}^{\prime\prime}\in\Gamma^{tr}(N,P)$ and $d(i_{\mathcal{N}(s_{1}%
^{\prime\prime})}|_{_{{M(s}^{\prime\prime}{_{1})}_{j}^{\iota}}},$ $\Phi
(s_{1}^{\prime\prime})|_{_{{M(s}_{1}^{\prime\prime}{)}_{j}^{\iota}}})=0${\ for
all }$j$ and $\iota$ and that $S(s_{1}^{\prime\prime})^{\iota}\setminus
V(C)\neq\emptyset$ for all $\iota$. By applying Theorem 2.4 to the section
$s_{1}^{\prime\prime}$, we see that there exists a homotopy ${s}_{\lambda
}^{\prime\prime\prime}$ relative to {a neighborhood of }$U(C)$ in
$\Gamma(N,P)$ with ${s}_{0}^{\prime\prime\prime}={s}_{1}^{\prime\prime}$ such
that there exists a fold-map $f:N\rightarrow P$ with $j^{2}f=s_{1}%
^{\prime\prime\prime}$. Thus we obtain a required homotopy $s_{\lambda}$ by
pasting ${s}_{\lambda}^{\prime}$, ${s}_{\lambda}^{\prime\prime}$ and
${s}_{\lambda}^{\prime\prime\prime}$. This proves Theorem 0.5.

Next we prove Theorem 2.1 by using the above notation. Since $Q(s)$ is
trivial, $Q(s_{1}^{\prime})$ is also trivial by (2) above. It follows from
Proposition 2.5 (4) that $Q(s_{1}^{\prime\prime})$ is trivial. Therefore,
$Q(s_{1}^{\prime\prime\prime})$ is also trivial by Theorem 2.4 (3). This
proves Theorem 2.1.
\end{proof}

\section{Proof of Theorems 0.1 and 0.3}

In this section we prove Theorems 0.1 and 0.3.

The following lemma is a consequence of the fundamental properties of
fold-maps ([E1, 3.8 and 3.9] and [Sa1, Lemma 3.1]).

\begin{lemma}
Let $n\geq p\geq1$. Let $f:N\rightarrow P$ be a fold-map. Assume that either
$n-p+1$ is odd, or if $n-p+1$ is even, then $Q(j^{2}f)$ is a trivial bundle.
Then there exists an epimorphism $TN\oplus\theta_{N}\rightarrow TP$ covering
$f$.
\end{lemma}

\begin{proof}
Recall the quadratic form $q(j^{2}f):S^{2}(K(j^{2}f))\rightarrow Q(j^{2}f)$.
If $n-p+1$ is odd, then we can provide $Q(j^{2}f)$ with an orientation so that
this index is less than $(n-p+1)/2$. Consequently, $Q(j^{2}f)$ is a trivial
bundle. Choose a Riemannian metric on $P$. Since $Q(j^{2}f)$ is the cokernel
bundle of $df|_{S(j^{2}f)}:TN|_{S(j^{2}f)}\rightarrow TP$, we may regard
$Q(j^{2}f)$ as a subbundle of $f^{\ast}(TP)$ and obtain a monomorphism
$\varphi^{\prime}:Q(j^{2}f)\rightarrow TP$ covering $f|S(j^{2}f)$. Let
$\varphi:\theta_{N}\rightarrow TP$ covering $f$ be any extended homomorphism
such that $\varphi|Q(j^{2}f)=\varphi^{\prime}$. Then define the homomorphism
$h:TN\oplus\theta_{N}\rightarrow TP$ by $h(\mathbf{v}_{1}\oplus\mathbf{v}%
_{2})=df(\mathbf{v}_{1})+\varphi(\mathbf{v}_{2})$. It is easy to see that $h$
is an epimorphism.
\end{proof}

If $n-p+1$ is even, then $Q(j^{2}f)$ is not necessarily trivial. The
motivation of Theorem 0.1 is [E1, 3.10], [An3, Theorem 4.8] and the above
lemma. We prove the following refined form of Theorem 0.1.

\begin{theorem}
Let $n\geq p\geq2$. Let $f:N\rightarrow P$ be a continuous map. Assume that
there exists an epimorphism $h:TN\oplus\theta_{N}\rightarrow TP$ covering $f$.
Then there exists a fold-map $g:N\rightarrow P$ homotopic to $f$ such that
$Q(j^{2}g)$ is a trivial bundle.
\end{theorem}

\begin{proof}
[Proof of Theorem 3.2]We set $\xi=$Ker$(h)$. Let $i_{\xi}:\xi\rightarrow
TN\oplus\theta_{N}$ be the inclusion. Let $\pi_{\theta_{N}}:$ $TN\oplus
\theta_{N}\rightarrow\theta_{N}$ be the canonical projection. Then we have the
homomorphism $\pi_{\theta_{N}}\circ i_{\xi}:\xi\rightarrow\theta_{N}$, which
we regard as a section of Hom($\xi,\theta_{N})$ over $N$. We deform $h$ by
homotopy to a smooth epimorphism, denoted by the same letter $h$, such that
$\pi_{\theta_{N}}\circ i_{\xi}$ is transverse to the zero-section of
Hom($\xi,\theta_{N})$. Let $V$ be the inverse image of the zero-section of
Hom($\xi,\theta_{N})$ by $\pi_{\theta_{N}}\circ i_{\xi}$. Then the normal
bundle of $V$ in $N$ is isomorphic to Hom($\xi|_{V},\theta_{V})$. We now
consider the homomorphism $h|TN:TN\rightarrow TP$.

Let $x\in N\setminus V$. For any vector $\mathbf{v\in}T_{f(x)}P$, let
$\mathbf{v=}h(\mathbf{w}_{1}\oplus\mathbf{w}_{2})$ with $\mathbf{w}_{1}\in
T_{x}N$ and $\mathbf{w}_{2}\in(\theta_{N})_{x}$. Since $\pi_{\theta_{N}}\circ
i_{\xi}:\xi_{x}\rightarrow(\theta_{N})_{x}$ is surjective, there exists a
vector $\mathbf{v}_{1}\in\xi_{x}$ with $\pi_{\theta_{N}}\circ i_{\xi
}(\mathbf{v}_{1})=\mathbf{w}_{2}$. Let $\mathbf{v}_{1}=\mathbf{w}_{3}%
\oplus\mathbf{w}_{2}$ with $\mathbf{w}_{3}\in T_{x}N$. Then we have
\begin{align*}
h(\mathbf{w}_{1}-\mathbf{w}_{3})  &  =h(\mathbf{w}_{1}\oplus\mathbf{w}%
_{2}-\mathbf{w}_{3}\oplus\mathbf{w}_{2})\\
&  =\mathbf{v}-h(\mathbf{v}_{1})\\
&  =\mathbf{v}\text{.}%
\end{align*}
Hence, $h|T_{x}N$ is surjective onto $T_{f(x)}P$.

Let $x\in V$. Since $\pi_{\theta_{N}}\circ i_{\xi}$ is a null-homomorphism and
since the kernel of $\pi_{\theta_{N}}$ is $T_{x}N$, we have $\xi_{x}\subset
T_{x}N$. Since \textrm{Ker}$(h)=\xi$, we have \textrm{Ker}$(h|T_{x}N)=\xi_{x}%
$. Since $\xi$ is of dimension $n-p+1$, we obtain that $h|T_{x}N$ is of rank
$p-1$. Therefore, $h|TN:TN\rightarrow TP$ induces a homomorphism
$H:TN\rightarrow f^{\ast}(TP)$, which is regarded as a section $s^{\prime
}:N\rightarrow\Omega^{n-p+1}(N,P)$. Next we show that the line bundle
\textrm{Cok}$(h|TN)|_{V}$ is isomorphic to $\theta_{V}$. In fact, the exact
sequence
\[
0\rightarrow\xi|_{V}\rightarrow TN|_{V}\rightarrow f^{\ast}(TP)|_{V}%
\rightarrow\mathrm{Cok}(h|TN)|_{V}\rightarrow0
\]
yields
\begin{align*}
W_{1}(\mathrm{Cok}(h|TN)|_{V})  &  =W_{1}(f^{\ast}(TP)|_{V})-W_{1}%
(TN|_{V})+W_{1}(\xi|_{V})\\
&  =W_{1}(\theta_{V}),
\end{align*}
where $W_{1}$ refers to the first Stiefel-Whitney class. Since $\mathrm{Cok}%
(h|TN)|_{V}$ and $\theta_{V}$ are line bundles, we obtain the assertion. By
choosing a Riemannian metric on $P$, we regard $\mathrm{Cok}(h|TN)|_{V}$ as a
subbundle of $f^{\ast}(TP)|_{V}$.

We take any non-singular symmetric map $\mathbf{b}:\xi|_{V}\otimes\xi
|_{V}\rightarrow\mathrm{Cok}(h|TN)|_{V}$ $(\subset f^{\ast}(TP)|_{V}$). Let
$\widetilde{\mathbf{b}}:S^{2}(TN)\rightarrow f^{\ast}(TP)$ be any extended
homomorphism of $\mathbf{b}$. Under the identification (1.2), we define the
section $s:N\rightarrow\Omega^{n-p+1,0}(N,P)$ as follows. For $x\in N$, set
\[
s(x)=\{x,f(x),H_{x},\widetilde{\mathbf{b}}_{x}\}\text{.}%
\]
If $x\notin V$, then $H_{x}$ is nonsingular and if $x\in V$, then
Ker$(H_{x})=\xi_{x}$ and $\mathbf{b}_{x}$ is nonsingular.

By Theorem 2.1 there exists a fold-map $g:N\rightarrow P$ such that $j^{2}g$
and $s$ are homotopic as sections in $\Gamma(N,P)$ and that $Q(j^{2}g)$ is a
trivial bundle. In particular, $g$ is homotopic to $f$, and the theorem is
proved. Hence, Theorem 0.1 is also proved.
\end{proof}

Next we prepare tools necessary for the proof of Theorem 0.3. Let
$\Delta(d_{1},\ldots,d_{k})$ be the diagonal $k\times k$ matrix with diagonal
components $(d_{1},\ldots,d_{k})$. In $\mathbf{R}^{k}$ for any number $k$, the
vector ${}^{t}(\overbrace{0,\ldots,0}^{i-1},1,\overbrace{0,\ldots,0}^{k-i})$
is denoted by $\mathbf{e}_{i}$. For a point $y=$ $(y_{1},\ldots,y_{k}%
)\in\mathbf{R}^{k}$, let $\mathcal{T}^{k}(y)$ denote the $k\times k$ matrix
$(\delta_{ij}-2y_{i}y_{j})_{1\leq i,j\leq k}$, where $\delta_{ij}$ refers to
the Kronecker delta. Let $\mathcal{T}^{k}:\theta_{\mathbf{R}^{k}}%
^{k}\rightarrow\theta_{\mathbf{R}^{k}}^{k}$ denote the homomorphism defined by
$\mathcal{T}^{k}(y,\mathbf{v})=(y,\mathcal{T}^{k}(y)(\mathbf{v}))$. The
formula $\det\mathcal{T}^{k}(y)=1-2(y_{1}^{2}+\cdots+y_{k}^{2})$ is well known
as an exercise in linear algebra. For simplicity we write $\mathcal{T}(y)$ and
$\mathcal{T}$ for $\mathcal{T}^{n}(y)$ and $\mathcal{T}^{n}$ respectively. Let
$\Delta\lbrack a;j]$ be $\Delta(\overbrace{a,\ldots,a}^{n-j},\overbrace
{-a,\ldots,-a}^{j})$. In particular, for $\varepsilon=1$ or $0$\ we denote
$\Delta\lbrack1;\varepsilon]$ by $I_{-\varepsilon}^{n}$. Let $\mathcal{T}%
I_{-\varepsilon}^{n}$ (resp. $I_{-1}^{n}\mathcal{T}$)$:\theta_{\mathbf{R}^{n}%
}^{n}\rightarrow\theta_{\mathbf{R}^{n}}^{n}$ be the homomorphism defined by
$\mathcal{T}I_{-\varepsilon}^{n}(y,\mathbf{v})=(y,\mathcal{T}%
(y)I_{-\varepsilon}^{n}(\mathbf{v}))$ (resp. by $I_{-1}^{n}\mathcal{T}%
(y,\mathbf{v})=(y,I_{-1}^{n}\mathcal{T}(y)(\mathbf{v}))$) and $\mathcal{T}%
(y)I_{-\varepsilon}^{n}(_{p}^{1})$ the $p\times n$ matrix consisting of the
first $p$ row vectors of $\mathcal{T}(y)I_{-\varepsilon}^{n}$. Then we have
the homomorphisms $\mathcal{T}I_{-\varepsilon}^{n}(_{p}^{1}):\theta
_{\mathbf{R}^{n}}^{n}\rightarrow\theta_{\mathbf{R}^{p}}^{p}$ and
$\mathbf{T}I_{-\varepsilon}^{n}\mathcal{(}_{p}^{1}):\theta_{\mathbf{R}^{n}%
}^{n}\rightarrow\theta_{\mathbf{R}^{n}}^{p}$\ defined by $\mathcal{T}%
I_{-\varepsilon}^{n}\mathcal{(}_{p}^{1})(y,\mathbf{v})=((y_{1},\ldots
,y_{p}),(\mathcal{T}(y)I_{-\varepsilon}^{n}(_{p}^{1}))\mathbf{v})$ and
$\mathbf{T}I_{-\varepsilon}^{n}\mathcal{(}_{p}^{1})(y,\mathbf{v}%
)=(y,(\mathcal{T}(y)I_{-\varepsilon}^{n}(_{p}^{1}))\mathbf{v})$ respectively,
where $y=(y_{1},\ldots,y_{n})\in\mathbf{R}^{n}$. Let $\mathbf{0}_{n-p}$ be the
null vector of degree $n-p$. The subset of $\mathbf{R}^{n}$ consisting of all
points $y$ such that $\mathcal{T}(y)I_{-\varepsilon}^{n}(_{p}^{1})$ is of rank
$p-1$ coincides with $S_{1/\sqrt{2}}^{p-1}\times\mathbf{0}_{n-p}$, which is
the subspace consisting of all points $y=$ $(y_{1},\ldots,y_{n})$ such that
$2(y_{1}^{2}+\cdots+y_{p}^{2})=1$ and $y_{p+1}=\cdots=y_{n}=0$ (see Lemma 3.3
below). We denote the intrinsic derivatives of $\mathbf{T}I_{-\varepsilon}%
^{n}\mathcal{(}_{p}^{1})$ for $n>p$ and $I_{-1}^{n}\mathcal{T}$ for $n=p$ due
to I. R. Porteous (see [B, Lemma 7.4]) by $\mathbf{d(T}I_{-\varepsilon}%
^{n}\mathcal{(}_{p}^{1}))$ and $\mathbf{d}(I_{-1}^{n}\mathcal{T)}$
respectively on $S_{1/\sqrt{2}}^{p-1}\times\mathbf{0}_{n-p}$. In the following
lemma we set $\mathbf{c}={}^{t}(c_{1},\ldots,c_{p},0,\ldots,0)$,
$\overset{\bullet}{\mathbf{c}}={}^{t}(c_{1},\ldots,c_{p})$ and $\mathbf{c}%
^{-}={}^{t}(c_{1},\ldots,c_{n-1},-c_{n})$.

\begin{lemma}
$(1)$ $\mathcal{T}(y)I_{-\varepsilon}^{n}\mathcal{(}_{p}^{1})$ $($resp.
$I_{-1}^{n}\mathcal{T}(y))$ is of rank $p$ $($resp. rank $n)$ if and only if
$y\notin S_{1/\sqrt{2}}^{p-1}\times\mathbf{0}_{n-p}$ $($resp. $y\notin
S_{1/\sqrt{2}}^{n-1}).$

\noindent$(2)$ If $c\in S_{1/\sqrt{2}}^{p-1}\times\mathbf{0}_{n-p}$, then we have

$(2$\textrm{-}$\mathrm{i)}$ $\mathcal{T}(c)I_{-\varepsilon}^{n}(_{p}^{1})$
$($resp. $I_{-1}^{n}\mathcal{T}(c))$ is of rank $p-1$ $($resp. $n-1)$,

$(2$\textrm{-}$\mathrm{ii)}$ the kernel of $\mathcal{T}(c)I_{-\varepsilon}%
^{n}\mathcal{(}_{p}^{1})$ is generated by the vectors $\mathbf{c}$,
$\mathbf{e}_{p+1}$,..., $\mathbf{e}_{n}$ $(n>p)$,

$(2$\textrm{-}$\mathrm{iii})$ the kernel of $I_{-1}^{n}\mathcal{T}(c)$ is
generated by the vector $\mathbf{c}$ $(n=p)$,

$(2$\textrm{-}$\mathrm{iv)}$ the cokernel of $\mathcal{T}(c)I_{-\varepsilon
}^{n}\mathcal{(}_{p}^{1})$ is generated by $\overset{\bullet}{\mathbf{c}}$
$(n>p)$,

$(2$\textrm{-}$\mathrm{v})$ the cokernel of $I_{-1}^{n}\mathcal{T}(c)$ is
generated by $\mathbf{c}^{-}$ $(n=p)$,

$(2\mathrm{-vi)}$ the intrinsic derivatives satisfy
\[%
\begin{array}
[c]{ll}%
(\mathbf{d}_{c}(\mathbf{T}I_{-\varepsilon}^{n}\mathcal{(}_{p}^{1}%
))\mathcal{(}\mathbf{c}))(\mathbf{c})=2(-\overset{\bullet}{\mathbf{c}}), & \\
(\mathbf{d}_{c}(\mathbf{T}I_{-\varepsilon}^{n}\mathcal{(}_{p}^{1}%
))\mathcal{(}\mathbf{c}))(\mathbf{e}_{k})=\mathbf{0} & \text{for }n\geq k\geq
p+1,\\
(\mathbf{d}_{c}(\mathbf{T}I_{-\varepsilon}^{n}\mathcal{(}_{p}^{1}%
))\mathcal{(}\mathbf{e}_{k}))(\mathbf{e}_{\ell})=2\delta_{k\ell}%
(-\overset{\bullet}{\mathbf{c}}) & \text{for }n>k\geq p+1\text{ and }n\geq
\ell\geq p+1,\\
(\mathbf{d}_{c}(\mathbf{T}I_{-\varepsilon}^{n}\mathcal{(}_{p}^{1}%
))\mathcal{(}\mathbf{e}_{k}))(\mathbf{e}_{\ell})=(-1)^{\varepsilon\delta_{kn}%
}2\delta_{k\ell}(-\overset{\bullet}{\mathbf{c}}) & \text{for }n=k\geq\ell\geq
p+1,\\
(\mathbf{d}_{c}(I_{-1}^{n}\mathcal{T})\mathcal{(}\mathbf{c}))(\mathbf{c}%
)=2(-\mathbf{c}^{-}) & \text{for }n=p.
\end{array}
\]
\end{lemma}

\begin{proof}
Suppose that $n>p$ and one of $y_{p+1},\ldots,y_{n}$ is not $0$, say
$y_{p+1}\neq0$. Then we have
\[
\text{rank}\left(
\begin{array}
[c]{lll}%
&  & -2y_{1}y_{p+1}\\
(\delta_{ij}-2y_{i}y_{j}) &  & \text{ \ \ \ \ }\vdots\\
&  & -2y_{p}y_{p+1}%
\end{array}
\right)  =\text{rank}\left(
\begin{array}
[c]{lll}%
&  & y_{1}\\
E_{p} &  & \vdots\\
&  &  y_{p}%
\end{array}
\right)  =p\text{,}%
\]
where $E_{p}$ is the unit matrix of rank $p$. Next suppose that $y_{p+1}%
=\cdots=y_{n}=0$. Since $\det\mathcal{T}^{p}(y_{1},\ldots,y_{p})=1-2(y_{1}%
^{2}+\cdots+y_{p}^{2})$, $\mathcal{T}(y)I_{-\varepsilon}^{n}\mathcal{(}%
_{p}^{1})$ is of rank $p$ if and only if $y\notin S_{1/\sqrt{2}}^{p-1}%
\times\mathbf{0}_{n-p}$.

If $c\in S_{1/\sqrt{2}}^{p-1}\times\mathbf{0}_{n-p}$, then one of
$c_{1},\ldots,c_{p}$ is not $0$, say $c_{p}\neq0$. Then we have
\[
\det\mathcal{T}^{p-1}(c_{1},\ldots,c_{p-1})=1-2(c_{1}^{2}+\cdots+c_{p-1}%
^{2})=2c_{p}^{2}\neq0\text{.}%
\]
Hence, $\mathcal{T}(c)\mathcal{(}_{p}^{1})$ is of rank $p-1$. This proves (1)
and (2-i) for $n>p$. The assertions (1) and (2-i) for $n=p$ and (2-ii, iii,
iv, v) follow from a direct calculation.

We prove the assertion (2-vi). Let $\mathbf{0}_{i\times j}$ be the $i\times j$
null matrix. If $p+1\leq k\leq n$ and $p+1\leq\ell\leq n$,\ then we have
\[
\partial/\partial y_{k}(\mathcal{T}(y)I_{-\varepsilon}^{n}(_{p}^{1}))=\left(
\begin{array}
[c]{lll}%
\mathbf{0}_{p\times(k-1)} &
\begin{array}
[c]{l}%
-2(-1)^{\varepsilon\delta_{kn}}y_{1}\\
\text{ \ \ \ }\vdots\\
-2(-1)^{\varepsilon\delta_{kn}}y_{p}%
\end{array}
& \mathbf{0}_{p\times(n-k)}%
\end{array}
\right)  ,
\]
and hence,
\[
(\partial/\partial y_{k}(\mathcal{T}(y)I_{-\varepsilon}^{n}(_{p}%
^{1})))(\mathbf{e}_{\ell})={}^{t}(-2(-1)^{\varepsilon\delta_{kn}}\delta
_{k\ell}y_{1},\ldots,-2(-1)^{\varepsilon\delta_{kn}}\delta_{k\ell}y_{p}).
\]
Hence, we have
\[
(d_{c}(\mathbf{T}I_{-\varepsilon}^{n}(_{p}^{1}))(\mathbf{e}_{k}))(\mathbf{e}%
_{\ell})=2(-1)^{\varepsilon\delta_{kn}}\delta_{k\ell}(-\overset{\bullet
}{\mathbf{c}}).
\]
Furthermore, consider the curve $tc=t(c_{1},\ldots,c_{p},0,\ldots,0)$
parametrized by $t$. Then we have
\[
\partial/\partial t(\mathcal{T}(tc)I_{-\varepsilon}^{n}(_{p}^{1}))=\left(
\begin{array}
[c]{ll}%
\left(  -4tc_{i}c_{j}\right)  _{1\leq i\text{, }j\leq p} & \mathbf{0}%
_{p\times(n-p)}%
\end{array}
\right)  ,
\]
and hence, $(\mathbf{d}_{c}(\mathbf{T}I_{-\varepsilon}^{n}\mathcal{(}_{p}%
^{1}))\mathcal{(}\mathbf{c}))(\mathbf{c})=2(-\overset{\bullet}{\mathbf{c}})$
and $(\mathbf{d}_{c}(\mathbf{T}I_{-\varepsilon}^{n}\mathcal{(}_{p}%
^{1}))\mathcal{(}\mathbf{c}))(\mathbf{e}_{k})=\mathbf{0}$ for $p+1\leq k\leq
n$ follow from%
\[%
\begin{array}
[c]{ll}%
(\partial/\partial t(\mathcal{T}(tc)I_{-\varepsilon}^{n}(_{p}^{1}%
))|_{t=1})(\mathbf{c})=-4\Vert\mathbf{c\Vert}^{2}\overset{\bullet}{\mathbf{c}%
}=2(-\overset{\bullet}{\mathbf{c}}), & \text{and}\\
(\partial/\partial t(\mathcal{T}(tc)I_{-\varepsilon}^{n}(_{p}^{1}%
))|_{t=1})(\mathbf{e}_{k})=\mathbf{0} & \text{for }p+1\leq k\leq n,
\end{array}
\]
respectively.

If $n=p$, then we have, for $I_{-1}^{n}\mathcal{T}$,
\[
\partial/\partial t(I_{-1}^{n}\mathcal{T}(tc))=\left(
\begin{array}
[c]{l}%
\text{ }(-4tc_{i}c_{j})_{\substack{1\leq i\leq n-1\\1\leq\text{ }j\leq n}}\\
4tc_{n}c_{1},\cdots,4tc_{n}c_{n}%
\end{array}
\right)
\]
and hence $(\mathbf{d}_{c}(I_{-1}^{n}\mathcal{T})\mathcal{(}\mathbf{c}%
))(\mathbf{c})=2(-\mathbf{c}^{-})$ follows from%
\[
(\partial/\partial t(I_{-1}^{n}\mathcal{T}(tc))|_{t=1})(\mathbf{c}%
)=-4\Vert\mathbf{c\Vert}^{2}\mathbf{c}^{-}=2(-\mathbf{c}^{-}).
\]

This shows the assertion.
\end{proof}

It will be plausible to search fold-maps which have spheres as singularities
as in Theorem 0.3.

\begin{proof}
[Proof of Theorem 0.3]We may assume that $N$ is connected. Let $\widetilde
{h^{\prime}}:TN\oplus\theta_{N}\rightarrow f^{\ast}(TP)\oplus\xi\oplus
\theta_{N}$ over $N$ be an isomorphism. Take a local chart $U_{P}$ of $P$
diffeomorphic to\textbf{\ }$\mathbf{R}^{p}$ with $f^{-1}(U_{P})\neq\emptyset$,
where we have a trivialization $TP|_{U_{P}}\cong\theta_{U_{P}}^{p}$. Take a
point $o\in f^{-1}(U_{P})\subset N$\ and an embedding $e:\mathbf{R}%
^{n}\rightarrow N$ with $e(0)=o$ such that $e(\mathbf{R}^{n})$ is contained in
a contractible neighborhood $U$ of $o$ in $f^{-1}(U_{P})$, which is a local
chart of $N$ and $\xi$. Set $\frak{D}=e(D_{1}^{n})$, where $D_{1}^{n}$ is the
unit disk in $\mathbf{R}^{n}$. Since $O(n+1)/(O(n)\times1)$ is $(n-1)$%
-connected, there exists a bundle map $h^{\prime}:TN|_{N\setminus
\text{Int}\frak{D}}\rightarrow(f^{\ast}(TP)\oplus\xi)|_{N\setminus
\text{Int}\frak{D}}$ such that $h^{\prime}\oplus id_{\theta_{N\setminus
\text{Int}\frak{D}}}$ is homotopic to $\widetilde{h^{\prime}}|_{N\setminus
\text{Int}\frak{D}}$ by the obstruction theory. If $\pi_{f^{\ast}(TP)}%
:f^{\ast}(TP)\oplus\xi\rightarrow f^{\ast}(TP)$ and $\pi_{\xi}:f^{\ast
}(TP)\oplus\xi\rightarrow\xi$ are the canonical projections, then we have
$h^{\prime}=(\pi_{f^{\ast}(TP)}\circ h^{\prime})\oplus(\pi_{\xi}\circ
h^{\prime}).$

By applying the Phillips Submersion Theorem ([P]) for $\pi_{f^{\ast}(TP)}\circ
h^{\prime}$, we obtain a submersion $g:N\setminus$\textrm{Int}$\frak{D}%
\rightarrow P$ and $dg:T(N\setminus$\textrm{Int}$\frak{D})\rightarrow TP$
induces a bundle map $\widetilde{dg}:T(N\setminus$\textrm{Int}$\frak{D}%
)\rightarrow g^{\ast}(TP)$ such that $\pi_{f^{\ast}(TP)}\circ h^{\prime}$ and
$\widetilde{dg}$ are homotopic. Hence, $\widetilde{dg}\oplus(\pi_{\xi}\circ
h^{\prime}):TN|_{N\setminus\text{Int}\frak{D}}\rightarrow g^{\ast}%
(TP)\oplus\xi|_{N\setminus\text{Int}\frak{D}}$ is homotopic to $h^{\prime}$.
Therefore, we may assume the following by a routine argument:

(i) $f|(N\setminus$Int$\frak{D)}=g$.

(ii) We have a bundle map $h:TN|_{N\setminus\text{Int}\frak{D}}\rightarrow
g^{\ast}(TP)\oplus\xi|_{N\setminus\text{Int}\frak{D}}$ and an isomorphism
$\widetilde{h}:TN\oplus\theta_{N}\rightarrow f^{\ast}(TP)\oplus\xi\oplus
\theta_{N}$ such that $h=\widetilde{dg}\oplus(\pi_{\xi}\circ h^{\prime})$
(note that $(\pi_{\xi}\circ h^{\prime})=(\pi_{\xi}\circ h)$) and
$\widetilde{h}|_{N\setminus\text{Int}\frak{D}}=h\oplus id_{\theta
_{N\setminus\text{Int}\frak{D}}}$.

Since $U$ is contractible, we have trivializations $t_{N}:TN|_{U}%
\rightarrow\theta_{U}^{n}$, $t_{P}:f^{\ast}(TP)|_{U}\rightarrow\theta_{U}^{p}$
and $t_{\xi}:\xi|_{U}\rightarrow\theta_{U}^{n-p}$. We consider the bundle map
\[
\Psi=(t_{P}\oplus t_{\xi})\circ h\circ t_{N}^{-1}|_{U\setminus
\text{\textrm{Int}}\frak{D}}:\theta_{U\setminus\text{\textrm{Int}}\frak{D}%
}^{n}\rightarrow\theta_{U\setminus\text{\textrm{Int}}\frak{D}}^{n}\text{,}%
\]
which we regard as\ a map $U\setminus$\textrm{Int}$\frak{D}\rightarrow
GL^{+}(n)$ by choosing appropriate orientations of $TN|_{U}$ and $f^{\ast
}(TP)|_{U}$ ($GL^{+}(n)$ refers to the group of all $n\times n$ regular
matrices with positive determinants). Note that $\Psi\oplus id_{\theta
_{U\setminus\text{\textrm{Int}}\frak{D}}}$ is extended to the bundle map
\[
\widetilde{\Psi}=(t_{P}\oplus t_{\xi}\oplus id_{\theta_{U}})\circ\widetilde
{h}\circ(t_{N}\oplus id_{\theta_{U}})^{-1}\text{,}%
\]
which is regarded as a map $U\rightarrow GL^{+}(n+1)$. Let $i_{GL^{+}}%
:GL^{+}(n)\rightarrow GL^{+}(n)\times1\subset GL^{+}(n+1)$ be the inclusion.
By construction, we have
\begin{equation}
\widetilde{\Psi}|_{U\setminus\text{\textrm{Int}}\frak{D}}=i_{GL^{+}}\circ\Psi.
\end{equation}
In the following, $[\ast]$ expresses the homotopy class represented by $\ast$.
It follows from (3.1) that $(i_{GL^{+}})_{\ast}:\pi_{n-1}(GL^{+}%
(n))\rightarrow\pi_{n-1}(GL^{+}(n+1))$ maps $[\Psi|_{\partial\frak{D}}]$ to
$0$.

We use the map $\mathbf{i}^{j}:\partial D_{1}^{n}\rightarrow\partial D_{1}%
^{n}$ defined by $(y_{1},\ldots,y_{n})\rightarrow(y_{1},\ldots,y_{n-1}%
,(-1)^{j}y_{n})$ ($0\leq j\leq n-p+1)$. In particular, we have that
\[
(I_{-1}^{n}\mathcal{T})\circ\mathbf{i}^{1}(y)=I_{-1}^{n}I_{-1}^{n}%
\mathcal{T}(y)I_{-1}^{n}=\mathcal{T}(y)I_{-1}^{n}.
\]
This implies that if $j$ is odd, then%
\[
\lbrack\mathcal{T}I_{-1}^{n}|_{\partial D_{1}^{n}}]=[(I_{-1}^{n}%
\mathcal{T})\circ\mathbf{i}^{j}|_{\partial D_{1}^{n}}]=-[I_{-1}^{n}%
\mathcal{T}|_{\partial D_{1}^{n}}]
\]
in $\pi_{n-1}(SO(n))$. Let $p^{S}:SO(n)\rightarrow S^{n-1}$ be the map defined
by $p^{S}(T)=T\mathbf{e}_{n}$. We have the homomorphisms $\partial_{n}:\pi
_{n}(S^{n-1})\rightarrow\pi_{n-1}(SO(n-1))$ and $p_{\ast}^{S}:\pi
_{n-1}(SO(n))\rightarrow\pi_{n-1}(S^{n-1})$. Then the kernel of $(i_{GL^{+}%
})_{\ast}$ is described as follows (see [Ste, 23.2 Theorem and 23.4 Theorem]
and [W, Proposition 4]). If $n=1,3,7$, then $(i_{GL^{+}})_{\ast}$ is an
isomorphism and $\pi_{2}(SO(3))\cong\pi_{6}(SO(7))\cong\{0\}$. If $n$ is odd
and $n\neq1,3,7$ (resp. even), then the kernel of $(i_{GL^{+}})_{\ast}$ is
generated by
\[
\lbrack I_{-1}^{n}\mathcal{T}|_{\partial D_{1}^{n}}]=-[\mathcal{T}I_{-1}%
^{n}|_{\partial D_{1}^{n}}],
\]
which is of order two (resp. $\infty$, and is mapped onto the twice of a
generator of $\pi_{n-1}(S^{n-1})$ by $p_{\ast}^{S}$). Hence, by identifying
$\partial\frak{D}$ and $\partial D_{1}^{n}$\ with $S^{n-1},$ we have an
integer $m$ such that $[\Psi|_{\partial\frak{D}}]=m[I_{-1}^{n}\mathcal{T}%
|_{\partial D_{1}^{n}}]$, where $m=0$ for $n=3,7$ and $m=0$ or $1$\ for other
odd numbers $n$ greater than $2$.

Let $\ell$ be a non-negative integer such that $\ell\geq\lbrack(n-p+1)/2]$,
which will be defined below. By identifying \textrm{Int}$\frak{D}$ with
$\mathbf{R}^{n}$ coordinated by $(z_{1},\ldots,z_{n})$, we choose disjointly
embedded $n$-disks $\mathbf{E}_{j}$ ($0\leq j\leq\ell$) with center $o_{j}$
and radius $1$. We represent a point $x\in\mathbf{E}_{j}$ by the coordinates
\[
{y}(x)=(y_{1}(x),\ldots,y_{n}(x))=(z_{1}(x),\ldots,z_{n}(x))-o_{j}\text{.}%
\]

Next we extend $\widetilde{dg}$ to a homomorphism $H:TN\rightarrow f^{\ast
}(TP)$ such that $H$ is of rank $\geq p-1$ and is of rank $p-1$ exactly on a
finite number of spheres embedded in $\mathbf{E}_{j}$'s. Consider the
following cases to define $\Psi^{\mathbf{E}}:\cup_{j=0}^{\ell}\mathbf{E}%
_{j}\rightarrow{\mathrm{Hom}}(\mathbf{R}^{n},\mathbf{R}^{n})$, where
$\mathbf{G}$ refers to $[(n-p+1)/2]$ for short.

\textbf{Case I}: $n>p\geq2$ and $n$ is odd.

(\textbf{IA}: $m=0$, $\mathbf{G}\equiv1(2)$ and $\ell=\mathbf{G}$\textbf{),}

(\textbf{IB}: $m=0$, $\mathbf{G}\equiv0(2)$ and $\ell=\mathbf{G}+1$\textbf{),}

(\textbf{IC}: $m=1$, $\mathbf{G}\equiv0(2)$ and $\ell=\mathbf{G}$\textbf{),}

(\textbf{ID}: $m=1$, $\mathbf{G}\equiv1(2)$ and $\ell=\mathbf{G}+1$\textbf{).}

In the cases (\textbf{IA, IB, IC, ID)} we set
\[%
\begin{array}
[c]{ll}%
\Psi_{x}^{\mathbf{E}}=(I_{-1}^{n}\mathcal{T})\circ\mathbf{i}^{j}(y(x)) &
\text{for }x\in\mathbf{E}_{j}\text{ (}0\leq j\leq\ell\text{).}%
\end{array}
\]

\textbf{Case II}: $n>p\geq2$ and $n$ is even.

\textbf{(IIA}: $m\geq0,$ $\mathbf{G}\equiv0(2)$ and $\ell=\mathbf{G}%
+m+1$\textbf{). }
\[
\Psi_{x}^{\mathbf{E}}=\left\{
\begin{array}
[c]{ll}%
(I_{-1}^{n}\mathcal{T})\circ\mathbf{i}^{j}(y(x)) & \text{for }x\in
\mathbf{E}_{j}\text{ (}0\leq j\leq\mathbf{G}+1\text{),}\\
(I_{-1}^{n}\mathcal{T})(y(x)) & \text{for }x\in\mathbf{E}_{j}\text{
(}\mathbf{G}+1<j\leq\mathbf{G}+m+1\text{).}%
\end{array}
\right.
\]

(\textbf{IIB}: $m\geq0$, $\mathbf{G}\equiv1(2)$ and $\ell=\mathbf{G}%
+m$\textbf{). }
\[
\Psi_{x}^{\mathbf{E}}=\left\{
\begin{array}
[c]{ll}%
(I_{-1}^{n}\mathcal{T})\circ\mathbf{i}^{j}(y(x)) & \text{for }x\in
\mathbf{E}_{j}\text{ (}0\leq j\leq\mathbf{G}\text{),}\\
(I_{-1}^{n}\mathcal{T})(y(x)) & \text{for }x\in\mathbf{E}_{j}\text{
(}\mathbf{G}<j\leq\mathbf{G}+m\text{).}%
\end{array}
\right.
\]

(\textbf{IIC}: $m<0$, $\mathbf{G}\equiv0(2)$ and $\ell=\mathbf{G}%
+|m|+1$\textbf{).}
\[
\Psi_{x}^{\mathbf{E}}=\left\{
\begin{array}
[c]{ll}%
(I_{-1}^{n}\mathcal{T})\circ\mathbf{i}^{j}(y(x)) & \text{for }x\in
\mathbf{E}_{j}\text{ (}0\leq j\leq\mathbf{G}+1\text{),}\\
(I_{-1}^{n}\mathcal{T})\circ\mathbf{i}^{1}(y(x)) & \text{for }x\in
\mathbf{E}_{j}\text{ (}\mathbf{G}+1<j\leq\mathbf{G}+|m|+1\text{).}%
\end{array}
\right.
\]

(\textbf{IID}: $m<0$, $\mathbf{G}\equiv1(2)$ and $\ell=\mathbf{G}%
+|m|$\textbf{).}
\[
\Psi_{x}^{\mathbf{E}}=\left\{
\begin{array}
[c]{ll}%
(I_{-1}^{n}\mathcal{T})\circ\mathbf{i}^{j}(y(x)) & \text{for }x\in
\mathbf{E}_{j}\text{ (}0\leq j\leq\mathbf{G}\text{),}\\
(I_{-1}^{n}\mathcal{T})\circ\mathbf{i}^{1}(y(x)) & \text{for }x\in
\mathbf{E}_{j}\text{ (}\mathbf{G}<j\leq\mathbf{G}+|m|\text{).}%
\end{array}
\right.
\]

\textbf{Case III}: $n=p\geq2$ and $n$ is an odd number except for $3$ and $7$.

\textbf{(IIIA}: $m=0$ and $\ell=1$)
\[%
\begin{array}
[c]{ll}%
\Psi_{x}^{\mathbf{E}}=(I_{-1}^{n}\mathcal{T)}(y(x)) & \text{for }%
x\in\mathbf{E}_{j}\text{ (}j=0\text{, }1\text{).}%
\end{array}
\]

\textbf{(IIIB}: $m=1$ and $\ell=0$)
\[%
\begin{array}
[c]{ll}%
\Psi_{x}^{\mathbf{E}}=(I_{-1}^{n}\mathcal{T)}(y(x)) & \text{for }%
x\in\mathbf{E}_{0}.
\end{array}
\]

\textbf{Case IV:} $n=p$, $n=3$, $7$, $m=0$ and $\ell=0$.
\[%
\begin{array}
[c]{ll}%
\Psi_{x}^{\mathbf{E}}=(I_{-1}^{n}\mathcal{T})(y(x)) & \text{for }%
x\in\mathbf{E}_{0}.
\end{array}
\]

\textbf{Case V}: $n=p\geq2$ and $n$ is even.

(\textbf{Case VA}: $m>0$ and $\ell=m-1$)%
\[%
\begin{array}
[c]{ll}%
\Psi_{x}^{\mathbf{E}}=(I_{-1}^{n}\mathcal{T)}(y(x)) & \text{for }%
x\in\mathbf{E}_{j}\text{ (}0\leq j\leq m-1\text{).}%
\end{array}
\]

(\textbf{Case VB}: $m=0$ and $\ell=1$)%
\[%
\begin{array}
[c]{ll}%
\Psi_{x}^{\mathbf{E}}=(I_{-1}^{n}\mathcal{T})\circ\mathbf{i}^{j}(y(x)) &
\text{for }x\in\mathbf{E}_{j}\text{ (}0\leq j\leq1\text{).}%
\end{array}
\]

(\textbf{Case VC}: $m<0$ and $\ell=|m|-1$)%
\[%
\begin{array}
[c]{ll}%
\Psi_{x}^{\mathbf{E}}=(I_{-1}^{n}\mathcal{T)}\circ\mathbf{i}^{1}(y(x)) &
\text{for }x\in\mathbf{E}_{j}\text{ (}0\leq j\leq|m|-1\text{).}%
\end{array}
\]

We can extend $\Psi$ and $\Psi^{\mathbf{E}}$\ to a homomorphism $\Psi
_{U}:\theta_{U}^{n}\rightarrow\theta_{U}^{n}$ so that $\Psi_{U}|_{\mathbf{E}%
_{j}}:\mathbf{E}_{j}\rightarrow{\mathrm{Hom}}(\mathbf{R}^{n},\mathbf{R}^{n})$
coincides with $\Psi^{\mathbf{E}}|_{\mathbf{E}_{j}}$. In fact, in
Int$\frak{D}$,\ connect the base point $(0,\ldots,0,1)$ of $\mathbf{E}_{0}$
with the base points $(0,\ldots,0,1)$ of the other $\mathbf{E}_{j}$'s
($j\neq0$) by paths, which are disjoint except at $(0,\ldots,0,1)$ of
$\mathbf{E}_{0}$. By the definition of $\Psi^{\mathbf{E}}$ we have that the
sum $\Psi^{\mathbf{E}}|_{\partial\mathbf{E}_{0}}+\cdots+\Psi^{\mathbf{E}%
}|_{\partial\mathbf{E}_{\ell}}$ and $\Psi^{\mathbf{E}}|_{\partial\frak{D}}$
are homotopic to each other when regarded as maps $S^{n-1}\rightarrow
GL^{+}(n)$. Then it is an elementary consequence of the homotopy theory that
$\Psi$ and $\Psi^{\mathbf{E}}$ are extended to $\Psi_{U}$ as maps
$U\rightarrow{\mathrm{Hom}}(\mathbf{R}^{n},\mathbf{R}^{n})$ so that $\Psi
_{U}|(U\setminus(\cup_{j=0}^{\ell}$Int$\mathbf{E}_{j}))$ is a map into
$GL^{+}(n)$. For $n\geq p\geq2$ let $\pi_{\mathbf{R}^{p}}:\mathbf{R}%
^{n}\rightarrow\mathbf{R}^{p}$ be the canonical projection onto the first $p$
components. Then we set%
\begin{equation}%
\begin{array}
[c]{ll}%
H(\mathbf{v}_{x})=(x,d_{x}g(\mathbf{v}_{x})) & \text{for }\mathbf{v}_{x}\in
T_{x}(N\setminus\text{Int}\frak{D}),\\
H(t_{N}^{-1}(x,\mathbf{v}))=t_{P}^{-1}(x,\pi_{\mathbf{R}^{p}}\circ\Psi
_{x}(\mathbf{v})) & \text{for }x\in U\text{, }\mathbf{v\in R}^{n}.
\end{array}
\end{equation}
This is well defined by the property (ii). In particular, we have
$H(t_{N}^{-1}(x,\mathbf{v}))=t_{P}^{-1}(x,\pi_{\mathbf{R}^{p}}\circ\Psi
_{x}^{\mathbf{E}}(\mathbf{v}))$ for $x\in\mathbf{E}_{j}$ $(0\leq j\leq\ell)$,
$\mathbf{v\in R}^{n}$. We note that if $n>p$, then $\pi_{\mathbf{R}^{p}%
}(I_{-1}^{n}\mathcal{T}(y)(\mathbf{v}))=\pi_{\mathbf{R}^{p}}(\mathcal{T}%
(y)(\mathbf{v}))$ and $\pi_{\mathbf{R}^{p}}(I_{-1}^{n}I_{-1}^{n}%
\mathcal{T}(y)I_{-1}^{n}(\mathbf{v}))=\pi_{\mathbf{R}^{p}}(\mathcal{T}%
(y)I_{-1}^{n}(\mathbf{v}))$ for $y\in\mathbf{R}^{n}$, $\mathbf{v\in R}^{n}$.

In $\mathbf{E}_{j}$, let $S_{\mathbf{E}_{j}}^{p-1}$ denote the $(p-1)$-sphere
consisting of all points \ $(y_{1},\ldots,y_{p},0,\ldots,0)$ with $2(y_{1}%
^{2}+\cdots+y_{p}^{2})=1$. By Lemma 3.3, we have that $H$ is of rank $p$
outside of $\cup_{j=0}^{\ell}S_{\mathbf{E}_{j}}^{p-1}$ and that $H$ is of rank
$p-1$ on $\cup_{j=0}^{\ell}S_{\mathbf{E}_{j}}^{p-1}$. Hence, $H$ defines a
section $s_{H}:N\rightarrow\Omega^{n-p+1}(N,P)$.

Now we are ready to define a section $s\in\Gamma^{tr}(N,P)$ satisfying the
assumption of Theorem 2.4. Under the identification (1.2), we first define
$s(x)=(x,f(x);H_{x},B_{x})$ for $x\in\mathbf{E}_{j}$ as follows.

Let $n>p\geq2$. Define the integer $I(j)$ to be $I(j)=j$ for $j\leq\mathbf{G}$
and $I(j)=0$ for $j>\mathbf{G}$. Then we set
\[
B_{x}=-(y_{1}(x)\Delta\lbrack4;I(j)],\ldots,y_{p}(x)\Delta\lbrack4;I(j)]),
\]
where if $x\in S_{\mathbf{E}_{j}}^{p-1}$, then we use the notation
$c(x)=(c_{1}(x),\ldots,c_{n}(x))$ in place of $y(x)=(y_{1}(x),\ldots
,y_{n}(x))$ and write $B_{x}=-(c_{1}(x)\Delta\lbrack4;I(j)],\ldots
,c_{p}(x)\Delta\lbrack4;I(j)]).$ Recall that
\begin{align*}
H_{x}  &  =\pi_{\mathbf{R}^{p}}\circ(I_{-1}^{n}\mathcal{T)}\circ\mathbf{i}%
^{i}(y(x))\\
&  =\left\{
\begin{array}
[c]{ll}%
\pi_{\mathbf{R}^{p}}\circ(I_{-1}^{n}\mathcal{T})(y(x)) & \text{if }i\text{ is
even},\\
\pi_{\mathbf{R}^{p}}\circ(\mathcal{T}I_{-1}^{n})(y(x)) & \text{if }i\text{ is
odd}.
\end{array}
\right.
\end{align*}
By Lemma 3.3, if $x\in S_{\mathbf{E}_{j}}^{p-1}$, then \textrm{Ker}$H_{x}$ is
generated by the vectors $\mathbf{c}(x)$, $\mathbf{e}_{p+1}$,...,
$\mathbf{e}_{n}$ and \textrm{Cok}$H_{x}$ is generated and oriented by
$-\overset{\bullet}{\mathbf{c}}(x)$. Furthermore, we have
\begin{equation}%
\begin{array}
[c]{ll}%
q(s(x))(\mathbf{c}(x),\mathbf{c}(x))=2(-\overset{\bullet}{\mathbf{c}}(x)), &
\\
q(s(x))(\mathbf{e}_{k},\mathbf{c}(x))=0 & \text{ for }p<k\leq n,\\
q(s(x))(\mathbf{e}_{k},\mathbf{e}_{\ell})=\left\{
\begin{array}
[c]{l}%
4\delta_{k\ell}(-\overset{\bullet}{\mathbf{c}}(x))\\
-4\delta_{k\ell}(-\overset{\bullet}{\mathbf{c}}(x))
\end{array}
\right.  &
\begin{array}
[c]{l}%
\text{for }p<k\leq n-I(j)\text{ and }p<\ell\leq n,\\
\text{for }n-I(j)<k\leq n\text{ and }p<\ell\leq n.
\end{array}
\end{array}
\end{equation}
Namely, the symmetric matrix associated to $q(s(x))$ has the index $I(j)$
under the basis $\mathbf{c}(x)$, $\mathbf{e}_{p+1}$,...,$\mathbf{e}_{n}$.

Let $n=p\geq2$. When $m\leq0$ and $n$ is even, we have to use $(I_{-1}%
^{n}\mathcal{T)}\circ\mathbf{i}^{1}(y(x))$ on $\mathbf{E}_{j}$ for some $j$.
By Lemma 3.3 (2), Ker($I_{-1}^{n}\mathcal{T}(c(x)))$ (resp. $(I_{-1}%
^{n}\mathcal{T)}\circ\mathbf{i}^{1}(c(x))$) is generated by $\mathbf{c}(x)$
(resp. $I_{-1}^{n}(\mathbf{c}(x))$), and Cok$(I_{-1}^{n}\mathcal{T}(c(x)))$
(resp. Cok$((I_{-1}^{n}\mathcal{T)}\circ\mathbf{i}^{1}(y(x)))$) is generated
and oriented by $-\mathbf{c}^{-}(x)$ (resp. $-\mathbf{c}(x)$). In the cases
III, IV, VA, and VB with $j=0$ we set%
\[
B_{x}=-(y_{1}(x)\Delta\lbrack4;0],\ldots,y_{n-1}(x)\Delta\lbrack
4;0],-y_{n}(x)\Delta\lbrack4;0])\text{ \ \ \ \ for }x\in\mathbf{E}_{j}.
\]
Then it follows that%
\begin{equation}%
\begin{array}
[c]{ll}%
q(s(x))(\mathbf{c}(x),\mathbf{c}(x))=2(-\mathbf{c}^{-}(x)) & \text{on
}S_{\mathbf{E}_{j}}^{p-1}.
\end{array}
\end{equation}
If either $m=0$ and $j=1$, or $m<0$, then we set%
\[
B_{x}=-(y_{1}(x)\Delta\lbrack4;0],\ldots,y_{n-1}(x)\Delta\lbrack
4;0],y_{n}(x)\Delta\lbrack4;0])\text{ \ \ \ \ for }x\in\mathbf{E}_{j}.
\]
Then it follows that%
\begin{equation}%
\begin{array}
[c]{ll}%
q(s(x))(I_{-1}^{n}\mathbf{c}(x),I_{-1}^{n}\mathbf{c}(x))=2(-\mathbf{c}(x)) &
\text{on }S_{\mathbf{E}_{j}}^{p-1}.
\end{array}
\end{equation}

For $x\in S_{\mathbf{E}_{j}}^{p-1}$, $q(s(x))$ is nonsingular and hence
$s(x)\in\Omega_{x,f(x)}^{n-p+1,0}(N,P)$. If $x\notin\cup_{j=0}^{\ell
}S_{\mathbf{E}_{j}}^{p-1}$, then $H_{x}$ is nonsingular and we can extend $H$
to a section $s\in\Gamma(N,P),$ since the fiber of $\pi_{1}^{2}|\Omega
^{n-p+1,0}(N,P)$ over $s(x)$ is contractible.

By construction, $S(s)^{\iota}$ is nonempty for each $0\leq\iota\leq
\lbrack(n-p+1)/2]$.

Next we see that $s$ lies in $\Gamma^{tr}(N,P)$. By definition, $d^{2}%
(s)_{\thicksim}\circ d(s)_{\thicksim}:TN|_{S(s)}\rightarrow{\mathrm{Hom}%
}(K(s),Q(s))$ on each $S_{\mathbf{E}_{j}}^{p-1}$ corresponds to the respective
second intrinsic derivative $\mathbf{d(T}I_{-\varepsilon}^{n}\mathcal{(}%
_{p}^{1}))$ or $\mathbf{d}(I_{-1}^{n}\mathcal{T)}\ $on each $S_{\mathbf{E}%
_{j}}^{p-1}$\ in Lemma 3.3 (we note that $\mathbf{i}^{1}$ is a
diffeomorphism). Hence, $s$ is transverse to $\Sigma^{n-p+1}(N,P)$ on each
$S_{\mathbf{E}_{j}}^{p-1}$. This implies the assertion.

We now consider the primary differences. Let $\mathcal{N}(s)$\ be the
orthogonal normal bundle of $\cup_{j=0}^{\ell}S_{\mathbf{E}_{j}}^{p-1}$. We
show that $d(i_{\mathcal{N}(s)}|_{S_{\mathbf{E}_{j}}^{p-1}},$ $\Phi
(s)|_{S_{\mathbf{E}_{j}}^{p-1}})=0$ for all $\mathbf{E}_{j}$ except for
\textbf{CaseVB}, $j=1$ and \textbf{CaseVC}, $0\leq j\leq|m|-1$. By Lemma 3.3
and (3.2) we have $\mathcal{N}(s)_{c(x)}=K(s)_{c(x)}$, which are generated by
$\mathbf{c}(x),\mathbf{e}_{p+1},\ldots,\mathbf{e}_{n}$ for $n>p$ and by
$\mathbf{c(}x)$\ for $n=p$.\ It follows from Lemma 3.3, (3.3), and (3.4) that
$Q(s)_{c(x)}$ is generated and oriented by $-\overset{\bullet}{\mathbf{c}}(x)$
for $n>p$ and by $-\mathbf{c}^{-}(x)$ as described above for $n=p$.
Furthermore, by the above fact about $d^{2}(s)_{\thicksim}\circ
d(s)_{\thicksim},$ we have, in the case $n>p\geq2,$%
\[%
\begin{array}
[c]{l}%
d_{c(x)}H(\Sigma_{i=1}^{p}c_{i}(x)\partial/\partial y_{i})(\mathbf{c(}%
x))=2(-\overset{\bullet}{\mathbf{c}}(x)),\\
d_{c(x)}H(\Sigma_{i=1}^{p}c_{i}(x)\partial/\partial y_{i})(\mathbf{e}%
_{k})=\mathbf{0}\text{ \ \ \ \ \ \ \ \ \ \ \ \ \ \ \ \ for }p<k\leq n,\\
d_{c(x)}H(\partial/\partial y_{k})(\mathbf{e}_{\ell})=\left\{
\begin{array}
[c]{ll}%
2\delta_{k\ell}(-\overset{\bullet}{\mathbf{c}}(x)) & \text{for }p<k\leq
n-I(j)\text{ and }p<\ell\leq n,\\
-2\delta_{k\ell}(-\overset{\bullet}{\mathbf{c}}(x)) & \text{for }n-I(j)<k\leq
n\text{ and }p<\ell\leq n.
\end{array}
\right.
\end{array}
\]
Hence, we have that $i_{\mathcal{N}(s)}|_{S_{\mathbf{E}_{j}}^{p-1}}$ and
$\Phi(s)|_{S_{\mathbf{E}_{j}}^{p-1}}$ are homotopic as monomorphisms
$\mathcal{N}(s)_{c(x)}\rightarrow T_{c(x)}N$ by (3.3). The case $n=p$ is similar.

Finally we consider \textbf{CaseVB}, $j=1$ and \textbf{CaseVC}, $0\leq
j\leq|m|-1$. The kernel of $(I_{-1}^{n}\mathcal{T)}\circ\mathbf{i}^{1}(y(x))$
is generated by $I_{-1}^{n}(\mathbf{c}(x))$ and the map $\partial
\mathbf{E}_{j}\rightarrow S^{n-1}$ defined by $x\longmapsto\sqrt{2}I_{-1}%
^{n}(\mathbf{c}(x))$ is of degree $-1$. Since the map $\partial\mathbf{E}%
_{j}\rightarrow S^{n-1}$ defined by $x\longmapsto\sqrt{2}(\mathbf{c}(x))$ is
of degree $1$, we have that $d(i_{\mathcal{N}(s)}|_{S_{\mathbf{E}_{j}}^{p-1}%
},$ $\Phi(s)|_{S_{\mathbf{E}_{j}}^{p-1}})=1-(-1)=2$. We have a way to remove
this obstruction by using Remark 4.5 and Lemma 5.6. Let $g_{n}:\mathbf{R}%
^{n}\rightarrow\mathbf{R}^{n}$ be the fold map defined by $g_{n}(x_{1}%
,\ldots,x_{n})=(x_{1},\ldots,x_{n-1},x_{n}^{2})$. By Remark 4.5, we can choose
two distinct points $c_{j}^{1}(x)$, $c_{j}^{2}(x)$ of $S_{\mathbf{E}_{j}%
}^{p-1}$ for all such $j$'s and respective small ball neighborhoods
$U(c_{j}^{1}(x))$, $V(c_{j}^{1}(x))$ and $U(c_{j}^{2}(x))$, $V(c_{j}^{2}(x))$
with $V(c_{j}^{1}(x))\subset$Int$U(c_{j}^{1}(x))$ and $V(c_{j}^{2}(x))\subset
$Int$U(c_{j}^{2}(x)),$ which do not intersect with the other $\mathbf{E}_{k}$.
Furthermore, there exists a homotopy $s_{\lambda}\in\Gamma^{tr}(N,P)$ such
that for all $j$ in \textbf{CaseVB}, $j=1$ and \textbf{CaseVC}, $0\leq
j\leq|m|-1$,

(1) $s_{0}=s$ and $S(s_{\lambda})=S(s),$

(2) $s_{\lambda}|(N\setminus(U(c_{j}^{1}(x))\cup U(c_{j}^{2}%
(x))))=s|(N\setminus(U(c_{j}^{1}(x))\cup U(c_{j}^{2}(x)))),$

(3) $s_{1}|V(c_{j}^{1}(x))=j^{2}g_{n}|V(c_{j}^{1}(x))$ and $s_{1}|V(c_{j}%
^{2}(x))=j^{2}g_{n}|V(c_{j}^{2}(x))$ hold under suitable coordinates on
$U(c_{j}^{1}(x))$, $U(c_{j}^{2}(x))$, $\pi_{P}\circ s_{1}(U(c_{j}^{1}(x)))$
and $\pi_{P}\circ s_{1}(U(c_{j}^{2}(x)))$, where the $n$-disk $D_{10}^{n}$
with radius $10$ in $\mathbf{R}^{n}$ corresponds to those points inside
Int$V(c_{j}^{1}(x))$ and Int$V(c_{j}^{2}(x))$.

We replace $s_{1}=j^{2}g_{n}$ with the section $\varpi^{+}$ given in Lemma
$5.6$ on $V(c_{j}^{1}(x))$ and $V(c_{j}^{2}(x))$ for each $j$. We denote this
new section by $\frak{s}$. Then by the additive property of the primary
difference we have $d(i_{\mathcal{N}(\frak{s})}|_{S_{\mathbf{E}_{j}}^{p-1}},$
$\Phi(\frak{s})|_{S_{\mathbf{E}_{j}}^{p-1}})=0$.

This implies that $s$ for the cases except for \textbf{CaseVB}, $j=1$ and
\textbf{CaseVC}, $0\leq j\leq|m|-1$ and $\frak{s}$ satisfy the assumption of
Theorem 2.4. Therefore, the assertion of Theorem 0.3 follows from Theorem 2.4.
This completes the proof.
\end{proof}

We obtain another version of a theorem due to \`{E}lia\v{s}berg[E2, 5.1D and 5.4].

\begin{corollary}
Let $n\geq p\geq2$. Any continuous map $f:S^{n}\rightarrow S^{p}$ is homotopic
to a fold-map which folds exactly on the boundaries of a finite number of
\ disjointly embedded disks of dimension $p$ within an embedded disk of
dimension $p$ in $S^{n}$.
\end{corollary}

After this paper was submitted, the author was informed of the result [Sp,
Theorem 5.1], which is related to Theorem 0.3. The approach of the arguments
to his result seems quite different from ours.

\section{Proof of Theorem 2.4}

Let $\delta$ be a small positive smooth function on a manifold $X$ equipped
with a Riemannian metric and let $E$ be a subbundle of $TX$. In this paper
$D_{\delta}(E)$ always denotes the associated disk bundle of $E$ with radius
$\delta$ such that $\exp_{X,x}|D_{\delta}(E)_{x}$ is an embedding for any
$x\in X$.

The method of the proof of Theorem 2.4 is technically the same as the proof of
[An3, Proposition 4.6] for the case $n=p$. Nevertheless, it will be necessary
in understanding Theorem 0.5 to write it noticing the arguments and tools for
the case $n>p$.

Let $s$ be a section of {$\Gamma^{tr}(N,P)$ and} let $M$ be any one of
$M(s)_{j}^{\iota}$'s. If $d(i_{\mathcal{N}(s)}|_{_{{M}}},\Phi(s)|_{_{{M}}}%
)=0$, then there exists a homotopy $\varphi^{M}(s)_{\lambda}:\mathcal{N}%
(s)|_{M}\rightarrow TN|_{M}$ relative to $M\setminus$Int$N_{1}$ in
\textrm{Mono}$(\mathcal{N}(s)|_{M},TN|_{M})$ such that $\varphi^{M}%
(s)_{0}=i_{\mathcal{N}(s)}|_{M}$ and $\varphi^{M}(s)_{1}$$=\Phi(s)|_{M}$ by
the definition of the primary difference. Let $\mathrm{Iso}(TN|_{M},TN|_{M}%
)$\ denote the space consisting of all isomorphisms of $T_{c}N$, $c\in M$.\ We
define the restriction map
\[
r_{M}:\mathrm{Iso}(TN|_{M},TN|_{M})\rightarrow\mathrm{Mono}(\mathcal{N}%
(s)|_{M},TN|_{M})
\]
by $r_{M}(h)=h|\mathcal{N}(s)_{c}$, where $h:T_{c}N\rightarrow T_{c}N$. Then
$r_{M}$ induces a structure of a fiber bundle with fiber \textrm{Iso}%
$(\mathbf{R}^{p-1},\mathbf{R}^{p-1})\times{\mathrm{Hom}}(\mathbf{R}%
^{p-1},\mathbf{R}^{n-p+1})$. By applying the covering homotopy property of the
fiber bundle $r_{M}$ to the sections $id_{TN|_{M}}$ and the homotopy
$\varphi^{M}(s)_{\lambda},$ we obtain a homotopy $\Phi^{M}(s)_{\lambda}:$
$TN|_{M}\rightarrow TN|_{M}$ such that $\Phi^{M}(s)_{0}=id_{TN|_{M}}$,
$\Phi^{M}(s)_{\lambda}|_{c}=id_{T_{c}N}$ for any $c\in M\setminus$Int$N_{1}$
and $r_{M}\circ\Phi^{M}(s)_{\lambda}=\varphi^{M}(s)_{\lambda}$. We define
$\Phi(s)_{\lambda}:$ $TN|_{S(s)}\rightarrow TN|_{S(s)}$ by $\Phi(s)_{\lambda
}|_{M}=\Phi^{M}(s)_{\lambda}$ and $\Phi(s)_{\lambda}|_{c}=id_{T_{c}N}$ for any
$c\in U(C)\cap S(s)$. In the proof of the following lemma, $\Phi(s)_{\lambda
}|_{c}$ ($c\in S(s)$) is regarded as a linear isomorphism of $T_{c}N$.

Let $r_{0}$ be a small positive real number with $r_{0}<1/10$.

\begin{lemma}
Let $s\in\Gamma^{tr}(N,P)$ be a section satisfying the hypotheses of Theorem
$2.4$. Then there exists a homotopy $s_{\lambda}$ relative to $N\setminus
\mathrm{Int}N_{2-4r_{0}}$ in $\Gamma^{tr}(N,P)$ with $s_{0}=s$ satisfying

$(4.1.1)$ for any $\lambda$, $S(s_{\lambda})=S(s)$,

$(4.1.2)$ for any point $c\in S(s_{1})$, we have $K(s_{1})_{c}=\mathcal{N}%
(s)_{c}$ and $d^{2}(s_{1})_{\sim}\circ{d(s}_{1}{)}_{\sim}|K(s_{1})=d^{2}%
s_{1}|K(s_{1})$.
\end{lemma}

\begin{proof}
Recall the exponential map $\exp_{N,x}:T_{x}N\rightarrow N$ defined near
$\mathbf{0}\in T_{x}N$. We write an element of $\mathcal{N}(s)_{c\text{ }}$as
$(c,\mathbf{v})$. There exists a small positive number $\delta$ such that the
map
\[
e:D_{\delta}(\mathcal{N}(s))|_{S(s)\cap N_{2}}\rightarrow N
\]
defined by $e(c,\mathbf{v})=\exp_{N,c}(c,\mathbf{v})$ is an embedding, where
$c\in S(s)\cap N_{2}$ and $(c,\mathbf{v})\in D_{\delta}(\mathcal{N}(s)_{c})$
(note that $e|(S(s)\cap N_{2})$ is the inclusion). Let $\psi:[0,\infty
)\rightarrow\lbrack0,1]$ be a decreasing smooth function such that $\psi(t)=1$
if $t\leq\delta/10$ and $\psi(t)=0$ if $t\geq\delta$.

If we represent $s(x)\in\Omega^{n-p+1,0}(N,P)$ by a jet $j_{x}^{2}\sigma_{x}$
for a germ $\sigma_{x}:(N,x)\rightarrow(P,\pi_{P}^{2}\circ s(x))$, then we
define the homotopy $s_{\lambda}$ of $\Gamma^{tr}(N,P)$ using $\Phi
(s)_{\lambda}$ by
\begin{equation}
\left\{
\begin{array}
[c]{ll}%
\begin{array}
[c]{l}%
s_{\lambda}(e(c,\mathbf{v}))\\
=j_{e(c,\mathbf{v})}^{2}(\sigma_{e(c,\mathbf{v})}\circ\exp_{N,c}%
\circ\widetilde{\Phi(s)}_{\psi(\Vert\mathbf{v}\Vert)\lambda}|_{c}\circ
\exp_{N,c}^{-1})
\end{array}
& \text{\textrm{if} $c\in S(s)\cap N_{2}$ \textrm{and} }\Vert\text{$\mathbf{v}%
\Vert\leq\delta,$}\\%
\begin{array}
[c]{l}%
s_{\lambda}(x)=s(x)
\end{array}
& \text{\textrm{if} $x\notin\mathrm{Im}(e).$}%
\end{array}
\right.
\end{equation}
Here, $\widetilde{\Phi(s)}_{\psi(\Vert\mathbf{v}\Vert)\lambda}|_{c}$\ refers
to $\ell(\mathbf{v})\circ(\Phi(s)_{\psi(\Vert\mathbf{v}\Vert)\lambda}%
|_{c})\circ\ell(-\mathbf{v})$, where $\ell(\ast)$ is defined in \S 1. If
$\Vert\mathbf{v}\Vert\geq\delta$, then $\Phi(s)_{\psi(\Vert\mathbf{v}%
\Vert)\lambda}|_{c}=\Phi(s)_{0}|_{c}$, and if $c\in S(s)\setminus
$Int$N_{2-4r_{0}}$, then $\Phi(s)_{\lambda}|_{c}=\Phi(s)_{0}|_{c}$. In these
cases we have
\begin{align*}
s_{\lambda}(e(c,\mathbf{v})) &  =j_{e(c,\mathbf{v})}^{2}(\sigma
_{e(c,\mathbf{v})}\circ\exp_{N,c}\circ\widetilde{\Phi(s)}_{\psi(\Vert
\mathbf{v}\Vert)\lambda}|_{c}\circ\exp_{N,c}^{-1})\\
&  =j_{e(c,\mathbf{v})}^{2}(\sigma_{e(c,\mathbf{v})})\\
&  =s(e(c,\mathbf{v})).
\end{align*}
Hence, $s_{\lambda}$ is well defined. Furthermore, we have that

(1) $\pi_{P}^{2}\circ s_{\lambda}(x)=\pi_{P}^{2}\circ s(x)$,

(2) $S(s_{\lambda})=S(s)$,

(3) if $c\in S(s)$, then we have that $\Phi(s)_{1}(\mathcal{N}(s)_{c}%
)=K(s)_{c}$,

(4) $s_{\lambda}$ is transverse to $\Sigma^{n-p+1,0}(N,P)$.

\noindent The property (4.1.2) is satisfied for $s_{1}$ by (4.1) and (3).
\end{proof}

For a vector bundle $\mathcal{F}$ over $\Sigma$ and a map $j:S\rightarrow
\Sigma$, the induced bundle map $j^{\ast}(\mathcal{F})\rightarrow\mathcal{F}$
over $j$ is denoted by $j^{\mathcal{F}}$ in the following.

The following lemma is proved similarly as in [An3, Lemma 6.2] by using the
Hirsch Immersion Theorem [H1], and its proof is left to the reader.

\begin{lemma}
Let $s$ be a section of $\Gamma^{tr}(N,P)$ satisfying the property $(4.1.2)$
for $s$ $($in place of $s_{1}$$)$ of Lemma $4.1$. Then there exists a homotopy
$s_{\lambda}$ relative to $N\setminus\mathrm{Int}N_{2-4r_{0}}$ in $\Gamma
^{tr}(N,P)$ with $s_{0}=s$ such that

$(4.2.1)$ $S(s_{\lambda})=S(s)$ for any $\lambda$,

$(4.2.2)$ $\pi_{P}^{2}\circ s_{1}|S(s_{1})$ is an immersion into $P$ such that
$d(\pi_{P}^{2}\circ s_{1}|S(s_{1})):TS(s_{1})\rightarrow TP$ is equal to
$(\pi_{P}^{2}\circ s_{1})^{TP}\circ d^{1}s_{1}|TS(s_{1})$, where $(\pi_{P}%
^{2}\circ s_{1})^{TP}:(\pi_{P}^{2}\circ s_{1})^{\ast}TP\rightarrow TP$ is the
canonical induced bundle map.
\end{lemma}

Here we give two lemmas necessary for the proof of Theorem 2.4. Their proofs
are elementary and so are left to the reader.

\begin{lemma}
Let $\pi:E\rightarrow S$ be a smooth $(n-p+1)$-dimensional vector bundle with
metric over a $(p-1)$-dimensional manifold, where $S$ is identified with the
zero-section. Let $f_{i}:E\rightarrow P$ $(i=1,2)$ be fold-maps which fold
exactly on $S$ such that

$(\mathrm{i})$ $f_{1}|S=f_{2}|S$,

$(\mathrm{ii})$ $d_{c}f_{1}=d_{c}f_{2}$ and $d_{c}^{2}f_{1}=d_{c}^{2}f_{2}$
for any $c\in S,$

$(\mathrm{iii})$ $K(j^{2}f_{i})_{c}$ are tangent to $\pi^{-1}(c)$ for any
$c\in S$.

\noindent Let $\eta:S\rightarrow\lbrack0,1]$ be any smooth function. Then
there exists a positive function $\varepsilon:S\rightarrow\mathbf{R}$ such
that $\exp_{P,f_{1}(c)}((1-\eta(c))\exp_{P,f_{1}(c)}^{-1}(f_{1}(\mathbf{v}%
_{c}))+\eta(c)\exp_{P,f_{2}(c)}^{-1}(f_{2}(\mathbf{v}_{c})))$, which is
denoted simply by $((1-\eta)f_{1}+\eta f_{2})(\mathbf{v}_{c})$, is contained
in a convex neighborhood of $f_{1}(c)=f_{2}(c)$ in $P$ for any $c\in S$ and
$\mathbf{v}_{c}\in\pi^{-1}(c)$ with $\Vert\mathbf{v}_{c}\Vert\leq
\varepsilon(c),$ and that the map $(1-\eta)f_{1}+\eta f_{2}$ is a fold-map
which folds exactly on $S$ and satisfies $d_{c}((1-\eta)f_{1}+\eta
f_{2})=d_{c}f_{i}$ and $d_{c}^{2}((1-\eta)f_{1}+\eta f_{2})=d_{c}^{2}f_{i}$.
\end{lemma}

\begin{lemma}
Let $\pi:E\rightarrow S$ be a smooth $(n-p+1)$-dimensional vector bundle with
metric over a $(p-1)$-dimensional manifold, where $S$ is identified with the
zero-section, and let $(\Omega,\Sigma)$ be a pair of a smooth manifold and its
submanifold of codimension $n-p+1$. Let $\varepsilon:S\rightarrow\mathbf{R}$
be a positive smooth function. Let $h_{i}:D_{\varepsilon}(E)\rightarrow
(\Omega,\Sigma)$ $(i=0,1)$ be smooth maps such that $S=h_{0}^{-1}%
(\Sigma)=h_{1}^{-1}(\Sigma)$, $h_{0}|S=h_{1}|S$ and that $h_{i}$ are
transverse to $\Sigma$. Assume that for any $c\in S$, the isomorphisms
$T_{c}E/T_{c}S\rightarrow T_{h_{i}(c)}{\Omega}/T_{h_{i}(c)}{\Sigma}$ induced
from $d_{c}(h_{i})$ coincide with each other. Then for a sufficiently small
positive function $\varepsilon:S\rightarrow\mathbf{R}$, there exists a
homotopy $h_{\lambda}:(D_{\varepsilon}(E),S)\rightarrow(\Omega,\Sigma)$
between $h_{0}$ and $h_{1}$ such that

$(1)$ $h_{\lambda}|S=h_{0}|S$, $h_{\lambda}^{-1}(\Sigma) =h_{0}^{-1}(\Sigma)$
for any $\lambda$,

$(2)$ $h_{\lambda}$ is smooth and is transverse to $\Sigma$ for any $\lambda$.
\end{lemma}

\begin{proof}
[Proof of Theorem 2.4.]By Lemmas 4.1 and 4.2 we may assume that $s$ satisfies
(4.1.2) and (4.2.2), where $s_{1}$ is replaced by $s$. In particular, we have
$\mathcal{N}(s)=K(s)$. Since $s$ is an embedding, we can choose a Riemannian
metric on $\Omega^{n-p+1,0}(N,P)$ so that the induced metric by $s$ coincides
with the metric on $N$ near $S(s)$. We set $E(S(s))=\exp_{N}(D_{\delta\circ
s}(K(s)))$, where $\delta:\Sigma^{n-p+1,0}(N,P)\rightarrow\mathbf{R}$ is a
sufficiently small positive function such that $\delta\circ s|(S(s)\cap
N{_{2}})$ is constant. This becomes a tubular neighborhood of $S(s)$.

We first prove the following assertion:

(A) {There exists a homotopy $h_{\lambda}$ relative to $N\setminus
\mathrm{Int}N_{2-2r_{0}}$ in $\Gamma^{tr}(N,P)$ with $h_{0}=s$ satisfying the
following. }

{$(1)$ $h_{\lambda}\in\Gamma^{tr}(N,P)$ and $S(h_{\lambda})^{\iota
}=S(s)^{\iota}$ for any $\lambda$ and }$\iota${. }

{$(2)$ We have a fold-map $G:(N\setminus\mathrm{Int}N_{2-2r_{0}})\cup
E(S(s))\rightarrow P$ with $j^{2}G=h_{1}.$ }

Take any Riemannian metric on $P$. If we identify $Q(s)$ with the orthogonal
normal line bundle to the immersion $\pi_{P}^{2}\circ s|S(s):S(s)\rightarrow
P$, then $\exp_{P}\circ(\pi_{P}^{2}\circ s|S(s))^{TP}|D_{\gamma}(Q(s))$ is an
immersion for some positive function $\gamma$. In the proof we express a point
of $E(S(s))$ as $\mathbf{v}_{c}$, where $c\in S(s)$, $\mathbf{v}_{c}\in
K(s)_{c}$\ and $\Vert\mathbf{v}_{c}\Vert\leq\delta(s(c))$. In the proof we say
that a smooth homotopy
\[
k_{\lambda}:(E(S(s)),\partial E(S(s)))\rightarrow(\Omega^{n-p+1,0}%
(N,P),\Sigma^{n-p}(N,P))
\]
has the property (C) if it satisfies that for any $\lambda$

(C-1) $k_{\lambda}^{-1}(\Sigma^{n-p+1,0}(N,P))=S(s)$, and $k_{\lambda
}|S(s)=k_{0}|S(s)$,

(C-2) $k_{\lambda}$ is smooth and transverse to $\Sigma^{n-p+1,0}(N,P)$.

We have the quadratic form $q(s):S^{2}K(s)\rightarrow Q(s)$ defined in \S1
such that
\[
q(s)_{c}(\mathbf{v}_{c},\mathbf{w}_{c})=(d_{c}^{2}s(\mathbf{v}_{c}%
))(\mathbf{w}_{c})\text{ \ \ for any point }c\in S(s).
\]
If we choose $\delta$ sufficiently small compared with $\gamma$, then we can
define the fold-map $g_{0}:E(S(s))\rightarrow P$ by
\begin{equation}
g_{0}(\mathbf{v}_{c})=\exp_{P,\pi_{P}^{2}\circ s(c)}\circ(\pi_{P}^{2}\circ
s|S(s))^{TP}(q(s)_{c}(\mathbf{v}_{c},\mathbf{v}_{c})).
\end{equation}
Now we need to modify $g_{0}$ by using Lemma 4.3 so that $g_{0}$ is compatible
with $g$. Let $\eta:S(s)\rightarrow\mathbf{R}$ be a smooth function such that

(i) $0\leq\eta(c)\leq1$ for $c\in S(s)$,

(ii) $\eta(c)=0$ for $c\in S(s)\setminus\mbox{{\rm Int}}N_{2-3r_{0}}$,

(iii) $\eta(c)=1$ for $c\in S(s)\cap N_{2-4r_{0}}$.

\noindent Then consider the map $G:(N\setminus\mathrm{Int}N_{2-2r_{0}})\cup
E(S(s))\rightarrow P$ defined by
\[
\left\{
\begin{array}
[c]{lll}%
G(x)=g(x) & \text{if $x\in N\setminus\mathrm{Int}N_{2-2r_{0}}$}, & \\
G(\mathbf{v}_{c})=(1-\eta(c))g(\mathbf{v}_{c})+\eta(c)g_{0}(\mathbf{v}_{c}) &
\text{if }\mathbf{v}_{c}\text{$\in E(S(s))$}. &
\end{array}
\right.
\]
Without loss of generality we may assume that $G$ is well-defined on a
neighborhood of $(N\setminus\mathrm{Int}N_{2-2r_{0}})\cup E(S(s))$ by
replacing $\delta$ and $E(S(s))$ by smaller ones. It follows from Lemma 4.3
that $G$ is a fold-map defined on a neighborhood of $(N\setminus\mbox{{\rm
Int}}N_{2-2r_{0}})\cup E(S(s))$, that $G|E(S(s))$ folds exactly on $S(s)$, and
that $d_{c}(G)=d_{c}(g_{0})$ and $d_{c}^{2}(G)=d_{c}^{2}(g_{0})$ for all
\ $c\in S(s)$. Furthermore, we note that if $c\in S(s)\setminus N_{2-2r_{0}},$
then $G(\mathbf{v}_{c})=g(\mathbf{v}_{c})$.

Next we construct a homotopy $H_{\lambda}$ relative to $N\setminus
\mathrm{Int}N_{2-2r_{0}}$\ in $\Gamma^{tr}((N\setminus\mathrm{Int}N_{2-2r_{0}%
})\cup E(S(s)),P)$ satisfying the property (C) such that $H_{0}=s|(N\setminus
\mathrm{Int}N_{2-2r_{0}})\cup E(S(s))$ and $H_{1}=j^{2}G$ on $(N\setminus
\mathrm{Int}N_{2-2r_{0}})\cup E(S(s))$. This construction is quite similar to
that in [An3, Proof of Proposition 4.6], while we write it down for completeness.

Set $\exp_{\Omega}=\exp_{\Omega^{n-p+1,0}(N,P)}$ for simplicity. By applying
Lemma 4.4 to the sections $s$ and $\exp_{\Omega}\circ ds\circ(\exp
_{N}|D_{\delta}(K(s)))^{-1}$ on $E(S(s))$, we first obtain a homotopy
$h_{\lambda}^{\prime}\in\Gamma^{tr}(E(S(s)),P)$ with $h_{0}^{\prime}=s$ and
$h_{1}^{\prime}=\exp_{\Omega}\circ ds\circ(\exp_{N}|D_{\delta}(K(s)))^{-1}$ on
$E(S(s))$ satisfying the properties (1) and (2) of Lemma 4.4. Let $p_{\Omega
}:T\Omega^{n-p+1,0}(N,P)|_{\Sigma^{n-p+1,0}(N,P)}\rightarrow T\Omega
^{n-p+1,0}(N,P)|_{\Sigma^{n-p+1,0}(N,P)}/T\Sigma^{n-p+1,0}(N,P)$ be the
canonical projection. We have $K(s)\linebreak =\mathcal{N}(s)$ and
$d^{2}(s)_{\thicksim}\circ d(s)_{\thicksim}|K(s)=d^{2}s|K(s)$. This yields
that $p_{\Omega}\circ ds|K(s)$ and $p_{\Omega}\circ(s|S(s))^{\mathbf{K}%
}:K(s)\rightarrow T\Omega^{n-p+1,0}(N,P)|_{\Sigma^{n-p+1,0}(N,P)}%
/T\Sigma^{n-p+1,0}(N,P)$ coincide with each other. By Lemma 4.4 we can
construct a homotopy $h_{\lambda}^{\prime\prime}$ in $\Gamma^{tr}(E(S(s)),P)$
with the property (C) such that $h_{0}^{\prime\prime}=h_{1}^{\prime}$ and
$h_{1}^{\prime\prime}=\exp_{\Omega}\circ(s|S(s))^{\mathbf{K}}\circ(\exp
_{N}|D_{\delta}(K(s)))^{-1}$ on $E(S(s))$. By pasting $h_{\lambda}^{\prime}$
and $h_{\lambda}^{\prime\prime}$ we obtain a homotopy $h_{\lambda}^{1}%
\in\Gamma^{tr}(E(S(s)),P)$ satisfying the property (C) with $h_{0}^{1}=s$ and
$h_{1}^{1}=\exp_{\Omega}\circ(s|S(s))^{\mathbf{K}}\circ(\exp_{N}|D_{\delta
}(K(s)))^{-1}$ on $E(S(s))$.

We now recall the additive structure of $J^{2}(N,P)$ in (1.2). Then we have
the homotopy $j_{\lambda}:S(s)\rightarrow J^{2}(N,P)$ defined by
\[
j_{\lambda}(c)=(1-\lambda)s(c)+\lambda j^{2}G(c)\quad\text{covering }\pi
_{P}^{2}\circ s|S(s):S(s)\rightarrow P,
\]
where $\pi_{P}^{2}\circ s|S(s)$ is the immersion as in $(4.2.2)$. Since
$K(s)_{c}=K(j^{2}G)_{c}$ and $Q(s)_{c}=Q(j^{2}G)_{c}$ by the construction of
$\pi_{P}^{2}\circ s|S(s)$ and the fold-map $G$, it follows that for any $c\in
S(s)$ we have $K(j_{\lambda})_{c}=K(s)_{c}$ and $Q(j_{\lambda})_{c}=Q(s)_{c}$.
Hence, we have that
\[
d_{c}^{i}(j_{\lambda})=(1-\lambda)d_{c}^{i}(s)+\lambda d_{c}^{i}(j^{2}%
G)=d_{c}^{i}(s)=d_{c}^{i}(j^{2}G).
\]
This implies that $j_{\lambda}$ is a map of $S(s)$ into $\Sigma^{n-p+1,0}%
(N,P)$. Therefore, the homotopy of bundle maps $(j_{\lambda})^{\mathbf{K}%
}:K(s)\rightarrow(\mathbf{K}\subset)T\Omega^{n-p+1,0}(N,P)$ induces the
homotopy $h_{\lambda}^{2}$ satisfying the property (C) defined by
\[
h_{\lambda}^{2}=\exp_{\Omega}\circ(j_{\lambda})^{\mathbf{K}}\circ(\exp
_{N}|D_{\delta}(K(s)))^{-1}|E(S(s))
\]
such that $h_{0}^{2}=h_{1}^{1}=\exp_{\Omega}\circ(s|S(s))^{\mathbf{K}}%
\circ(\exp_{N}|D_{\delta}(K(s)))^{-1}$ and $h_{1}^{2}=\exp_{\Omega}\circ
(j^{2}G|S(s))^{\mathbf{K}}\circ(\exp_{N}|D_{\delta}(K(s)))^{-1}$ on $E(S(s))$.

By applying Lemma 4.4 to $j^{2}G|E(S(s))$ similarly as $h_{\lambda}^{1}$, we
have a homotopy $h_{\lambda}^{3}$ satisfying the property (C) such that
$h_{0}^{3}=h_{1}^{2}=\exp_{\Omega}\circ(j^{2}G|S(s))^{\mathbf{K}}\circ
(\exp_{N}|D_{\delta}(K(s)))^{-1}$ and $h_{1}^{3}=j^{2}G$ on $E(S(s))$.

Furthermore, we need to modify the homotopy obtained by pasting $h_{\lambda
}^{1}$, $h_{\lambda}^{2}$ and $h_{\lambda}^{3}$ on $(N_{2}\setminus
\mathrm{Int}N_{2-4r_{0}})\cap E(S(s))$ to a homotopy ${H}_{\lambda}^{\prime}$
in $\Gamma^{tr}(N_{2}\cap E(S(s)),P)$\ satisfying the property (C) such that
${H}_{1}^{\prime}$\ coincides with $j^{2}G$ on $N_{2}\cap E(S(s))$ and
${H}_{\lambda}^{\prime}|((N_{2}\setminus\mathrm{Int}N_{2-2r_{0}})\cap
E(S(s)))=s|((N_{2}\setminus\mathrm{Int}N_{2-2r_{0}})\cap E(S(s)))$. Then we
can extend ${H}_{\lambda}^{\prime}$ to the homotopy $H_{\lambda}$ in
$\Gamma^{tr}((N\setminus\mathrm{Int}N_{2-2r_{0}})\cup E(S(s)),P)$ satisfying
the property (C) so that $H_{\lambda}(x)=s(x)$ for $x\in N\setminus
\mathrm{Int}N_{2-2r_{0}}$.

By applying the homotopy extension property to $s$ and $H_{\lambda}$, we
obtain an extended homotopy with $h_{0}=s$,
\[
h_{\lambda}:(N,S(s))\rightarrow(\Omega^{n-p+1,0}(N,P),\Sigma^{n-p+1,0}%
(N,P))\text{.}%
\]

\noindent Hence, $h_{\lambda}$ is a required homotopy in $\Gamma^{tr}(N,P)$ in
the assertion (A).

We apply Theorem 2.2 for the section $\pi_{1}^{2}\circ h_{1}$ and $G$. Since
$J^{1}(N,P)$ is canonically identified with \textrm{Hom}$((\pi_{N})^{\ast
}(TN),(\pi_{P})^{\ast}(TP))$, we may regard $\pi_{1}^{2}\circ h_{1}$ as a
homomorphism in $\frak{m}(N,P;S(s)^{0},\ldots,S(s)^{[(n-p+1)/2]}%
,N\setminus\mathrm{Int}N_{2-r_{0}},G)$. By Theorem 2.2 we obtain a homotopy
$B_{\lambda}:TN\rightarrow TP$ relative to $N\setminus\mathrm{Int}N_{2-r_{0}}$
of $\frak{m}(N,P;S(s)^{0},\ldots,S(s)^{[(n-p+1)/2]},N\setminus\mathrm{Int}%
N_{2-r_{0}},G)\ $and\ a fold-map $f:N\rightarrow P$ with $S(f)=S(h_{1})$ such
that $B_{0}=\pi_{1}^{2}\circ h_{1}$ and $B_{1}=df$. The homomorphism
$B_{\lambda}:TN\rightarrow TP$ in $\frak{m}(N,P;S^{0},\ldots,S^{[(n-p+1)/2]}%
,N\setminus\mathrm{Int}N_{2-r_{0}},G)$ satisfies, by definition, that for any
point $c\in S(s)^{\iota}$ and any $\lambda$\ there exist a small neighborhood
$U_{c}^{\lambda}$ and a fold-map $f_{U_{c}^{\lambda}}:U_{c}^{\lambda
}\rightarrow P$\ such that $df_{U_{c}^{\lambda}}=B_{\lambda}|TU_{c}^{\lambda}$
and $S(j^{2}f_{U_{c}^{\lambda}})^{\iota}=S(s)^{\iota}\cap U_{c}^{\lambda}%
$.\ Then this homotopy is lifted to a homotopy $\widetilde{B_{\lambda}}$
relative to $N\setminus\mathrm{Int}N_{2-r_{0}}$ in $\Gamma(N,P)$ such that

(1) $\widetilde{B_{0}}=h_{1},$

(2) $\widetilde{B_{1}}=j^{2}f,$

(3) $\pi_{1}^{2}\circ\widetilde{B_{\lambda}}=B_{\lambda}$,

(4) for any point $c\in S(s)$, $d_{c}^{2}(\widetilde{B_{\lambda}}%
(c))=d_{c}^{2}f_{U_{c}^{\lambda}}$,

\noindent since any fiber of $\pi_{1}^{2}:\Omega^{n-p+1,0}(N,P)\setminus
\Sigma^{n-p+1,0}(N,P)\rightarrow J^{1}(N,P)\setminus\Sigma^{n-p+1}(N,P)$ is
contractible and so on. Finally we obtain the required homotopy $s_{\lambda}$
by pasting $h_{\lambda}$ and $\widetilde{B_{\lambda}}$.

By the construction, we have that $S(s_{\lambda})=S(s)$ and $Q(s_{\lambda
})=Q(s)$ for all $\lambda$. This proves Theorem 2.4 (3).
\end{proof}

\begin{remark}
Let $g^{\iota}:\mathbf{R}^{n}\rightarrow\mathbf{R}^{p}$ be the fold map
defined by $g^{\iota}(x_{1},\ldots,x_{n})=(x_{1},\ldots,x_{p-1},\linebreak
x_{p}^{2}+\cdots+x_{n-\iota}^{2}-x_{n-\iota+1}^{2}-\cdots-x_{n}^{2}).$ Let
$s\in\Gamma^{tr}(N,P)$ satisfy the same assumption of Proposition $2.5$. Then
by the argument in the above proof $($in particular, by $(4.2))$ we can prove
the following. Given a point $c$ of $M(s)_{j}^{\iota}$, there exist small ball
neighborhoods $U(c)$, $V(c)$ of $c$ with $V(c)\subset\mathrm{Int}U(c)$ and a
homotopy $s_{\lambda}\in\Gamma^{tr}(N,P)$ such that $U(c)$ does not intersect
all the other $M(s)_{j}^{\iota}$'s and that

$(1)$ $s_{0}=s$ and $S(s_{\lambda})=S(s),$

$(2)$ $s_{\lambda}|(N\setminus U(c))=s|(N\setminus U(c)),$

$(3)$ $s_{1}|V(c)=j^{2}g^{\iota}|V(c)$ under suitable coordinates of $U(c)$
and $\pi_{P}^{2}\circ s_{1}(U(c))$.
\end{remark}

\section{Proof of Proposition 2.5}

Let $X$ be an $i\times j$ matrix and $Y$ be a $k\times\ell$ matrix. Then
$X\dotplus Y$ denotes the $(i+k)\times(j+\ell)$ matrix, which is the direct
sum of $X$ and $Y$. Let $D_{r}^{k}$ be the $k$-disk of radius $r$ in
$\mathbf{R}^{k}$. For a point $x=(x_{1},\ldots,x_{n})\in\mathbf{R}^{n}$
$(n>p)$, set\ $\overset{\bullet}{x}$ $=(x_{1},\ldots,x_{p})$ and
$\overset{\bigstar}{x}=(x_{p+1},\ldots,x_{n})$.\ Let $\varkappa
:[0\mathbf{,\infty)}\rightarrow\lbrack0,1]$ be a decreasing smooth function
such that $\varkappa(t)=1$ for $0\leq t\leq1/10$, $0<\varkappa(t)<1$ for
$1/10<t<1$, and $\varkappa(t)=0$ for $t\geq1$. Using the identification
$J^{2}(n,p)\cong\mathrm{Hom}(\mathbf{R}^{n},\mathbf{R}^{p})\oplus
\mathrm{Hom}(S^{2}\mathbf{R}^{n},\mathbf{R}^{p})$ in (1.4), we express an
element of $J^{2}(n,p)$ as $(\mathbf{A}_{p\times n};\mathbf{B}_{n\times n}%
^{1},\ldots,\mathbf{B}_{n\times n}^{p})$, where $\mathbf{A}_{p\times n}$ is a
$p\times n$ matrix and $\mathbf{B}_{n\times n}^{j}$ are $n\times n$ symmetric
matrices. Similarly, $(A_{p\times p};B_{p\times p}^{1},\ldots,B_{p\times
p}^{p})$ expresses an element of $J^{2}(p,p)$. Set $A_{p\times p}^{\varkappa
}(x)=\varkappa(\Vert\overset{\bigstar}{x}\Vert)A_{p\times p}+(1-\varkappa
(\Vert\overset{\bigstar}{x}\Vert))((E_{p-1})\dotplus(2x_{p}))$ and
$\mathbf{B}_{n\times n}^{p}=B_{p\times p}^{p}\dotplus\Delta(\overbrace
{2,\ldots,2}^{n-p-\iota},\overbrace{-2,\ldots,-2}^{\iota})$ with $0\leq
\iota\leq\lbrack(n-p+1)/2]$. Set
\begin{align*}
\mathbf{A}_{p\times n}^{\varkappa,\iota}(x)  &  =\left(  A_{p\times
p}^{\varkappa}(x)%
\begin{array}
[c]{c}%
\mathbf{0}_{(p-1)\times(n-p)}\\
2x_{p+1},\ldots,2x_{n-\iota},-2x_{n-\iota+1},\ldots,-2x_{n}%
\end{array}
\right)  \in\mathrm{Hom}(\mathbf{R}^{n},\mathbf{R}^{p}),\\
\mathbf{B}_{n\times n}^{\varkappa,j}(x)  &  =\varkappa(\Vert\overset{\bigstar
}{x}\Vert)(B_{p\times p}^{j}\dotplus\mathbf{0}_{(n-p)\times(n-p)}%
)\text{\ \ \ \ \ \ \ \ \ \ \ \ \ \ \ \ \ \ \ \ \ \ \ \ \ \ for\ \ }1\leq j\leq
p-1,\\
\mathbf{B}_{n\times n}^{\varkappa,\iota,p}(x)  &  =\varkappa(\Vert
\overset{\bigstar}{x}\Vert)(\mathbf{B}_{n\times n}^{p})+(1-\varkappa
(\Vert\overset{\bigstar}{x}\Vert))\Delta(\overbrace{0,\ldots,0}^{p-1}%
\overbrace{2,\ldots,2}^{n-p+1-\iota},\overbrace{-2,\ldots,-2}^{\iota}).
\end{align*}
For $n>p\geq3$, we define the map $\mathcal{J}_{p,n}^{\iota}(x):J^{2}%
(p,p)\rightarrow J^{2}(n,p)$ by
\begin{equation}
\mathcal{J}_{p,n}^{\iota}(x)(A_{p\times p};B_{p\times p}^{1},\ldots,B_{p\times
p}^{p})=(\mathbf{A}_{p\times n}^{\varkappa,\iota}(x);\mathbf{B}_{n\times
n}^{\varkappa,1}(x),\ldots,\mathbf{B}_{n\times n}^{\varkappa,p-1}%
(x),\mathbf{B}_{n\times n}^{\varkappa,\iota,p}(x)).
\end{equation}

Let $\mu:\mathbf{R}^{p}\rightarrow\mathbf{R}$ be a smooth map such that
$\mu^{-1}(0)$ is connected and that $\mu(\overset{\bullet}{x})=2x_{p}$ outside
of $D_{4}^{p}$. Let $\mu$ have $0$ as a regular value for a while. Then
$\mu^{-1}(0)$ is a $(p-1)$-dimensional submanifold of $\mathbf{R}^{p}$ such
that $(\mathbf{R}^{p}\setminus D_{4}^{p})\cap\mu^{-1}(0)=\mathbf{R}%
^{p-1}\times0\setminus D_{4}^{p}$. Let $N(\mu^{-1}(0))$ be the orthogonal
normal line bundle of $\mu^{-1}(0)$ in $\mathbf{R}^{p}$, which is oriented by
$\mathrm{grad}\mu$, and let $\mathbf{e}(N(\mu^{-1}(0))_{c})$ be the unit
vector consistent with the given orientation of $N(\mu^{-1}(0))_{c}$. We often
identify $\mu^{-1}(0)$, $N(\mu^{-1}(0))$ and $\mathbf{e}(N(\mu^{-1}(0))_{c})$
with $\mu^{-1}(0)\times\mathbf{0}_{n-p}$, the normal bundle of $\mu
^{-1}(0)\times\mathbf{0}_{n-p}$ in $\mathbf{R}^{p}\times\mathbf{0}_{n-p}$ and
its unit vector consistent with the given orientation respectively. We orient
$\mu^{-1}(0)$\ so that the juxtaposition of an oriented basis of $T_{c}%
\mu^{-1}(0)$ and $\mathbf{e}(N(\mu^{-1}(0))_{c})$ is compatible with the
canonical orientation of $\mathbf{R}^{p}$. We have the map $\mathbf{e}%
(\mu^{-1}(0)):\mu^{-1}(0)\rightarrow S^{p-1}$ defined by $\mathbf{e(}\mu
^{-1}(0))(c)=\mathbf{e}(N(\mu^{-1}(0))_{c})$, and the induced homomorphism
\[
\mathbf{e}(\mu^{-1}(0))^{\ast}:H^{p-1}(S^{p-1},\mathbf{e}_{p};\mathbf{Z}%
)\rightarrow H^{p-1}(\mu^{-1}(0),\mu^{-1}(0)\setminus D_{4}^{p};\mathbf{Z}%
)\cong\mathbf{Z}\text{.}%
\]
Let $[\mu^{-1}(0)]^{c}\in H^{p-1}(\mu^{-1}(0),\mu^{-1}(0)\setminus D_{4}%
^{p};\mathbf{Z})$ and $[S^{p-1}]^{c}$\ $\in H^{p-1}(S^{p-1},\mathbf{e}%
_{p};\mathbf{Z})$ denote the dual classes of the fundamental classes
$[\mu^{-1}(0)]$ and $[S^{p-1}]$ respectively. The integer $\mathrm{deg}%
(\mathbf{e}(\mu^{-1}(0)))$ such that $\mathbf{e}(\mu^{-1}(0))^{\ast}%
([S^{p-1}]^{c})=\mathrm{deg}(\mathbf{e}(\mu^{-1}(0)))[\mu^{-1}(0)]^{c}$ is
called the degree of $\mathbf{e}(\mu^{-1}(0))$.

We give another interpretation of $\mathrm{deg}(\mathbf{e}(\mu^{-1}(0)))$.
There are two monomorphisms $i_{N(\mu^{-1}(0))}$ and $\mathbf{n}_{\mu^{-1}%
(0)}$ in \textrm{Mono}$(N(\mu^{-1}(0)),T\mathbf{R}^{p}|_{\mu^{-1}(0)})$, where
$i_{N(\mu^{-1}(0))}$ is the inclusion and $\mathbf{n}_{\mu^{-1}(0)}$ is
defined by $\mathbf{n}_{\mu^{-1}(0)}(\mathbf{e}(N(\mu^{-1}(0))_{c}%
))=(c,\mathbf{e}_{p})$ for $c\in\mu^{-1}(0)$. It is clear that if $c\in
\mu^{-1}(0)\setminus D_{4}^{p}$, then $i_{N(\mu^{-1}(0))}(\mathbf{e}%
(N(\mu^{-1}(0))_{c}))=\mathbf{n}_{\mu^{-1}(0)}(\mathbf{e}(N(\mu^{-1}%
(0))_{c}))=(c,\mathbf{e}_{p})$. We regard $i_{N(\mu^{-1}(0))}$ and
$\mathbf{n}_{\mu^{-1}(0)}$ as sections of the fiber bundle \textrm{Mono}%
$(N(\mu^{-1}(0)),T\mathbf{R}^{p}|_{\mu^{-1}(0)})$ over $\mu^{-1}(0)$. Then in
the obstruction theory ([Ste, \S37.5]), the primary difference \linebreak
$d(i_{N(\mu^{-1}(0))},\mathbf{n}_{\mu^{-1}(0)})$ in $H^{p-1}(\mu^{-1}%
(0),\mu^{-1}(0)\setminus D_{4}^{p};\pi_{p-1}(\mathrm{Mono}(\mathbf{R}%
,\mathbf{R}^{p})))$ is defined to be the unique obstruction so that
$i_{N(\mu^{-1}(0))}$ and $\mathbf{n}_{\mu^{-1}(0)}$ are homotopic as sections
of \linebreak ${\mathrm{Mono}}(N(\mu^{-1}(0)),T\mathbf{R}^{p}|_{\mu^{-1}(0)})$
relative to $\mu^{-1}(0)\setminus D_{4}^{p}$.

We prove the following lemma.

\begin{lemma}
Under the identification $H^{p-1}(\mu^{-1}(0),\mu^{-1}(0)\setminus D_{4}%
^{p};\pi_{p-1}(\mathrm{Mono}(\mathbf{R},\mathbf{R}^{p})))\cong H^{p-1}%
(\mu^{-1}(0),\mu^{-1}(0)\setminus D_{4}^{p};\mathbf{Z})\cong\mathbf{Z}$, the
primary difference $d(i_{N(\mu^{-1}(0))},\mathbf{n}_{\mu^{-1}(0)})$ coincides
with $\mathrm{deg}(\mathbf{e}(\mu^{-1}(0))).$
\end{lemma}

\begin{proof}
Consider a triangulation of $\mu^{-1}(0)$ such that $\mu^{-1}(0)\cap D_{4}%
^{p}$ is a subcomplex and take a base point of $\mu^{-1}(0)$ outside
$D_{4}^{p}$. Let $Z$ be a $(p-1)$-simplex. Then there exists an isomorphism
$\phi_{Z}:Z\times(\mathbf{0}_{p-1}\times\mathbf{R)}\rightarrow N(\mu
^{-1}(0))|_{Z}$ induced from a curve connecting $Z$ and the base point. Since
$N(\mu^{-1}(0))$ is a trivial line bundle, we may assume that $\phi_{Z}$ is
defined by $\phi_{Z}(c,\mathbf{e}_{p})=\mathbf{e}(N(\mu^{-1}(0))_{c})$ for any
$c\in Z$. Then $\{i_{N(\mu^{-1}(0))}\circ\phi_{Z}\}$ and $\{\mathbf{n}%
_{\mu^{-1}(0)}\circ\phi_{Z}\}$ induce the elements $[i_{N(\mu^{-1}(0))}]$ and
$[\mathbf{n}_{\mu^{-1}(0)}]$ of $H^{p-1}(\mu^{-1}(0),\mu^{-1}(0)\setminus
D_{4}^{p};\pi_{p-1}(\mathrm{Mono}(\mathbf{R},\mathbf{R}^{p})))$ respectively.
It is obvious that $[\mathbf{n}_{\mu^{-1}(0)}]$ is the null element and that
$[i_{N(\mu^{-1}(0))}]$ corresponds to $\mathrm{deg}(\mathbf{e}(\mu
^{-1}(0)))[\mu^{-1}(0)]^{c}$. Consider the following exact sequence of
cohomology groups with the coefficient $\pi_{p-1}(\mathrm{Mono}(\mathbf{R}%
,\mathbf{R}^{p}))$:
\begin{align*}
H^{p-1}((\mu^{-1}(0),\mu^{-1}(0)\setminus D_{4}^{p})\times\lbrack0,1])  &
\rightarrow H^{p-1}((\mu^{-1}(0),\mu^{-1}(0)\setminus D_{4}^{p})\times
\{0,1\})\\
\overset{\delta}{\rightarrow}H^{p}((\mu^{-1}(0),\mu^{-1}(0)\setminus D_{4}%
^{p})\times([0,1],\{0,1\}))  &  \cong H^{p-1}(\mu^{-1}(0),\mu^{-1}(0)\setminus
D_{4}^{p}).
\end{align*}
Here, $\delta([i_{N(\mu^{-1}(0))}]\oplus(-[\mathbf{n}_{\mu^{-1}(0)}]))$
corresponds to $d(i_{N(\mu^{-1}(0))},\mathbf{n}_{\mu^{-1}(0)})$. This proves
the lemma.
\end{proof}

We set
\begin{align*}
g_{p}(x_{1},\ldots,x_{p})  &  =(x_{1},\ldots,x_{p-1},x_{p}^{2}),\\
g^{\iota}(x_{1},\ldots,x_{n})  &  =(x_{1},\ldots,x_{p-1},x_{p}^{2}%
+\cdots+x_{n-\iota}^{2}-x_{n-\iota+1}^{2}-\cdots-x_{n}^{2}).
\end{align*}
In the following definitions (5.2) and (5.3), $\mu$ may not have $0$ as a
regular value. Under the identification (1.6), define the section $\omega
(\mu):\mathbf{R}^{p}\rightarrow\Omega^{1,0}(\mathbf{R}^{p},\mathbf{R}^{p})$
by
\begin{equation}
\omega(\mu)(\overset{\bullet}{x})=(\overset{\bullet}{x},g_{p}(\overset
{\bullet}{x});E_{p-1}\dot{+}(\mu(\overset{\bullet}{x}));\overbrace
{\mathbf{0}_{p\times p},\ldots,\mathbf{0}_{p\times p}}^{p-1},\Delta
(\overbrace{0,\ldots,0}^{p-1},2)).
\end{equation}
Using (1.6), define the section $\mathcal{J}^{\iota}(\omega(\mu))\in
\Gamma(\mathbf{R}^{n},\mathbf{R}^{p})$ by%
\begin{equation}%
\begin{array}
[c]{ll}%
\mathcal{J}^{\iota}(\omega(\mu))(x)=(x,{g}^{\iota}(x);\mathcal{J}_{p,n}%
^{\iota}(x)(\pi_{\Omega}\circ\omega(\mu)(\overset{\bullet}{x}))) & \text{for
}n>p>2,\\
\mathcal{J}^{\iota}(\omega(\mu))=\omega(\mu) & \text{for }n=p>2.
\end{array}
\end{equation}
If we set $\mu^{\varkappa}(x)=\varkappa(\Vert\overset{\bigstar}{x}\Vert
)\mu(\overset{\bullet}{x})+2(1-\varkappa(\Vert\overset{\bigstar}{x}%
\Vert))x_{p}$, then it follows that

(i) if either $\Vert\overset{\bigstar}{x}\Vert\geq1$ or $\Vert\overset
{\bullet}{x}\Vert\geq4$, then $\mu^{\varkappa}(x)=2x_{p},$

(ii) if $\Vert\overset{\bigstar}{x}\Vert\leq1/10$, then $\mu^{\varkappa
}(x)=\mu(\overset{\bullet}{x})$.

\noindent By (5.1), $\mathcal{J}_{p,n}^{\iota}(x)(\pi_{\Omega}\circ\omega
(\mu)(\overset{\bullet}{x}))$ is written as $(\mathbf{A}_{p\times
n}^{\varkappa,\iota}(x);\mathbf{B}^{1},\ldots,\mathbf{B}^{p})$, where
\begin{align}
\mathbf{A}_{p\times n}^{\varkappa,\iota}(x)  &  =\left(
\begin{array}
[c]{ll}%
E_{p-1} & \text{ \ \ \ \ \ \ \ \ \ \ \ \ \ \ \ \ \ \ }\mathbf{0}%
_{(p-1)\times(n-p+1)}\\
0,\ldots,0, & \mu^{\varkappa}(x),2x_{p+1},\ldots,2x_{n-\iota},-2x_{n-\iota
+1},\ldots,-2x_{n}%
\end{array}
\right)  \in\mathrm{Hom}(\mathbf{R}^{n},\mathbf{R}^{p}),\nonumber\\
\mathbf{B}^{j}  &  =\mathbf{0}_{n\times n}\text{
\ \ \ \ \ \ \ \ \ \ \ \ \ \ \ \ \ \ \ \ \ \ \ \ \ \ \ \ \ \ \ \ \ \ \ \ \ \ \ \ \ \ \ \ \ \ \ \ \ \ \ \ \ \ \ \ \ \ for
}1\leq j\leq p-1,\nonumber\\
\mathbf{B}^{p}  &  =\Delta(\overbrace{0,\ldots,0}^{p-1},\overbrace{2,\ldots
,2}^{n-p+1-\iota},\overbrace{-2,\ldots,-2}^{\iota}).
\end{align}
Therefore, we have

(1) if $\Vert x\Vert\geq10$, then $\mu^{\varkappa}(x)=2x_{p}$ and
$\mathcal{J}^{\iota}(\omega(\mu))(x)=j^{2}g^{\iota}(x),$

(2) $S(\mathcal{J}^{\iota}(\omega(\mu)))=\mu^{-1}(0)\times\mathbf{0}_{n-p}$,
by (ii) and the fact that if $\Vert\overset{\bigstar}{x}\Vert\neq0$, then
rank$\mathbf{A}_{p\times n}^{\varkappa,\iota}(x) \linebreak =p$.

On any point $c\in S(\mathcal{J}^{\iota}(\omega(\mu)))$, the $2$-jet
\[
\pi_{\Omega}\circ\mathcal{J}^{\iota}(\omega(\mu))(c)=(E_{p-1}\dotplus
(\mathbf{0}_{n-p+1});\mathbf{0}_{n\times n},\ldots,\mathbf{0}_{n\times
n},\mathbf{B}^{p})
\]
is represented by the germ ${g}^{\iota}:(\mathbf{R}^{n},0)\rightarrow
(\mathbf{R}^{p},0)$. Hence, $K(\mathcal{J}^{\iota}(\omega(\mu)))_{c}$ and
$Q(\mathcal{J}^{\iota}(\omega(\mu)))_{c}$ are generated and oriented by
$\mathbf{e}_{p},\ldots,\mathbf{e}_{n}$ and $\mathbf{e}_{p}$ respectively.
Therefore, we have a canonical isomorphism $\mathrm{Hom}(K(\mathcal{J}^{\iota
}(\omega(\mu)))_{c},Q(\mathcal{J}^{\iota}(\omega(\mu)))_{c})\cong
\mathbf{R}^{n-p+1}$.

We choose $\mathcal{N}(\omega(\mu))$ to be $N(\mu^{-1}(0))$, and then
$i_{N(\mu^{-1}(0))}=i_{\mathcal{N}(\omega(\mu))}$. Let%
\begin{align*}
m(p,n)  &  :\mathrm{Mono}(N(\mu^{-1}(0)),T\mathbf{R}^{p}|_{\mu^{-1}%
(0)})\rightarrow\\
&  \mathrm{Mono}(N(\mu^{-1}(0))\oplus(\mathbf{0}_{p}\times T\mathbf{R}%
^{n-p}\mathbf{)|}_{\mu^{-1}(0)\times\mathbf{0}_{n-p}},T\mathbf{R}^{n}%
|_{\mu^{-1}(0)\times\mathbf{0}_{n-p}})
\end{align*}
for $p>2$ be the fiber map defined by
\[
m(p,n)(h)(c,\mathbf{v}+a_{p+1}\mathbf{e}_{p+1}+\cdots+a_{n}\mathbf{e}%
_{n})=(c,h(\overset{\bullet}{c},\mathbf{v})+a_{p+1}\mathbf{e}_{p+1}%
+\cdots+a_{n}\mathbf{e}_{n})\text{,}%
\]
where $h\in\text{Mono}(N(\mu^{-1}(0))_{\overset{\bullet}{c}},T_{\overset
{\bullet}{c}}\mathbf{R}^{p})$, $\overset{\bullet}{c}\in\mu^{-1}(0)$ and
$\mathbf{v}\in N(\mu^{-1}(0))_{\overset{\bullet}{c}}$.

Set
\begin{align*}
\frak{B}^{1}(\pi_{p-1})  &  =\frak{B(}\pi_{p-1}(\mathrm{Mono}(N(\mu
^{-1}(0))_{\overset{\bullet}{c}},T_{\overset{\bullet}{c}}\mathbf{R}%
^{p})))\text{, for }\overset{\bullet}{c}\in\mu^{-1}(0),\\
\frak{B}^{2}(\pi_{p-1})  &  =\frak{B(}\pi_{p-1}(\mathrm{Mono}((N(\mu
^{-1}(0))\oplus(\mathbf{0}_{p}\times T\mathbf{R}^{n-p}\mathbf{))}_{c}%
,T_{c}\mathbf{R}^{n})))\text{, for }c\in\mu^{-1}(0)\times\mathbf{0}_{n-p}.
\end{align*}
Then $m(p,n)$ induces the fiberwise homomorphism $m(p,n)^{\prime}:\frak{B}%
^{1}(\pi_{p-1})\rightarrow\frak{B}^{2}(\pi_{p-1})$ of the trivial local
coefficients, and the epimorphism
\begin{align*}
m(p,n)_{\ast}^{\prime}  &  :H^{p-1}(\mu^{-1}(0),\mu^{-1}(0)\setminus D_{4}%
^{p};\frak{B}^{1}(\pi_{p-1}))\\
&  \rightarrow H^{p-1}(\mu^{-1}(0)\times\mathbf{0}_{n-p},(\mu^{-1}(0)\setminus
D_{4}^{p})\times\mathbf{0}_{n-p};\frak{B}^{2}(\pi_{p-1})).
\end{align*}

By the definition of the primary difference we have the following lemma.

\begin{lemma}
Let $n>p>2$. We have%
\[
m(p,n)_{\ast}^{\prime}(d(i_{N(\mu^{-1}(0))},\mathbf{n}_{\mu^{-1}%
(0)}))=d(i_{\mathcal{N}(\mathcal{J}^{\iota}(\omega(\mu)))},\Phi(\mathcal{J}%
^{\iota}(\omega(\mu)))).
\]
\end{lemma}

\begin{proof}
Recall the monomorphism $i_{\mathcal{N}(\mathcal{J}^{\iota}(\omega(\mu)))}$
defined in $\S2$. It is easy to see that \linebreak $m(p,n)(i_{N(\mu^{-1}%
(0))})=i_{\mathcal{N}(\mathcal{J}^{\iota}(\omega(\mu)))}.$ Next we show that
$m(p,n)(\mathbf{n}_{\mu^{-1}(0)})$ is homotopic to $\Phi(\mathcal{J}^{\iota
}(\omega(\mu)))$. In fact, $m(p,n)(\mathbf{n}_{\mu^{-1}(0)})(\mathbf{e}%
(N(\mu^{-1}(0))_{\overset{\bullet}{c}}))=\mathbf{e}_{p}$, $m(p,n)(\mathbf{n}%
_{\mu^{-1}(0)})(\mathbf{e}_{\ell})=\mathbf{e}_{\ell}$ for $\ell>p$ and
$\Phi(\mathcal{J}^{\iota}(\omega(\mu)))(\mathbf{e}_{\ell})=\mathbf{e}_{\ell}$
for $\ell>p$. Hence, it suffices to show \linebreak $\Phi(\mathcal{J}^{\iota
}(\omega(\mu)))(\mathbf{e}(N(\mu^{-1}(0))_{\overset{\bullet}{c}}))=(\Vert
$grad$\mu(\overset{\bullet}{c})\Vert/2)m(p,n)(\mathbf{n}_{\mu^{-1}%
(0)})(\mathbf{e}(N(\mu^{-1}(0))_{\overset{\bullet}{c}}))$.

In its proof we use the notation introduced in [B, \S1]. Let $Z_{j,\ell}$ be
the coordinate corresponding to the $(j,\ell)$ component of a $p\times n$
matrix. We have, around $c\in S(\mathcal{J}^{\iota}(\omega(\mu)))$,
\begin{align*}
&  \left(  d(\mathcal{J}^{\iota}(\omega(\mu)))\left(  \sum_{i=1}^{n}a_{i}%
\frac{\partial}{\partial x_{i}}\right)  \right)  (Z_{j,\ell})\\
&  =\left(  \sum_{i=1}^{n}a_{i}\frac{\partial}{\partial x_{i}}\right)
(Z_{j,\ell}\circ\mathcal{J}^{\iota}(\omega(\mu)))\\
&  =\left\{
\begin{array}
[c]{ll}%
(\sum_{i=1}^{n}a_{i}\frac{\partial}{\partial x_{i}})(\text{constant})=0 &
\text{for }0\leq j\leq p-1\text{ and }j=p\text{, }1\leq\ell\leq p-1,\\
(\sum_{i=1}^{n}a_{i}\frac{\partial}{\partial x_{i}})(\mu(\overset{\bullet}%
{x})) & \text{for }j=\ell=p,\\
(\sum_{i=1}^{n}a_{i}\frac{\partial}{\partial x_{i}})(2x_{\ell})=2a_{\ell} &
\text{for }j=p\text{ and }p<\ell\leq n-\iota,\\
(\sum_{i=1}^{n}a_{i}\frac{\partial}{\partial x_{i}})(-2x_{\ell})=-2a_{\ell} &
\text{for }j=p\text{ and }n-\iota<\ell\leq n.
\end{array}
\right.
\end{align*}
From the definition of the intrinsic derivative, it follows that%
\begin{equation}%
\begin{array}
[c]{l}%
d^{2}(\mathcal{J}^{\iota}(\omega(\mu)))_{\thicksim}\circ d(\mathcal{J}^{\iota
}(\omega(\mu)))_{\thicksim}(\text{grad}\mu(\overset{\bullet}{c}))(\mathbf{e}%
_{p})=\Vert\text{grad}\mu(\overset{\bullet}{c})\Vert^{2}\mathbf{e}_{p},\\
d^{2}(\mathcal{J}^{\iota}(\omega(\mu)))_{\thicksim}\circ d(\mathcal{J}^{\iota
}(\omega(\mu)))_{\thicksim}(\text{grad}\mu(\overset{\bullet}{c}))(\mathbf{e}%
_{\ell})=0\text{ \ for }p<\ell\leq n,\\
d^{2}(\mathcal{J}^{\iota}(\omega(\mu)))_{\thicksim}\circ d(\mathcal{J}^{\iota
}(\omega(\mu)))_{\thicksim}(\mathbf{e}_{k})(\mathbf{e}_{\ell})=\left\{
\begin{array}
[c]{ll}%
2\delta_{k\ell}\mathbf{e}_{p} & \text{for }p<k\leq n-\iota\text{ and }%
p<\ell\leq n,\\
-2\delta_{k\ell}\mathbf{e}_{p} & \text{for }n-\iota<k\leq n\text{ and }%
p<\ell\leq n.
\end{array}
\right.
\end{array}
\end{equation}
These formulas prove that $\mathcal{J}^{\iota}(\omega(\mu))$ is transverse to
$\Sigma^{n-p+1}(\mathbf{R}^{n},\mathbf{R}^{p})$ when $0$ is a regular value of
$\mu$. Since the symmetric matrix associated to $d^{2}(\mathcal{J}^{\iota
}(\omega(\mu)))_{c}:K(\mathcal{J}^{\iota}(\omega(\mu)))_{c}\rightarrow
\mathrm{Hom}(K(\mathcal{J}^{\iota}(\omega(\mu)))_{c},Q(\mathcal{J}^{\iota
}(\omega(\mu)))_{c})$ under the basis $\mathbf{e}_{p},\ldots,\mathbf{e}_{n}$
is equal to \linebreak $\Delta(\overbrace{2,\ldots,2}^{n-p+1-\iota}%
,\overbrace{-2,\ldots,-2}^{\iota})$ by (5.4), it follows from the definition
of $\Phi(\mathcal{J}^{\iota}(\omega(\mu)))$ in (2.3) that
\begin{align*}
&  \Phi(\mathcal{J}^{\iota}(\omega(\mu)))(\text{grad}\mu(\overset{\bullet}%
{c}))\\
&  =i_{K(\mathcal{J}^{\iota}(\omega(\mu)))}\circ(d^{2}(\mathcal{J}^{\iota
}(\omega(\mu)))|K(\mathcal{J}^{\iota}(\omega(\mu))))^{-1}\\
&  \text{
\ \ \ \ \ \ \ \ \ \ \ \ \ \ \ \ \ \ \ \ \ \ \ \ \ \ \ \ \ \ \ \ \ \ \ \ }\circ
d^{2}(\mathcal{J}^{\iota}(\omega(\mu)))_{\thicksim}\circ d(\mathcal{J}^{\iota
}(\omega(\mu)))_{\thicksim}(\text{grad}\mu(\overset{\bullet}{c}))\\
&  =i_{K(\mathcal{J}^{\iota}(\omega(\mu)))}((\Vert\text{grad}\mu
(\overset{\bullet}{c})\Vert^{2}/2)\mathbf{e}_{p})\\
&  =(\Vert\text{grad}\mu(\overset{\bullet}{c})\Vert^{2}/2)\mathbf{e}_{p}.
\end{align*}
This is what we want.
\end{proof}

We have proved the following lemma in [An3, Lemma 7.2].

\begin{lemma}
Let $p\geq3$. For $i=1,\ldots,p-1$, there exist functions $\mu_{\lambda}%
^{i}:\mathbf{R}^{p}\rightarrow\mathbf{\ R}$, ${\lambda}\in\mathbf{R,}$ which
are smooth with respect to the variables $x_{1},\ldots,x_{p}$ and $\lambda$
such that

{$(1)$ $\mu_{\lambda}^{i}(\overset{\bullet}{x})=2x_{p}$ if $\lambda\leq0$ or
$\Vert\overset{\bullet}{x}\Vert\geq4$, }

{$(2)$ $\mu_{\lambda}^{i}(\overset{\bullet}{x})=\mu_{1}^{i}(\overset{\bullet
}{x})$ if $\lambda\geq1$, }

{$(3)$ if $|\lambda-(1/2)|\geq1/2$, then $0$ is a regular value of
$\mu_{\lambda}^{i}$, }

{$(4)$ if $p\geq3$ and }{$1\leq i<p-1$ or if $p=3$ and $i=2$, then the
manifold $(\mu_{1}^{i})^{-1}(0)$ is a connected oriented manifold, }

{$(5)$ $\mu_{1}^{i}$ has a unique point $(0,\ldots,0,1)\in$}$\mathbf{R}^{p}${
such that $\mathbf{e}(S(\omega(\mu_{1}^{i})))(0,\ldots,0,1)=-\mathbf{e}_{p}$
and the degree of $\mathbf{e}(S(\omega(\mu_{1}^{i})))$ is equal to $(-1)^{i}$. }
\end{lemma}

\begin{remark}
We explain how $S(\omega(\mu_{1}^{i}))$ in Lemma 5.3 is constructed from
$\mathbf{R}^{p-1}\times0$. If $p\geq3$ and $1\leq i<p-1$ $($resp. $p=3$ and
$i=2)$, then $S(\omega(\mu_{1}^{i}))$ is constructed by doing the surgery by
an $i$-handle $($resp. three disjoint and non-linking $1$-handles$)$ attached
to $\mathbf{R}^{p-1}\times0$ within an embedded $p$-disk in {$\mathbf{R}^{p}$,
and hence }$S(\omega(\mu_{1}^{i}))$ is connected and oriented. When $p=3$ and
$i=2$, this process of the surgery is different from the usual method $($see
the details in $[$An$3$, Lemma $7.2])$.
\end{remark}

As for the case $p=2$, we have the following lemma.

Let $A^{+}$ (resp. $A^{+\prime}$) be $[-4,4]\times\lbrack-3,5]$ (resp.
$[-\pi,\pi]\times\lbrack-\sqrt{2}-1,4]$), and let $A^{-}$ (resp. $A^{-\prime}%
$) be $[-4,4]\times\lbrack-5,3]$ (resp. $[-\pi,\pi]\times\lbrack-4,\sqrt
{2}+1]$). Let $\sigma:$ $\mathbf{R}^{2}\rightarrow\mathbf{R}^{2}$ be the
fold-map defined by $\sigma(y_{1},y_{2})=(y_{1},y_{2}^{2})$.

\begin{figure}[ptbh]
\begin{center}
\setlength{\unitlength}{0.7mm} \begin{picture}(70,125)(-50,-40)

\put(-50,0){\vector(1,0){100}}
\put(0,-40){\vector(0,1){100}}
\put(-40,-30){\line(1,0){80}}
\put(-40,50){\line(1,0){80}}
\put(-31,-28){\line(1,0){62}}
\put(-31,40){\line(1,0){62}}
\put(-31,10){\line(1,0){62}}
\put(-31,0){\thicklines\line(1,0){62}}
\put(-40,-30){\line(0,1){80}}
\put(40,-30){\line(0,1){80}}
\put(-31,-28){\line(0,1){68}}
\put(31,-28){\line(0,1){68}}
\put(0,20){\circle{19.2}}

\put(0,20){\circle{14}}
\put(0,20){\circle*{2}}
\put(-25,20){\vector(1,0){18}}
\put(-30,22){\shortstack[r]{\scriptsize$C(0,2)$}}
\put(3,52){\shortstack[r]{\scriptsize5}}
\put(40,52){\shortstack[r]{\scriptsize$A^{+}$}}
\put(3,42){\shortstack[r]{\scriptsize4}}
\put(30,42){\shortstack[r]{\scriptsize${A^{+}}^{\prime}$}}
\put(3,32){\shortstack[r]{\scriptsize3}}
\put(3,20){\shortstack[r]{\scriptsize2}}
\put(3,5){\shortstack[r]{\scriptsize1}}
\put(3,-25){\shortstack[r]{\scriptsize$-\sqrt{2}-1$}}
\put(3,-35){\shortstack[r]{\scriptsize$-3$}}
\put(-48,2){\shortstack[r]{\scriptsize$-4$}}
\put(-39,2){\shortstack[r]{\scriptsize$-\pi$}}
\put(33,2){\shortstack[r]{\scriptsize$\pi$}}
\put(43,2){\shortstack[r]{\scriptsize4}}

\put(0,63){\shortstack[r]{\scriptsize$x_2$}}
\put(53,0){\shortstack[r]{\scriptsize$x_1$}}
\end{picture}
\end{center}
\end{figure}

\begin{lemma}
There exists a homotopy $\omega_{\lambda}^{\pm}$ relative to $\mathbf{R}%
^{2}\setminus D_{8}^{2}$ in $\Gamma(\mathbf{R}^{2},\mathbf{R}^{2})$ such that

$(1)$ $\omega_{0}^{\pm}=j^{2}\sigma$,

$(2)$ $\omega_{1}^{\pm}\in\Gamma^{tr}(\mathbf{R}^{2},\mathbf{R}^{2}),$

$(3)$ $S(\omega_{1}^{\pm})$ is the union of $\mathbf{R}\times0$ and the circle
$C(0,\pm2)$ centered at $(0,\pm2)$ with radius $1/\sqrt{2},$

$(4)$ $d(i_{N(\mathbf{R}\times0)},\Phi(\omega_{1}^{\pm})|_{\mathbf{R}\times
0})=\mp1$ and $d(i_{N(C(0,\pm2))},\Phi(\omega_{1}^{\pm})|_{C(0,\pm2)})=0$.
\end{lemma}

\begin{proof}
Let $GL^{+}(2)$ (resp. $GL^{-}(2)$) refer to the set of the regular $2\times2$
matrices with positive (resp. negative) determinants, which is provided with
the orientation induced from $S^{1}$ by the map $U\longmapsto U\mathbf{e}%
_{2}/{\Vert U\mathbf{e}_{2}\Vert}$.

Let $h^{\pm}:\mathbf{R}^{2}\rightarrow\mathbf{R}^{2}$ be the fold-map defined
by $h^{\pm}(y_{1},y_{2})=e^{((1/2)-y_{1}^{2}-y_{2}^{2})}(y_{1},-(\pm y_{2})).$
We have the Jacobian matrix%
\[
Jh^{+}(y_{1},y_{2})=e^{((1/2)-y_{1}^{2}-y_{2}^{2})}%
\begin{pmatrix}
1-2y_{1}^{2} & -2y_{1}y_{2}\\
2y_{1}y_{2} & 2y_{2}^{2}-1
\end{pmatrix}
=((1)\dotplus(-1))Jh^{-}(y_{1},y_{2}),
\]
whose determinant is equal to $e^{(1-2(y_{1}^{2}+y_{2}^{2}))}(2(y_{1}%
^{2}+y_{2}^{2})-1)$. Therefore, $h^{\pm}$ folds exactly on the circle
$S_{1/\sqrt{2}}^{1}$ with radius $1/\sqrt{2}$ and $h^{+}$ (resp. $h^{-}$)
preserves (resp. reverses) the orientation outside of $S_{1/\sqrt{2}}^{1}$ and
reverses (resp. preserves) the orientation inside of $S_{1/\sqrt{2}}^{1}$.
Since
\begin{align*}
Jh^{\pm}(\cos\theta,\sin\theta)  &  =e^{-1/2}%
\begin{pmatrix}
-\cos2\theta & -\sin2\theta\\
\pm\sin2\theta & \pm(-\cos2\theta)
\end{pmatrix}
,\\
Jh^{\pm}((1/\sqrt{2})(\cos\theta,\sin\theta))  &  =%
\begin{pmatrix}
\sin^{2}\theta & -\sin\theta\cos\theta\\
\pm\sin\theta\cos\theta & \pm(-\cos^{2}\theta)
\end{pmatrix}
,
\end{align*}
$Jh^{\pm}|_{S^{1}}:S^{1}\rightarrow GL^{\pm}(2)$ is of degree $-(\pm2)$. It
follows that $K(j^{2}h^{\pm})_{(1/\sqrt{2})(\cos\theta,\sin\theta)}$ is
generated by ${}^{t}(\cos\theta,\sin\theta)$ and that $Q(j^{2}h^{\pm
})_{(1/\sqrt{2})(\cos\theta,\sin\theta)}$ is generated and oriented by ${}%
^{t}(-\cos\theta,\pm\sin\theta)$. Indeed, since $d^{2}/dt^{2}(e^{1/2-t^{2}%
/2}(t/\sqrt{2})(\cos\theta,-(\pm\sin\theta)))=(1/\sqrt{2})(t^{3}%
-3t)e^{1/2-t^{2}/2}(\cos\theta,-(\pm\sin\theta))$, we have%
\begin{align}
\frac{d^{2}}{dt^{2}}(h^{\pm}((t/\sqrt{2})(\cos\theta,\sin\theta)))|_{t=1}  &
=\frac{d^{2}}{dt^{2}}(e^{1/2-t^{2}/2}(t/\sqrt{2})(\cos\theta,-(\pm\sin
\theta)))|_{t=1}\nonumber\\
&  =\sqrt{2}(-\cos\theta,\pm\sin\theta).
\end{align}

Let%
\begin{align*}
&  z^{+}(x_{1},x_{2})\\
&  =(((1-1/\sqrt{2})x_{2}+1/\sqrt{2})\cos(x_{1}+3\pi/2),((1-1/\sqrt{2}%
)x_{2}+1/\sqrt{2})\sin(x_{1}+3\pi/2)),\\
&  z^{-}(x_{1},x_{2})\\
&  =((-(1-1/\sqrt{2})x_{2}+1/\sqrt{2})\cos(-x_{1}+\pi/2),(-(1-1/\sqrt{2}%
)x_{2}+1/\sqrt{2})\sin(-x_{1}+\pi/2)).
\end{align*}
Then $z^{+}$ (resp. $z^{-}$) maps the segments $[-\pi,\pi]\times0$, $[-\pi
,\pi]\times1$ and $[-\pi,\pi]\times(-\sqrt{2}-1)$ (resp. $[-\pi,\pi]\times0$,
$[-\pi,\pi]\times(-1)$ and $[-\pi,\pi]\times(\sqrt{2}+1)$) to the circles
centered at $(0,0)$ with radii $1/\sqrt{2}$ and $1$, and to the point $(0,0)$ respectively.

Let $D(0,\pm2)$ be the $2$-disk centered at $(0,\pm2)$ with radius $1$. Let
$M(2)$ be the space consisting of all $2\times2$ matrices. We define maps
$\mathbf{a}^{\pm}:{A^{\pm}}^{\prime}\rightarrow M(2)$ so that $\mathbf{a}%
^{\pm}({A^{\pm}}^{\prime}\setminus(C(0,\pm2)\cup\lbrack-\pi,\pi]\times
0))\subset GL(2)$ by setting%
\[%
\begin{array}
[c]{ll}%
\mathbf{a}^{+}(x_{1},x_{2})=Jh^{+}(z^{+}(x_{1},x_{2})) & \text{for }%
(x_{1},x_{2})\in\lbrack-\pi,\pi]\times\lbrack-\sqrt{2}-1,1]\text{,}\\
\mathbf{a}^{+}(x_{1},x_{2})=Jh^{+}(x_{1},x_{2}-2) & \text{for }(x_{1}%
,x_{2})\in D(0,2),\\
\mathbf{a}^{+}(x_{1},x_{2})=e^{-1/2}E_{2} & \text{for}\left\{
\begin{array}
[c]{l}%
x_{1}=\pm\pi\text{ and }1\leq x_{2}\leq4\text{, or}\\
-\pi\leq x_{1}\leq\pi\text{ and }x_{2}=4\text{, or}\\
x_{1}=0\text{ and }3\leq x_{2}\leq4\text{,}%
\end{array}
\right. \\
\mathbf{a}^{-}(x_{1},x_{2})=Jh^{-}(z^{-}(x_{1},x_{2})) & \text{for }%
(x_{1},x_{2})\in\lbrack-\pi,\pi]\times\lbrack-1,\sqrt{2}+1]\text{,}\\
\mathbf{a}^{-}(x_{1},x_{2})=Jh^{-}(x_{1},x_{2}+2) & \text{for }(x_{1}%
,x_{2})\in D(0,-2),\\
\mathbf{a}^{-}(x_{1},x_{2})=e^{-1/2}((1)\dotplus(-1)) & \text{for}\left\{
\begin{array}
[c]{l}%
x_{1}=\pm\pi\text{ and }-4\leq x_{2}\leq-1\text{, or}\\
-\pi\leq x_{1}\leq\pi\text{ and }x_{2}=-4\text{, or}\\
x_{1}=0\text{ and }-4\leq x_{2}\leq-3\text{.}%
\end{array}
\right.
\end{array}
\]
We here give the properties for $\mathbf{a}^{\pm}$ :

(1) if $-\pi\leq x_{1}\leq\pi$, then%
\begin{align*}
\mathbf{a}^{+}(x_{1},1)  &  =Jh^{+}(z^{+}(x_{1},1))=Jh^{+}(\cos(x_{1}%
+3\pi/2),\sin(x_{1}+3\pi/2)),\\
\mathbf{a}^{-}(x_{1},-1)  &  =Jh^{-}(z^{-}(x_{1},-1))=Jh^{-}(\cos(-x_{1}%
+\pi/2),\sin(-x_{1}+\pi/2)),
\end{align*}

(2) if $-\sqrt{2}-1\leq x_{2}\leq1$ for $\mathbf{a}^{+}$ and $-1\leq x_{2}%
\leq\sqrt{2}+1$ for $\mathbf{a}^{-}$, then%
\begin{align*}
\mathbf{a}^{\pm}(\pi,x_{2})  &  =\mathbf{a}^{\pm}(-\pi,x_{2})=Jh^{\pm
}(0,(1-1/\sqrt{2})x_{2}\pm1/\sqrt{2})\\
&  =e^{a^{\pm}(x_{2})}((1)\dotplus(\pm2\{(1-1/\sqrt{2})x_{2}\pm1/\sqrt
{2}\}^{2}\mp1))
\end{align*}
with $a^{\pm}(x_{2})=1/2-\{(1-1/\sqrt{2})x_{2}\pm1/\sqrt{2}\}^{2}$,

(3) $\mathbf{a}^{\pm}(0,\pm1)=\mathbf{a}^{\pm}(\pi,\pm1)=\mathbf{a}^{\pm}%
(-\pi,\pm1)=\mathbf{a}^{\pm}(0,\pm3)=e^{-1/2}((1)\dotplus(\pm1))$,

(4) if $-\pi\leq x_{1}\leq\pi$, then%
\[
\mathbf{a}^{\pm}(x_{1},\mp(\sqrt{2}+1))=Jh^{\pm}(0,0)=e^{1/2}((1)\dotplus
(\mp1)).
\]
By noting these facts, we can extend this to the space $A^{+\prime}$ (resp.
$A^{-\prime}$) by using a suitable retraction of $[-\pi,\pi]\times
\lbrack1,4]\setminus$Int$D(0,2)$ onto $\partial D(0,2)\cup0\times\lbrack3,4]$
(resp. of $[-\pi,\pi]\times\lbrack-4,-1]\setminus$Int$D(0,-2)$ onto $\partial
D(0,-2)\cup0\times\lbrack-4,-3]$).

Furthermore, there exists an extended map $\widetilde{\mathbf{a}^{\pm}%
}:{A^{\pm}}\rightarrow M(2)$ of $\mathbf{a}^{\pm}$ as follows:

(i) if $(x_{1},x_{2})\in A^{\pm\prime}$, then $\widetilde{\mathbf{a}^{\pm}%
}(x_{1},x_{2})=\mathbf{a}^{\pm}(x_{1},x_{2}),$

(ii) if $(x_{1},x_{2})\notin A^{\pm\prime},$ then

(ii-1) $\widetilde{\mathbf{a}^{\pm}}(x_{1},x_{2})=((1)\dotplus(2x_{2}))$ on
$\partial A^{\pm}$,

(ii-2) $\widetilde{\mathbf{a}^{\pm}}(x_{1},x_{2})\in GL^{+}(2)$ if $x_{2}>0$,
and $\widetilde{\mathbf{a}^{\pm}}(x_{1},x_{2})\in GL^{-}(2)$ if $x_{2}<0$,

(ii-3) if $\pi\leq|x_{1}|\leq4,$ and $-\sqrt{2}-1<x_{2}<1$ for $\mathbf{a}%
^{+}$ and $-1<x_{2}<\sqrt{2}+1$ for $\mathbf{a}^{-}$, then
\begin{align*}
\widetilde{\mathbf{a}^{\pm}}(x_{1},x_{2})  &  =\varphi_{\lbrack\pi,4]}%
(|x_{1}|)\left(
\begin{array}
[c]{cc}%
1 & 0\\
0 & 2x_{2}%
\end{array}
\right) \\
&  +(1-\varphi_{\lbrack\pi,4]}(|x_{1}|))e^{a^{\pm}(x_{2})}\left(
\begin{array}
[c]{cc}%
1 & 0\\
0 & \pm2\{(1-1/\sqrt{2})x_{2}\pm1/\sqrt{2}\}^{2}\mp1
\end{array}
\right)  ,
\end{align*}
where $\varphi_{\lbrack\pi,4]}:[\pi,4]\rightarrow\lbrack0,1]$ is an increasing
smooth function such that $\varphi_{\lbrack\pi,4]}(t)=0$ if $t\leq\pi+1/10$,
and $\varphi_{\lbrack\pi,4]}(t)=1$ if $t\geq4-1/10$. In particular,
$\widetilde{\mathbf{a}^{\pm}}(x_{1},0)=(1)\dotplus(0)$ if $\pi\leq|x_{1}%
|\leq4$.

Note in (ii-3) that if we set $\phi^{\pm}(|x_{1}|,x_{2})=2x_{2}\varphi
_{\lbrack\pi,4]}(|x_{1}|)+(\pm2\{(1-1/\sqrt{2})x_{2}\pm1/\sqrt{2}\}^{2}%
\mp1)(1-\varphi_{\lbrack\pi,4]}(|x_{1}|))e^{a^{\pm}(x_{2})}$, then $\phi^{\pm
}(|x_{1}|,x_{2})$ and $x_{2}$ become positive, zero and negative at the same
time and $d/dx_{2}(\phi^{\pm}(|x_{1}|,x_{2}))|_{x_{2}=0}=2\varphi_{\lbrack
\pi,4]}(|x_{1}|)+2(\sqrt{2}-1)(1-\varphi_{\lbrack\pi,4]}(|x_{1}|))>0$.

We now define a section $\omega{^{\pm}}\in\Gamma^{tr}(\mathbf{R}%
^{2},\mathbf{R}^{2})$, which will be the required section $\omega_{1}^{\pm}$,
under the identification (1.6). By considering (5.6), we first set

(1) if $(x_{1},x_{2})\in{A^{\pm}}$, then $\pi_{1}^{2}\circ\omega{^{\pm}}%
(x_{1},x_{2})=((x_{1},x_{2}),\sigma(x_{1},x_{2});\widetilde{\mathbf{a}^{\pm}%
}(x_{1},x_{2}))$,

(2) if $(x_{1},x_{2})\in D(0,\pm2)$, then $\omega^{\pm}(x_{1},x_{2}%
)=((x_{1},x_{2}),\sigma(x_{1},x_{2});\pi_{\Omega}\circ j^{2}h^{\pm}%
(x_{1},x_{2}-(\pm2))$,

(3$^{+}$) if $(x_{1},x_{2})\in\lbrack-\pi,\pi]\times\lbrack-\sqrt{2}-1,1/2]$,
then
\begin{align*}
&  \omega^{+}(x_{1},x_{2})\\
&  =((x_{1},x_{2}),\sigma(x_{1},x_{2});\widetilde{\mathbf{a}^{+}}(x_{1}%
,x_{2});-\cos(x_{1}+3\pi/2)\Delta(2,2),\sin(x_{1}+3\pi/2)\Delta(2,2)),
\end{align*}

(3$^{-}$) if $(x_{1},x_{2})\in\lbrack-\pi,\pi]\times\lbrack-1/2,\sqrt{2}+1]$,
then
\begin{align*}
&  \omega^{-}(x_{1},x_{2})\\
&  =((x_{1},x_{2}),\sigma(x_{1},x_{2});\widetilde{\mathbf{a}^{-}}(x_{1}%
,x_{2});-\cos(-x_{1}+\pi/2)\Delta(2,2),-\sin(-x_{1}+\pi/2)\Delta(2,2)),
\end{align*}

(4) if $\pi\leq|x_{1}|\leq4$, then
\[
\omega{^{\pm}}(x_{1},0)=((x_{1},0),\sigma(x_{1},0);\widetilde{\mathbf{a}^{\pm
}}(x_{1},0);\Delta(0,0),\Delta(\rho(|x_{1}|),2)),
\]
where $\rho:[\pi,4]\rightarrow\lbrack0,2]$ is a decreasing smooth function
such that $\rho(t)=2$ for $t\leq\pi+1/10$, $0<\rho(t)<2$ for $\pi
+1/10<t<4-1/10$ and $\rho(t)=0$ for $4-1/10\leq t$,

(5) if $(x_{1},x_{2})\in\mathbf{R}^{2}\setminus{A^{\pm}}$, then $\omega^{\pm
}(x_{1},x_{2})=j^{2}\sigma(x_{1},x_{2})$,

\noindent and then extend this to a section $\omega{^{\pm}}:\mathbf{R}%
^{2}\rightarrow J^{2}(\mathbf{R}^{2},\mathbf{R}^{2})$ arbitrarily by using
$\pi_{1}^{2}$.

We have to show $\omega{^{\pm}}\in\Gamma^{tr}(\mathbf{R}^{2},\mathbf{R}^{2})$.
It is easy to see that $S(\omega^{\pm})$ consists of $\mathbf{R}\times0$ and
$C(0,\pm2)$. We first show $\omega{^{\pm}}\in\Gamma(\mathbf{R}^{2}%
,\mathbf{R}^{2})$.\ Let $(x_{1},0)\in\lbrack-\pi,\pi]\times0$. Since
$Q(\omega^{\pm})_{(x_{1},0)}$ is oriented so that $q(\omega^{\pm})$ is
positive definite, $K(\omega^{\pm})_{(x_{1},0)}$ is generated by
$\mathbf{v}^{\pm}\mathbf{(}x_{1})$ and $Q(\omega^{\pm})_{(x_{1},0)}$ is
generated and oriented by $\mathbf{w}^{\pm}\mathbf{(}x_{1})$, where
\[%
\begin{array}
[c]{ll}%
\mathbf{v}^{+}\mathbf{(}x_{1})={}^{t}(\cos(x_{1}+3\pi/2),\sin(x_{1}%
+3\pi/2)), & \mathbf{v}^{+}\mathbf{(}\pm\pi)=\text{ }^{t}(0,1),\\
\mathbf{w}^{+}\mathbf{(}x_{1})={}^{t}(-\cos(x_{1}+3\pi/2),\sin(x_{1}%
+3\pi/2)), & \mathbf{w}^{+}\mathbf{(}\pm\pi)=\text{ }^{t}(0,1),\\
\mathbf{v}^{-}\mathbf{(}x_{1})={}^{t}(-\cos(-x_{1}+\pi/2),-\sin(-x_{1}%
+\pi/2)), & \mathbf{v}^{-}\mathbf{(}\pm\pi)=\text{ }^{t}(0,1),\\
\mathbf{w}^{-}\mathbf{(}x_{1})={}^{t}(-\cos(-x_{1}+\pi/2),-\sin(-x_{1}%
+\pi/2)), & \mathbf{w}^{-}\mathbf{(}\pm\pi)=\text{ }^{t}(0,1)
\end{array}
\]
(see \S1). We have%
\begin{align}
&  (d_{(x_{1},0)}^{2}\omega^{+})(\mathbf{v}^{+}\mathbf{(}x_{1}))(\mathbf{v}%
^{+}\mathbf{(}x_{1}))\nonumber\\
&  =(-{}^{t}\mathbf{v}^{+}\mathbf{(}x_{1})\cos(x_{1}+3\pi/2)\Delta
(2,2)\mathbf{v}^{+}\mathbf{(}x_{1}),{}^{t}\mathbf{v}^{+}\mathbf{(}x_{1}%
)\sin(x_{1}+3\pi/2)\Delta(2,2)\mathbf{v}^{+}\mathbf{(}x_{1}))\nonumber\\
&  =2\mathbf{w}^{+}\mathbf{(}x_{1}),
\end{align}
and similarly $(d_{(x_{1},0)}^{2}\omega^{-})(\mathbf{v}^{-}\mathbf{(}%
x_{1}))(\mathbf{v}^{-}\mathbf{(}x_{1}))=2\mathbf{w}^{-}\mathbf{(}x_{1})$. Let
$(x_{1},0)\in([-4,-\pi]\cup\lbrack\pi,4])\times0$. Then $K(\omega^{\pm
})_{(x_{1},0)}$ is generated by $^{t}(0,1)$, $Q(\omega^{\pm})_{(x_{1},0)}$ is
generated and oriented by $^{t}(0,1)$ and $(d_{(x_{1},0)}^{2}\omega^{\pm})($
$^{t}(0,1))($ $^{t}(0,1))=2$ $^{t}(0,1)$. By the definition (2) of
$\omega^{\pm}$\ we have that $\omega^{\pm}(x_{1},x_{2})\in\Sigma
^{1,0}(\mathbf{R}^{2},\mathbf{R}^{2})$ for $(x_{1},x_{2})\in C(0,\pm2)$.
Hence, we have $\omega^{\pm}\in\Gamma(\mathbf{R}^{2},\mathbf{R}^{2})$.

Next we show the transversality. Since $(x_{1},x_{2})\mapsto h^{\pm}%
(x_{1},x_{2}-(\pm2))$ is a fold-map around $C(0,\pm2),$ and since $\pi_{1}%
^{2}\circ\pi_{\Omega}\circ\omega^{\pm}(x_{1},x_{2})=Jh^{\pm}(z^{\pm}%
(x_{1},x_{2}))$ and $z^{\pm}$ is an immersion around $[-\pi,\pi]\mathbf{\times
}0$, $\omega^{\pm}$ is transverse to $\Sigma^{1}(\mathbf{R}^{2},\mathbf{R}%
^{2})$ on $C(0,\pm2)$ and on $[-\pi,\pi]\mathbf{\times}0$. Let $(x_{1}%
,0)\in([-4,-\pi]\cup\lbrack\pi,4])\times0$. Then we have that the $(2,2)$
component of$\ d_{(x_{1},0)}(\widetilde{\mathbf{a}^{\pm}})(\partial/\partial
x_{2})$ is equal to $d/dx_{2}(\phi^{\pm}(|x_{1}|,x_{2}))|_{x_{2}=0}>0$ by the
comment exactly below (ii-3). By this observation and (5), $\omega^{\pm}$ is
transverse to $\Sigma^{1}(\mathbf{R}^{2},\mathbf{R}^{2})$\ at $(x_{1},0)$ for
$\left|  x_{1}\right|  >\pi$. Hence, $\omega^{\pm}\in\Gamma^{tr}%
(\mathbf{R}^{2},\mathbf{R}^{2})$.

We prove the assertion Lemma 5.5 (4) for $\omega^{\pm}.$\ Since $(x_{1}%
,x_{2})\mapsto h^{\pm}(x_{1},x_{2}-(\pm2))$ is a fold-map around $C(0,\pm2)$,
we have $d(i_{N(C(0,\pm2))},\Phi(\omega^{\pm})|_{C(0,\pm2)})=0.$ Next we show
$d(i_{N(\mathbf{R}\times0)},\Phi(\omega^{\pm})|_{\mathbf{R}\times0})=\mp1$. We
know that $i_{N([-\pi,\pi]\times0)}$ is the inclusion. If we prove that
$\Phi(\omega^{\pm})_{(x_{1},0)}(\mathbf{e}_{2})=i_{K(\omega^{\pm})}\circ
(d^{2}(\omega^{\pm})|K(\omega^{\pm}))^{-1}\circ(d^{2}(\omega^{\pm}%
)_{\thicksim}\circ d(\omega^{\pm})_{\thicksim}|_{(x_{1},0)})(\mathbf{e}_{2})$
is equal to $(\sqrt{2}-1)\mathbf{v}^{\pm}\mathbf{(}x_{1})$ for $\left|
x_{1}\right|  \leq\pi$, then we have%
\begin{equation}
d(i_{N([-\pi,\pi]\times0)},\Phi(\omega^{\pm})|_{[-\pi,\pi]\times0})=\mp1.
\end{equation}
In fact, let $p_{Q}$ be the orthogonal projection of $\mathbf{R}^{2}$ onto the
space generated by $\mathbf{w}^{\pm}\mathbf{(}x_{1})$. If we regard
$d_{(x_{1},0)}\mathbf{a}^{\pm}:T_{(x_{1},0)}\mathbf{R}^{2}\rightarrow
T_{\mathbf{a}^{\pm}(x_{1},0)}M(2)$ as a map $T_{(x_{1},0)}\mathbf{R}%
^{2}\rightarrow M(2)$, then we have, by the definition of the intrinsic
derivative, that $((d^{2}(\omega^{\pm})_{\thicksim}\circ d(\omega^{\pm
})_{\thicksim}|_{(x_{1},0)})(\mathbf{e}_{2}))(\mathbf{v}^{\pm}\mathbf{(}%
x_{1}))$ is equal to $p_{Q}((d_{(x_{1},0)}\mathbf{a}^{\pm}(\partial/\partial
x_{2}))\mathbf{v}^{\pm}\mathbf{(}x_{1}))$. In the following calculation let
$a=x_{1}+3\pi/2$ and $b=(1-1/\sqrt{2})x_{2}+1/\sqrt{2}$ for $\omega^{+},$ and
let $a=-x_{1}+\pi/2$ and $b=-(1-1/\sqrt{2})x_{2}+1/\sqrt{2}$ for $\omega^{-}$.
Since $\mathbf{v}^{\pm}\mathbf{(}x_{1})=\pm$ $^{t}(\cos a,\sin a)$ is
independent of $x_{2}$, $p_{Q}((d_{(x_{1},0)}\mathbf{a}^{\pm}(\partial
/\partial x_{2}))\mathbf{v}^{\pm}\mathbf{(}x_{1}))$ is equal to the value at
$(x_{1},0)$ of
\begin{align*}
&  p_{Q}(\partial/\partial x_{2}(\mathbf{a}^{\pm}(x_{1},x_{2})\mathbf{v}^{\pm
}\mathbf{(}x_{1})))\\
&  =p_{Q}(\partial/\partial x_{2}(Jh^{\pm}(z^{\pm}(x_{1},x_{2}))\mathbf{v}%
^{\pm}\mathbf{(}x_{1})))\\
&  =p_{Q}(\partial/\partial x_{2}(Jh^{\pm}(b\cos a,b\sin a)\mathbf{v}^{\pm
}\mathbf{(}x_{1})))\\
&  =p_{Q}\left(  {\displaystyle\frac{\partial}{\partial x_{2}}}\left(
e^{1/2-b^{2}}\left(
\begin{array}
[c]{cc}%
1-2b^{2}\cos^{2}a & -2b^{2}\sin a\cos a\\
\pm2b^{2}\sin a\cos a & \pm(2b^{2}\sin^{2}a-1)
\end{array}
\right)  \left(
\begin{array}
[c]{l}%
\pm\cos a\\
\pm\sin a
\end{array}
\right)  \right)  \right) \\
&  =\pm\partial/\partial x_{2}(e^{1/2-b^{2}}(-(2b^{2}-1)\cos a,\pm
(2b^{2}-1)\sin a))\\
&  =\pm(\partial/\partial x_{2}(e^{1/2-b^{2}}(2b^{2}-1)))\mathbf{w}^{\pm
}\mathbf{(}x_{1}).
\end{align*}
The coefficient of $\mathbf{w}^{\pm}\mathbf{(}x_{1})$ at $(x_{1},0)$ is equal
to $2(\sqrt{2}-1)$. By (3$^{\pm}$) we have $\Phi(\omega^{\pm})_{(x_{1}%
,0)}(\mathbf{e}_{2})=(\sqrt{2}-1)\mathbf{v}^{\pm}\mathbf{(}x_{1})$. This
proves the assertion Lemma 5.5 (4), since the map $(x_{1},0)\rightarrow
\mathbf{v}^{\pm}\mathbf{(}x_{1})$ is of degree $\pm1.$

Finally we construct a homotopy $\pi_{\Omega}\circ\omega^{\pm}(\lambda
):\mathbf{R}^{2}\rightarrow\Omega^{1,0}(2,2)$, $\lambda\in\lbrack0,1]$. First
we set $\pi_{\Omega}\circ\omega^{\pm}(t)=\pi_{\Omega}\circ j^{2}\sigma$,
$\pi_{\Omega}\circ\omega^{\pm}(1-t)=\pi_{\Omega}\circ\omega^{\pm}$ for
sufficiently small $t$ with $t\geq0,$ and $\pi_{\Omega}\circ\omega^{\pm
}(\lambda)|(\mathbf{R}^{2}\setminus{A^{\pm})}=\pi_{\Omega}\circ j^{2}%
\sigma|(\mathbf{R}^{2}\setminus{A^{\pm})}$. Since $\pi_{2}(\Omega
^{1,0}(2,2))\cong$ $\pi_{2}(SO(3))\cong\{0\}$ ([An2]), this can be extended to
the whole space $\mathbf{R}^{2}\times\lbrack0,1]$. Then we obtain the required
homotopy $\omega_{\lambda}^{\pm}\in\Gamma(\mathbf{R}^{2},\mathbf{R}^{2})$
defined by $\omega_{\lambda}^{\pm}(x_{1},x_{2})=((x_{1},x_{2}),\sigma
(x_{1},x_{2});\pi_{\Omega}\circ\omega^{\pm}(\lambda)(x_{1},x_{2}))$. This
completes the proof.
\end{proof}

For $n>2$, we can construct a section $\varpi^{\pm}\in\Gamma^{tr}%
(\mathbf{R}^{n},\mathbf{R}^{n})$ such as $\omega^{\pm}\in\Gamma^{tr}%
(\mathbf{R}^{2},\mathbf{R}^{2})$ by generalizing the method given in the proof
of Lemma 5.5 as follows. However, we cannot assert that $\varpi^{\pm}$ is
homotopic to $j^{2}g_{n}$ in $\Gamma(\mathbf{R}^{n},\mathbf{R}^{n})$. The
proof is left to the reader.

\begin{lemma}
For $n\geq2$, there exists a section $\varpi^{\pm}\in\Gamma^{tr}%
(\mathbf{R}^{n},\mathbf{R}^{n})$ satisfying the properties:

$(1)$ $\varpi^{\pm}|(\mathbf{R}^{n}\setminus D_{8}^{n})=j^{2}g_{n}%
|(\mathbf{R}^{n}\setminus D_{8}^{n})$ and $\pi_{\mathbf{R}^{n}}^{2}\circ
\varpi^{\pm}=g_{n}$,

$(2)$ $S(\varpi^{\pm})$ is the union of $\mathbf{R}^{n-1}\times0$ and the
$(n-1)$-sphere $C(\mathbf{0}_{n-1},\pm2)$ centered at $(\mathbf{0}_{n-1}%
,\pm2)$ with radius $1/\sqrt{2},$

$(3)$ $d(i_{N(\mathbf{R}^{n-1}\times0)},\Phi(\varpi^{\pm})|_{\mathbf{R}%
^{n-1}\times0})=\mp1$ and $d(i_{N(C(\mathbf{0}_{n-1},\pm2))},\Phi(\varpi^{\pm
})|_{C(\mathbf{0}_{n-1},\pm2)})=0$.
\end{lemma}

By applying $\mu_{\lambda}^{i}$ in Lemma 5.3\ to (5.3) for $p>2$ and by using
Lemma 5.5 for $p=2$, we prove the following.

\begin{proposition}
Let $n\geq p\geq2$. $(1)$ If $p\geq3$, then the homotopy $\mathcal{J}^{\iota
}(\omega(\mu_{\lambda}^{i}))$ $(1\leq i\leq p-1)$ relative to $\mathbf{R}%
^{n}\setminus\mathrm{Int}D_{10}^{n}$ in $\Gamma(\mathbf{R}^{n},\mathbf{R}%
^{p})$ satisfies that

$(1$-$\mathrm{i})$ $\mathcal{J}^{\iota}(\omega(\mu_{0}^{i}))=j^{2}{g}^{\iota}$,

$(1$-$\mathrm{ii})$ $\mathcal{J}^{\iota}(\omega(\mu_{1}^{i}))$ is a smooth
section transverse to $\Sigma^{n-p+1,0}(\mathbf{R}^{n},\mathbf{R}^{p})$ and
$S(\mathcal{J}^{\iota}(\omega(\mu_{1}^{i})))=(\mu_{1}^{i})^{-1}(0)\times
\mathbf{0}_{n-p}$ is connected, where $p\geq3$ and $1\leq i<p-1$, or $p=3$ and
$i=2$,

$(1$-$\mathrm{iii})$ if $p\geq3$ and $1\leq i<p-1$, then $d(i_{\mathcal{N}%
(\mathcal{J}^{\iota}(\omega(\mu_{1}^{i})))},\Phi(\mathcal{J}^{\iota}%
(\omega(\mu_{1}^{i}))))$ is equal to $(-1)^{i}$,

$(1$-$\mathrm{iv})$ if $p=3$ and $i=2$, then $d(i_{\mathcal{N}(\mathcal{J}%
^{\iota}(\omega(\mu_{1}^{i})))},\Phi(\mathcal{J}^{\iota}(\omega(\mu_{1}%
^{i}))))=1$.

\noindent$(2)$ If $p=2$, then there exists a homotopy $\mathcal{J}^{\iota
}(\omega_{\lambda}^{\pm})$ relative to $\mathbf{R}^{n}\setminus\mathrm{Int}%
D_{10}^{n}$ in $\Gamma(\mathbf{R}^{n},\mathbf{R}^{2})$ such that

$(2$-$\mathrm{i})$ $\mathcal{J}^{\iota}(\omega_{0}^{\pm})=j^{2}{g}^{\iota}$,

$(2$-$\mathrm{ii})$ $\mathcal{J}^{\iota}(\omega_{1}^{\pm})$ is transverse to
$\Sigma^{n-1,0}(\mathbf{R}^{n},\mathbf{R}^{2})$ and $S(\mathcal{J}^{\iota
}(\omega_{1}^{\pm}))$ consists of $\mathbf{R\times0}_{n-1}$ and the circle
$C(0,\pm2)\times\mathbf{0}_{n-2}$ $($denoted by $\mathbf{C}(0,\pm2)$ for short$),$

$(2$-$\mathrm{iii})$ $d(i_{\mathcal{N}(\mathcal{J}^{\iota}(\omega_{1}^{\pm}%
))}|_{\mathbf{R\times0}_{n-1}},\Phi(\mathcal{J}^{\iota}(\omega_{1}^{\pm
}))|_{\mathbf{R\times0}_{n-1}})=\mp1,$

$(2$-$\mathrm{iv})$ $d(i_{\mathcal{N}(\mathcal{J}^{\iota}(\omega_{1}^{\pm}%
))}|_{\mathbf{C}(0,\pm2)},\Phi(\mathcal{J}^{\iota}(\omega_{1}^{\pm
}))|_{\mathbf{C}(0,\pm2)})=0$.

\noindent$(3)$ The line bundles $Q(\mathcal{J}^{\iota}(\omega(\mu_{1}^{i})))$
for $p>2$ and $Q(\mathcal{J}^{\iota}(\omega_{1}^{\pm}))$ for $p=2$ are trivial
and oriented, whose induced orientations on $S(\mathcal{J}^{\iota}(\omega
(\mu_{1}^{i})))\setminus D_{10}^{n}$ and $S(\mathcal{J}^{\iota}(\omega
_{1}^{\pm}))\setminus D_{10}^{n}$ are given by $\mathbf{e}_{p}$.
\end{proposition}

\begin{proof}
The proof for the case $n=p$ is easier. Hence, we assume that $n>p\geq2$.

(1-i) By Lemma 5.3 we have $\mu_{0}^{i}(\overset{\bullet}{x})=2x_{p}.$ Hence,
we have $\mathcal{J}^{\iota}(\omega(\mu_{0}^{i}))${$=j^{2}{g}^{\iota}$.}

(1-ii) It is clear that $S(\mathcal{J}^{\iota}(\omega(\mu_{1}^{i})))=(\mu
_{1}^{i})^{-1}(0)\mathbf{\times0}_{n-p}$, which is connected by Lemma 5.3 (4).
The first assertion follows from the fact stated just below (5.5).

(1-iii, iv) Set $L=(\mu_{1}^{i})^{-1}(0)$. By Lemmas 5.1 and Lemma 5.3 we have
$d(i_{N(L)},\mathbf{n}_{L})=\mathrm{deg}(\mathbf{e}(L))=(-1)^{i}$. By Lemma
5.2 we have $d(i_{\mathcal{N}(\mathcal{J}^{\iota}(\omega(\mu_{1}^{i})))}%
,\Phi(\mathcal{J}^{\iota}(\omega(\mu_{1}^{i}))))\linebreak =m(p,n)_{\ast
}^{\prime}(d(i_{N(L)},\mathbf{n}_{L}))$, which is equal to $(-1)^{i}$.

(2) The proof is rather long. Under the identification (1.6), we define
$\mathcal{J}^{\iota}(\omega_{\lambda}^{\pm}):\mathbf{R}^{n}\rightarrow
\Omega^{n-1,0}(\mathbf{R}^{n},\mathbf{R}^{2})$ by
\[
\mathcal{J}^{\iota}(\omega_{\lambda}^{\pm})(x)=(x,g^{\iota}(x);\mathcal{J}%
_{2,n}^{\iota}(\pi_{\Omega}\circ\omega_{\lambda}^{\pm})(x,\lambda)),
\]
where $\mathcal{J}_{2,n}^{\iota}(\pi_{\Omega}\circ\omega_{\lambda}^{\pm
})(x,\lambda)=$ $(\mathbf{A}_{2\times n}^{\iota}(x,\lambda);\mathbf{B}%
_{n\times n}^{\iota,1}(x,\lambda),\mathbf{B}_{n\times n}^{\iota,2}%
(x,\lambda))$ is defined as follows.

First, we may assume that if $\lambda\geq9/10,$ then $\omega_{\lambda}^{\pm
}(\overset{\bullet}{x})=\omega_{1}^{\pm}(\overset{\bullet}{x})$. For the
oriented vector $\mathbf{e}(Q(\omega_{\lambda}^{\pm}))$ we can construct an
extended vector field $\mathbf{e}_{\omega^{\pm}}(\overset{\bullet}{x}%
,\lambda)={}^{t}(e_{\omega^{\pm}}^{1}(\overset{\bullet}{x},\lambda
),e_{\omega^{\pm}}^{2}(\overset{\bullet}{x},\lambda))\in T_{\overset{\bullet
}{x}}\mathbf{R}^{2}$ with $\overset{\bullet}{x}\in\mathbf{R}^{2}$ and
$\lambda\in\lbrack0,1]$ with the following properties.

$\blacktriangleright$ $\mathbf{e}_{\omega^{\pm}}(\overset{\bullet}{x}%
,\lambda)$ is continuous with respect to $\overset{\bullet}{x}$ and $\lambda.$

$\blacktriangleright$ If $\overset{\bullet}{c}\in S(\omega_{\lambda}^{\pm})$,
then $\mathbf{e}_{\omega^{\pm}}(\overset{\bullet}{c},\lambda)=\mathbf{e}%
(Q(\omega_{\lambda}^{\pm})_{\overset{\bullet}{c}})$. In particular, if
$\overset{\bullet}{c}\notin$Int$A^{\pm}$ in addition, then $\mathbf{e}%
_{\omega^{\pm}}(\overset{\bullet}{c},\lambda)=$ $^{t}(0,1)$.

$\blacktriangleright$ If either $\lambda=0$ or $\overset{\bullet}{x}\notin
$Int$A^{\pm}$, then $\mathbf{e}_{\omega^{\pm}}(\overset{\bullet}{x},\lambda)=$
$^{t}(0,1)$.

$\blacktriangleright$ $\mathbf{e}_{\omega^{\pm}}(\overset{\bullet}{x}%
,\lambda)$ is smooth around $\lambda=0.$

$\blacktriangleright$ If $\lambda\geq9/10,$ then $\mathbf{e}_{\omega^{\pm}%
}(\overset{\bullet}{c},\lambda)=\mathbf{e}_{\omega^{\pm}}(\overset{\bullet}%
{c},1)$.

\noindent This is possible since the subset consisting of all $(S(\omega
_{\lambda}^{\pm}),\lambda)$ in $\mathbf{R}^{2}\times\mathbf{R}$ is a closed
set. In the following argument we need to extend the domain $[0,1]$ of
$\lambda$ to $(-\infty,1]$ by setting $\omega_{\lambda}^{\pm}=\omega_{0}^{\pm
}$ and $\mathbf{e}_{\omega^{\pm}}(\overset{\bullet}{c},\lambda)=\mathbf{e}%
_{\omega^{\pm}}(\overset{\bullet}{c},0)$ for $\lambda<0$. Suppose that
$\pi_{\Omega}\circ\omega_{\lambda}^{\pm}(\overset{\bullet}{x})$ is written as
\begin{equation}
(A_{2\times2}(\overset{\bullet}{x},\lambda);B_{2\times2}^{1}(\overset{\bullet
}{x},\lambda),B_{2\times2}^{2}(\overset{\bullet}{x},\lambda))
\end{equation}
under the identification (1.6). We define the $2\times n$ matrix
$\mathbf{A}_{2\times n}^{\iota}(x,\lambda)$ for $x\in\mathbf{R}^{n}$ and
$\lambda\in\lbrack0,1]$ to be
\begin{equation}
\left(
\begin{array}
[c]{ll}%
A_{2\times2}(\overset{\bullet}{x},\lambda-\Vert\overset{\bigstar}{x}\Vert
^{2}) &
\begin{array}
[c]{l}%
e_{\omega^{\pm}}^{1}(\overset{\bullet}{x},\lambda-\Vert\overset{\bigstar}%
{x}\Vert^{2})(2x_{3},\ldots,2x_{n-\iota},-2x_{n-\iota+1},\ldots,-2x_{n})\\
e_{\omega^{\pm}}^{2}(\overset{\bullet}{x},\lambda-\Vert\overset{\bigstar}%
{x}\Vert^{2})(2x_{3},\ldots,2x_{n-\iota},-2x_{n-\iota+1},\ldots,-2x_{n})
\end{array}
\end{array}
\right)  ,
\end{equation}
which satisfies the following properties.

(\textbf{A}1) If $\overset{\bullet}{x}\notin$Int$A^{\pm}$, then%
\[
\mathbf{A}_{2\times n}^{\iota}(x,\lambda)=\left(
\begin{array}
[c]{ll}%
1 & \mathbf{0}_{1\times(n-1)}\\
0 & (2x_{2},\ldots,2x_{n-\iota},-2x_{n-\iota+1},\ldots,-2x_{n})
\end{array}
\right)  .
\]

(\textbf{A}2) If $x_{2}\neq0$ or $\Vert\overset{\bigstar}{x}\Vert\neq0$, and
$\overset{\bullet}{x}\notin$Int$A^{\pm}$, then rank$\mathbf{A}_{2\times
n}^{\iota}(x,\lambda)=2.$

(\textbf{A}3) If $\overset{\bullet}{x}\in{A^{\pm}}$ and $\Vert\overset
{\bigstar}{x}\Vert\neq0$, then rank$\mathbf{A}_{2\times n}^{\iota}%
(x,\lambda)=2$.

(\textbf{A}4) If $\overset{\bullet}{x}\in{A^{\pm}}$ and $\Vert\overset
{\bigstar}{x}\Vert=0$, then
\[
\mathbf{A}_{2\times n}^{\iota}(x,\lambda)=(A_{2\times2}(\overset{\bullet}%
{x},\lambda),\mathbf{0}_{2\times(n-2)}),
\]
which is of rank $1$ if and only if $A_{2\times2}(\overset{\bullet}{x}%
,\lambda)$ is of rank $1$.

(\textbf{A}5) If $\Vert\overset{\bigstar}{x}\Vert\geq1$, then $\pi_{\Omega
}\circ\omega_{\lambda-\Vert\overset{\bigstar}{x}\Vert^{2}}^{\pm}=\pi_{\Omega
}\circ\omega_{0}^{\pm}=\pi_{\Omega}\circ j^{2}\sigma$, and hence,
\[
\mathbf{A}_{2\times n}^{\iota}(x,\lambda)=\mathbf{A}_{2\times n}^{\iota
}(x,0)=\left(
\begin{array}
[c]{ll}%
1 & \mathbf{0}_{1\times(n-1)}\\
0 & (2x_{2},\ldots,2x_{n-\iota},-2x_{n-\iota+1},\ldots,-2x_{n})
\end{array}
\right)  .
\]
The properties (\textbf{A}2) and (\textbf{A}3) follow from the fact that if
rank $A_{2\times2}(\overset{\bullet}{x},\lambda-\Vert\overset{\bigstar}%
{x}\Vert^{2})=1$, then the two column vectors of $A_{2\times2}(\overset
{\bullet}{x},\lambda-\Vert\overset{\bigstar}{x}\Vert^{2})$ is orthogonal to
the nonzero vector $\mathbf{e}_{\omega^{\pm}}(\overset{\bullet}{x}%
,\lambda-\Vert\overset{\bigstar}{x}\Vert^{2}).$

We set, for $j=1$, $2$,
\begin{equation}
\mathbf{B}_{n\times n}^{\iota,j}(x,\lambda)=\left(
\begin{array}
[c]{ll}%
B_{2\times2}^{j}(\overset{\bullet}{x},\lambda-\Vert\overset{\bigstar}{x}%
\Vert^{2}) & \mathbf{0}_{2\times(n-2)}\\
\mathbf{0}_{(n-2)\times2} & e_{\omega^{\pm}}^{j}(\overset{\bullet}{x}%
,\lambda-\Vert\overset{\bigstar}{x}\Vert^{2})\Delta\lbrack2;\iota]
\end{array}
\right)  ,
\end{equation}
where $\Delta\lbrack2;\iota]$ is, here, an $(n-2)\times(n-2)$ matrix. These
two matrices have the following property.

(\textbf{B}1) If $\overset{\bullet}{x}\notin A^{\pm}$, or $\Vert
\overset{\bigstar}{x}\Vert\geq1$, then $\mathbf{B}_{n\times n}^{\iota
,1}(x,\lambda)=\mathbf{0}_{n\times n}$ and $\mathbf{B}_{n\times n}^{\iota
,2}(x,\lambda)=(0)\dotplus(2)\dotplus\Delta\lbrack2;\iota]$.

We show $\mathcal{J}^{\iota}(\omega_{\lambda}^{\pm})\in\Gamma(\mathbf{R}%
^{n},\mathbf{R}^{2})$. Take a point $(c,\lambda)$ such that rank$\mathbf{A}%
_{2\times n}^{\iota}(c,\lambda)=1$. If $\overset{\bullet}{c}\notin$%
Int${A^{\pm}}$, then $c=(c_{1},0,\ldots,0)$ by (\textbf{A}1), and if
$\overset{\bullet}{c}=(c_{1},c_{2})\in{A^{\pm}}$, then $c=(c_{1}%
,c_{2},0,\ldots,0)$ and rank$A_{2\times2}(\overset{\bullet}{c},\lambda)=1$ by
(\textbf{A}3) and (\textbf{A}4). Set $\mathbf{e}(\mathbf{A}_{2\times n}%
^{\iota}(c,\lambda))\linebreak =(\mathbf{e}(K(\omega_{\lambda}^{\pm
})_{\overset{\bullet}{c}}),0,\ldots,0)$. Here, let $\mathbf{e}^{\prime}$ be a
vector orthogonal to $K(\omega_{\lambda}^{\pm})_{\overset{\bullet}{c}}$. We
orient $K(\omega_{\lambda}^{\pm})_{\overset{\bullet}{c}}$ so that the
juxtapositions $(\mathbf{e}^{\prime},\mathbf{e}(K(\omega_{\lambda}^{\pm
})_{\overset{\bullet}{c}}))$ and $(A_{2\times2}(\overset{\bullet}{c}%
,\lambda)(\mathbf{e}^{\prime}),\mathbf{e}(Q(\omega_{\lambda}^{\pm}%
)_{\overset{\bullet}{c}}))$ induce the canonical orientation of $\mathbf{R}%
^{2}$. Then the kernel of $\mathbf{A}_{2\times n}^{\iota}(c,\lambda)$ is
generated by $\mathbf{e}(\mathbf{A}_{2\times n}^{\iota}(c,\lambda
)),\mathbf{e}_{3},\ldots,\mathbf{e}_{n}$. From the above definition of
$\mathcal{J}^{\iota}(\omega_{\lambda}^{\pm})$ and the fact that $\omega
_{\lambda}^{\pm}\in\Gamma(\mathbf{R}^{2},\mathbf{R}^{2})$ it follows that%
\begin{equation}
q(\mathcal{J}^{\iota}(\omega_{\lambda}^{\pm}))_{c}(\mathbf{e}(\mathbf{A}%
_{2\times n}^{\iota}(c,\lambda)),\mathbf{e}(\mathbf{A}_{2\times n}^{\iota
}(c,\lambda)\mathbf{))}=q(\omega_{\lambda}^{\pm})_{\overset{\bullet}{c}%
}(\mathbf{e(}K(\omega_{\lambda}^{\pm})_{\overset{\bullet}{c}}),\mathbf{e(}%
K(\omega_{\lambda}^{\pm})_{\overset{\bullet}{c}}))
\end{equation}
is a positive multiple of $\mathbf{e(}Q\mathbf{(}\omega_{\lambda}^{\pm
})_{\overset{\bullet}{c}})$\ (see the equality just below (5.15) for
$\lambda=1$). Therefore, $Q(\mathcal{J}^{\iota}(\omega_{\lambda}^{\pm}))_{c}$
is generated and oriented by $\mathbf{e}(Q(\omega_{\lambda}^{\pm}%
)_{\overset{\bullet}{c}})$. Furthermore, we have, by (5.11),%
\begin{equation}%
\begin{array}
[c]{ll}%
q(\mathcal{J}^{\iota}(\omega_{\lambda}^{\pm}))_{c}(\mathbf{e}(\mathbf{A}%
_{2\times n}^{\iota}\mathbf{(}c\mathbf{,\lambda)),e}_{j})=\mathbf{0} &
\text{for }2<j\leq n,\\
q(\mathcal{J}^{\iota}(\omega_{\lambda}^{\pm}))_{c}(\mathbf{e}_{k}%
\mathbf{,e}_{\ell})=2\delta_{k\ell}\mathbf{e(}Q\mathbf{(}\omega_{\lambda}%
^{\pm})_{\overset{\bullet}{c}}) & \text{for }2<k\leq n-\iota\text{ and }%
2<\ell\leq n,\\
q(\mathcal{J}^{\iota}(\omega_{\lambda}^{\pm}))_{c}(\mathbf{e}_{k}%
\mathbf{,e}_{\ell})=-2\delta_{k\ell}\mathbf{e(}Q\mathbf{(}\omega_{\lambda
}^{\pm})_{\overset{\bullet}{c}}) & \text{for }n-\iota<k\leq n\text{ and
}2<\ell\leq n.
\end{array}
\end{equation}
This shows that $q(\mathcal{J}^{\iota}(\omega_{\lambda}^{\pm}))_{c}$ is
nonsingular and that the index of $q(\mathcal{J}^{\iota}(\omega_{\lambda}%
^{\pm}))_{c}$ is equal to $\iota$. Hence, we have $\mathcal{J}^{\iota}%
(\omega_{\lambda}^{\pm})\in\Gamma(\mathbf{R}^{n},\mathbf{R}^{2})$.

Now we are ready to prove (2-i) to (2-iv). It follows from (\textbf{A}1, 5)
and (\textbf{B}1) that if either $\overset{\bullet}{x}\notin$Int$A^{\pm}$ or
$\Vert\overset{\bigstar}{x}\Vert\geq1$, then we have $\mathcal{J}^{\iota
}(\omega_{\lambda}^{\pm})(x)=j^{2}g^{\iota}(x)$ for any $\lambda\in
\lbrack0,1]$.

(2-i) From (\textbf{A}1, 4, 5) and (5.11) it follows that $\mathcal{J}^{\iota
}(\omega_{0}^{\pm})(x)=j^{2}g^{\iota}(x)$.

(2-ii) If $\Vert\overset{\bigstar}{x}\Vert\neq0$, then $\mathcal{J}^{\iota
}(\omega_{1}^{\pm})(x)$ is a non-singular jet, and hence $S(\mathcal{J}%
^{\iota}(\omega_{1}^{\pm}))$ consists of $\mathbf{R}\times\mathbf{0}_{n-1}$
and $C(0,\pm2)\times\mathbf{0}_{n-2}$ by Lemma 5.5 (3) and (5.9). In (5.10)
$A_{2\times2}(\overset{\bullet}{x},\lambda-\Vert\overset{\bigstar}{x}\Vert
^{2})$ and $\mathbf{e}_{\omega}(\overset{\bullet}{x},\lambda-\Vert
\overset{\bigstar}{x}\Vert^{2})$ do not depend on the variables $x_{3}%
,\ldots,x_{n}$ around $\lambda=1,$ where $\Vert\overset{\bigstar}{x}\Vert$ is
sufficiently small so that $\lambda-\Vert\overset{\bigstar}{x}\Vert^{2}>9/10$.
Hence $\mathcal{J}^{\iota}(\omega_{1}^{\pm})$ is transverse to $\Sigma
^{n-p+1}(\mathbf{R}^{n},\mathbf{R}^{2})$ by Lemma 5.5 (2). See also the proof
of (2-iv) below.

(2-iii) Note that $K(\mathcal{J}^{\iota}(\omega_{1}^{\pm}))_{c}$ and
$Q(\mathcal{J}^{\iota}(\omega_{1}^{\pm}))_{c}$ are identified with
$K(\omega_{1}^{\pm})_{\overset{\bullet}{c}}\times\mathbf{0}_{n-2}%
\oplus\mathbf{0}_{2}\times\mathbf{R}^{n-2}$ and $Q(\omega_{1}^{\pm}%
)_{\overset{\bullet}{c}}$ respectively. We define%
\[
m(2,n):\text{Mono}(K(\omega_{1}^{\pm})|_{\mathbf{R\times}0},T\mathbf{R}%
^{2}|_{\mathbf{R\times}0})\rightarrow\text{Mono}(K(\mathcal{J}^{\iota}%
(\omega_{1}^{\pm}))|_{\mathbf{R\times0}_{n-1}},T\mathbf{R}^{n}%
|_{\mathbf{R\times0}_{n-1}})
\]
by $m(2,n)(h)=h\oplus id_{\mathbf{0}_{2}\times\mathbf{R}^{n-2}}$, where $h\in
$Mono$(K(\omega_{1}^{\pm})|_{\mathbf{R\times}0},T\mathbf{R}^{2}%
|_{\mathbf{R\times}0})$. Furthermore, the canonical map
\[
\overline{m(2,n)}:\text{Mono}(\mathbf{R},\mathbf{R}^{2})\rightarrow
\text{Mono}(\mathbf{R\oplus R}^{n-2},\mathbf{R}^{n})
\]
defined similarly as $m(2,n),$ induces a surjection
\[
\overline{m(2,n)}_{\ast}:\pi_{1}(\text{Mono}(\mathbf{R},\mathbf{R}%
^{2}))\rightarrow\pi_{1}(\text{Mono}(\mathbf{R\oplus R}^{n-2},\mathbf{R}%
^{n})).
\]
From the proof of (5.8) it follows that
\[%
\begin{array}
[c]{ll}%
(i_{\mathcal{N}(\mathcal{J}^{\iota}(\omega_{1}^{\pm}))}|_{\mathbf{R\times
0}_{n-1}})(\mathbf{e}_{\ell})=\mathbf{e}_{\ell} & \text{for }2\leq\ell\leq
n,\\
(\Phi(\mathcal{J}^{\iota}(\omega_{1}^{\pm}))|_{\mathbf{R\times0}_{n-1}%
})(\mathbf{e}_{\ell})=\mathbf{e}_{\ell} & \text{for }2<\ell\leq n,\\
(\Phi(\mathcal{J}^{\iota}(\omega_{1}^{\pm}))|_{\mathbf{R\times0}_{n-1}%
})_{(x_{1},0)}(\mathbf{e}_{2})=(\sqrt{2}-1)\mathbf{v}^{\pm}(x_{1}) & \text{for
}|x_{1}|\leq\pi.
\end{array}
\]
Therefore, we obtain that
\[
d(i_{\mathcal{N}(\mathcal{J}^{\iota}(\omega_{1}^{\pm}))}|_{\mathbf{R\times
0}_{n-1}},\Phi(\mathcal{J}^{\iota}(\omega_{1}^{\pm}))|_{\mathbf{R\times
0}_{n-1}})=\mp1.
\]

(2-iv) Let $c^{\pm}(\theta)=((1/\sqrt{2})\cos\theta,(1/\sqrt{2})\sin\theta
\pm2,\mathbf{0}_{n-2})\in\mathbf{C}(0,\pm2)$, $\overset{\bullet}{c^{\pm}%
}(\theta)\linebreak =((1/\sqrt{2})\cos\theta,(1/\sqrt{2})\sin\theta\pm2)$ and
$\mathbf{e}(\theta)={}^{t}(\cos\theta,\sin\theta,\mathbf{0}_{n-2})$. If we
prove%
\begin{equation}%
\begin{array}
[c]{ll}%
\Phi(\mathcal{J}^{\iota}(\omega_{1}^{\pm}))_{c^{\pm}(\theta)}\left(
\mathbf{e}(\theta)\right)  =\mathbf{e}(\theta), & \\
\Phi(\mathcal{J}^{\iota}(\omega_{1}^{\pm}))_{c^{\pm}(\theta)}(\mathbf{e}%
_{\ell})=\mathbf{e}_{\ell} & \text{for }2<\ell\leq n,
\end{array}
\end{equation}
then we have $d(i_{\mathcal{N}(\mathcal{J}^{\iota}(\omega_{1}^{\pm}%
))}|_{\mathbf{C}(0,\pm2)},\Phi(\mathcal{J}^{\iota}(\omega_{1}^{\pm
}))|_{\mathbf{C}(0,\pm2)})=0$. In fact, $K(\omega_{1}^{\pm})_{\overset
{\bullet}{\mathbf{c}^{\pm}}(\theta)}$ is generated by the vector $^{t}%
(\cos\theta,\sin\theta),$ and $Q(\omega_{1}^{\pm})_{\overset{\bullet
}{\mathbf{c}^{\pm}}(\theta)}$ is generated and oriented by ${}^{t}(-\cos
\theta,\pm\sin\theta)$. If $p_{Q}$ is the orthogonal projection of
$\mathbf{R}^{2}$ onto the space generated by ${}^{t}(-\cos\theta,\pm\sin
\theta),$ and if we regard $d\mathbf{a}^{\pm}:T_{\overset{\bullet}%
{\mathbf{c}^{\pm}}(\theta)}\mathbf{R}^{2}\rightarrow TM(2)$ as a map
$T_{\overset{\bullet}{\mathbf{c}^{\pm}}(\theta)}\mathbf{R}^{2}\rightarrow
M(2)$, then we have%
\begin{align}
&  ((d^{2}(\mathcal{J}^{\iota}(\omega_{1}^{\pm}))_{\thicksim}\circ
d(\mathcal{J}^{\iota}(\omega_{1}^{\pm}))_{\thicksim}|_{c^{\pm}(\theta
)})(\mathbf{e}(\theta)/\sqrt{2}))(\mathbf{e}(\theta))\nonumber\\
&  =((d^{2}(\omega^{\pm})_{\thicksim}\circ d(\omega^{\pm})_{\thicksim
}|_{\overset{\bullet}{\mathbf{c}^{\pm}}(\theta)})((1/\sqrt{2})^{t}(\cos
\theta,\sin\theta)))(^{t}(\cos\theta,\sin\theta)),
\end{align}
which is equal to, as in (5.6),%
\[%
\begin{array}
[c]{l}%
p_{Q}((d/dt(\mathbf{a}^{\pm}(\overset{\bullet}{\mathbf{c}^{\pm}}%
(\theta)+((t-1)/\sqrt{2})(\cos\theta,\sin\theta)))|_{t=1})(^{t}(\cos
\theta,\sin\theta)))\\
=p_{Q}\left(  {\displaystyle\frac{d}{dt}}\left(  e^{1/2-t^{2}/2}\left(
\begin{array}
[c]{cc}%
1-t^{2}\cos^{2}\theta & -t^{2}\sin\theta\cos\theta\\
\pm(t^{2}\sin\theta\cos\theta) & \pm(t^{2}\sin^{2}\theta-1)
\end{array}
\right)  \right)  \left.  \left(
\begin{array}
[c]{l}%
\cos\theta\\
\sin\theta
\end{array}
\right)  \right|  _{t=1}\right)  \\
=2({}^{t}(-\cos\theta,\pm\sin\theta)).
\end{array}
\]
Since%
\[
((d^{2}(\mathcal{J}^{\iota}(\omega_{1}^{\pm}))_{\thicksim}\circ d(\mathcal{J}%
^{\iota}(\omega_{1}^{\pm}))_{\thicksim}|_{c^{\pm}(\theta)})(\mathbf{e}%
_{k}))(\mathbf{v})=p_{Q}((d\mathbf{A}_{2\times n}^{\iota}(c^{\pm}%
(\theta),1)(\partial/\partial x_{k})(\mathbf{v}))|_{c^{\pm}(\theta)})
\]
for $\mathbf{v}\in K(\mathcal{J}^{\iota}(\omega_{1}^{\pm}))$ and since
\[
\partial/\partial x_{k}(\mathbf{A}_{2\times n}^{\iota}(c^{\pm}(\theta
),1))=\left\{
\begin{array}
[c]{cc}%
\left(  \mathbf{0}_{2\times(k-1)}%
\begin{array}
[c]{c}%
2e_{\omega^{\pm}}^{1}(\overset{\bullet}{\mathbf{c}^{\pm}}(\theta),1)\\
2e_{\omega^{\pm}}^{2}(\overset{\bullet}{\mathbf{c}^{\pm}}(\theta),1)
\end{array}
\mathbf{0}_{2\times(n-k)}\right)   & \text{for }2<k\leq n-\iota,\\
\left(  \mathbf{0}_{2\times(k-1)}%
\begin{array}
[c]{c}%
-2e_{\omega^{\pm}}^{1}(\overset{\bullet}{\mathbf{c}^{\pm}}(\theta),1)\\
-2e_{\omega^{\pm}}^{2}(\overset{\bullet}{\mathbf{c}^{\pm}}(\theta),1)
\end{array}
\mathbf{0}_{2\times(n-k)}\right)   & \text{for }n-\iota<k\leq n,
\end{array}
\right.
\]
by (5.10), we have%
\[%
\begin{array}
[c]{cc}%
((d^{2}(\mathcal{J}^{\iota}(\omega_{1}^{\pm}))_{\thicksim}\circ d(\mathcal{J}%
^{\iota}(\omega_{1}^{\pm}))_{\thicksim}|_{c^{\pm}(\theta)})(\mathbf{e}%
_{k}))(^{t}(\cos\theta,\sin\theta))=\mathbf{0}\text{ } & \text{for }2<k\leq n,
\end{array}
\]%
\begin{align}
&  ((d^{2}(\mathcal{J}^{\iota}(\omega_{1}^{\pm}))_{\thicksim}\circ
d(\mathcal{J}^{\iota}(\omega_{1}^{\pm}))_{\thicksim}|_{c^{\pm}(\theta
)})(\mathbf{e}_{k}))(\mathbf{e}_{\ell})\nonumber\\
&  =\left\{
\begin{array}
[c]{cc}%
2\delta_{k\ell}({}^{t}(-\cos\theta,\pm\sin\theta)) & \text{for }2<k\leq
n-\iota\text{ and }2<\ell\leq n,\\
-2\delta_{k\ell}({}^{t}(-\cos\theta,\pm\sin\theta)) & \text{for }n-\iota<k\leq
n\text{ and }2<\ell\leq n.
\end{array}
\right.
\end{align}
Thus (5.12), (5.13), (5.15) and (5.16) prove (5.14) by considering the
definition of $q(s)$ via $d^{2}(s)$ in \S 1. This proves (2-iv).

(3) Let $\sigma_{\lambda}\in\Gamma(\mathbf{R}^{n},\mathbf{R}^{p})$ denote
$\mathcal{J}^{\iota}(\omega(\mu_{\lambda}^{i}))$ for $p\geq3$ and
$\mathcal{J}^{\iota}(\omega_{\lambda}^{\pm})$\ for $p=2$. Since $\sigma
_{0}=j^{2}g^{\iota}$, we orient $Q(\sigma_{0})=\theta_{\mathbf{R}^{p-1}%
\times\mathbf{0}_{n-p+1}}$ by $\mathbf{e}_{p}$. From the construction of
$\sigma_{1}$, (5.2), Remark 5.4 and Lemma 5.5, it follows that $S(\sigma_{1})$
is constructed in $\mathbf{R}^{p}\times\mathbf{0}_{n-p}$ by using $\omega
(\mu_{1}^{i})$ for $p>2$ and $\omega_{1}^{\pm}$ for $p=2$ so that
$S(\sigma_{1})=S(\omega(\mu_{1}^{i}))\times\mathbf{0}_{n-p}$ for $p>2$ and
$S(\sigma_{1})=S(\omega_{1}^{\pm})\times\mathbf{0}_{n-2}$ for $p=2$.
Furthermore, $Q(\sigma_{1})$ is canonically isomorphic to the trivial bundle
$Q(\omega(\mu_{1}^{i}))$ for $p>2$ or $Q(\omega_{1}^{\pm})$ for $p=2$, which
has the orientation on $\mathbf{R}^{p-1}\times\mathbf{0}_{n-p+1}\setminus
D_{10}^{n}$ given by $\mathbf{e}_{p}$.
\end{proof}

\begin{proof}
[Proof of Proposition 2.5.]For the given section $s$, let $M$ be any one of
${M(s)}_{j}^{\iota}$'s and set $m=d(i_{\mathcal{N}(s)}|_{M},\Phi(s)|_{M})$.
Suppose that $m\neq0$. By Remark 4.5 we may assume that we have a point $c\in
M\setminus V(C)$, its small neighborhood $U_{c}$ in $N\setminus V(C)$\ with
suitable coordinates $e:(\mathbf{R}^{p-1}\times\mathbf{R}^{n-p+1}%
,\mathbf{0})\rightarrow(U_{c},c),$ and a small neighborhood $V_{c}\ $of the
point $\pi_{P}^{2}\circ s(c)\in P$ with suitable coordinates such that $s\circ
e$ coincides with $j^{2}g^{\iota}$ on $U_{c}$. We take distinct points
$c_{\ell}\in U_{c}\cap M$ with coordinates $(t_{1}^{\ell},\ldots,t_{p-1}%
^{\ell},0,\ldots,0)$\ and disjoint embedding germs $e_{\ell}:(\mathbf{R}%
^{n},\mathbf{0})\rightarrow(U_{c},c_{\ell})$ such that $e_{\ell}%
(x)=(x_{1}+t_{1}^{\ell},\ldots,x_{p-1}+t_{p-1}^{\ell},x_{p},\ldots,x_{n})$
($1\leq\ell\leq|m|$). Let $\sigma_{\lambda}\in\Gamma(\mathbf{R}^{n}%
,\mathbf{R}^{p})$ denote $\mathcal{J}^{\iota}(\omega(\mu_{\lambda}^{i}))$ for
$p\geq3$ and $\mathcal{J}^{\iota}(\omega_{\lambda}^{\pm})$\ for $p=2$ (if
necessary, we have to change the scale) respectively in Proposition 5.7 such
that%
\[
d(i_{\mathcal{N}(\sigma_{1})},\Phi(\sigma_{1}))=-m/|m|.
\]
For each $e_{\ell}(\mathbf{R}^{n})$, we can construct the homotopy
$\sigma(e_{\ell})_{\lambda}\in\Gamma(e_{\ell}(\mathbf{R}^{n}),P)$ defined by
$\sigma(e_{\ell})_{\lambda}(x)=\sigma_{\lambda}(e_{\ell}^{-1}(x))$.

By using $\sigma(e_{\ell})_{\lambda}$'s for each ${M(s)}_{j}^{\iota}$, we have
a homotopy $s_{\lambda}$ in $\Gamma(N,P)$ defined by%
\[%
\begin{array}
[c]{ll}%
s_{\lambda}|e_{\ell}(\mathbf{R}^{n})=\sigma(e_{\ell})_{\lambda} & \text{on
each }e_{\ell}(\mathbf{R}^{n}),\\
s_{\lambda}|(U_{c}\setminus\cup_{\ell=1}^{|m|}e_{\ell}(\mathbf{R}%
^{n}))=s|(U_{c}\setminus\cup_{\ell=1}^{|m|}e_{\ell}(\mathbf{R}^{n})) &
\text{on each }U_{c}\text{ for }{M(s)}_{j}^{\iota},\\
s_{\lambda}=s & \text{outside of all }U_{c}\text{ for }{M(s)}_{j}^{\iota}.
\end{array}
\]
Then it is easy to see by the additive property of the primary difference that
\linebreak $d(i_{\mathcal{N}(s_{1})}|_{{M(s_{1})}_{j}^{\iota}},\Phi
(s_{1})|_{{M(s_{1})}_{j}^{\iota}})$ is equal to $0$ for all $j$ and $\iota$.
Thus we have proved (1), (2) and (3) of Proposition 2.5.

(4) The triviality of $Q(s_{1})$ follows from the above construction of
$s_{1}$ and Proposition 5.7 (3).
\end{proof}

\bigskip

Department of Mathematics, Faculty of Science,

Yamaguchi University, Yamaguchi 753-8512, Japan

e-mail: andoy@po.cc.yamaguchi-u.ac.jp
\end{document}